\documentclass[preprint,3p,12pt]{elsarticle}

\usepackage{graphicx}
\usepackage{amsmath}
\usepackage{amsfonts}
\usepackage{enumitem} 							
\usepackage[labelfont=bf,tableposition=top]{caption} 	
\usepackage{subfig}									
\usepackage{epstopdf}								
\usepackage{floatrow}
\newfloatcommand{capbtabbox}{table}[][\FBwidth]
\newsavebox{\bigleftbox}
\floatstyle{plaintop}
\restylefloat{table}
\usepackage[table,xcdraw]{xcolor}
\usepackage{soul}
\usepackage{mathrsfs}  
\usepackage{multirow}
\usepackage{amsmath}
\usepackage{amssymb}
\usepackage[capitalise]{cleveref}
\crefname{rmk}{Remark}{Remarks}
\usepackage{titlesec}
\usepackage{etoolbox}
\usepackage{lipsum}
\usepackage{chngpage}
\usepackage{wrapfig}
\usepackage{subfiles}

\usepackage{tikz}

\newcommand{\numvec}[1]{\underline{\mathrm{#1}}}
\newcommand{\numvecT}[1]{\underline{\mathrm{#1}}^{\mathrm{T}}}
\def\doubleunderline#1{\underline{\underline{#1}}}
\newcommand{\mat}[1]{\doubleunderline{\mathrm{#1}}}

\renewcommand{\vec}[1]{\boldsymbol{#1}}

\newcommand{\scalar}{\phi}
\newcommand{\scalartest}{v}

\newcommand{\scalarhspace}{V^h}
\newcommand{\Ndofs}{{N_{\text{dofs}}}}

\newcommand{\Dom}{\Omega}
\newcommand{\TDom}{(0,T)}
\newcommand{\DomEl}{K}

\newcommand{\bdy}[1]{\partial #1}
\newcommand{\bdyDom}{{\bdy{\Dom}}}
\newcommand{\bdyNeum}{{\bdy{\Dom}_{N}}}
\newcommand{\bdyDirch}{{\bdy{\Dom}_{D}}}
\newcommand{\IF}{\Gamma}
\newcommand{\normal}{\vec{n}}

\makeatletter
\def\H{\@ifnextchar\bgroup\H@arg\H@noarg}
\newcommand{\H@arg}[1]{H^1( #1 )}
\newcommand{\H@noarg}{H^1}
\def\Ho{\@ifnextchar\bgroup\Ho@arg\Ho@noarg}
\newcommand{\Ho@arg}[1]{H^1_0( #1 )}
\newcommand{\Ho@noarg}{H^1_0}
\def\BH{\@ifnextchar\bgroup\BH@arg\BH@noarg}
\newcommand{\BH@arg}[1]{\vec{H}^1( #1 )}
\newcommand{\BH@noarg}{\vec{H}^1}
\def\BHo{\@ifnextchar\bgroup\BHo@arg\BHo@noarg}
\newcommand{\BHo@arg}[1]{\vec{H}^1_{\vec{0}}( #1 )}
\newcommand{\BHo@noarg}{\vec{H}^1_{\vec{0}}}
\def\Hdiv{\@ifnextchar\bgroup\Hdiv@arg\Hdiv@noarg}
\newcommand{\Hdiv@arg}[1]{\vec{H}( \text{div}; #1 )}
\newcommand{\Hdiv@noarg}{\vec{H}(\text{div})}
\def\Hdivo{\@ifnextchar\bgroup\Hdivo@arg\Hdivo@noarg}
\newcommand{\Hdivo@arg}[1]{\vec{H}_{\vec{0}}( \text{div}; #1 )}
\newcommand{\Hdivo@noarg}{\vec{H}_{\vec{0}}(\text{div})}
\def\L{\@ifnextchar\bgroup\L@arg\L@noarg}
\newcommand{\L@arg}[1]{L^2( #1 )}
\newcommand{\L@noarg}{L^2}
\def\Lp{\@ifnextchar\bgroup\Lp@arg\Lp@noarg}
\newcommand{\Lp@arg}[1]{L^p( #1 )}
\newcommand{\Lp@noarg}{L^p}
\def\Cz{\@ifnextchar\bgroup\Cz@arg\Cz@noarg}
\newcommand{\Cz@arg}[1]{\mathcal{C}^0( #1 )}
\newcommand{\Cz@noarg}{\mathcal{C}^0}
\def\Cinf{\@ifnextchar\bgroup\Cinf@arg\Cinf@noarg}
\newcommand{\Cinf@arg}[1]{\mathcal{C}^\infty( #1 )}
\newcommand{\Cinf@noarg}{\mathcal{C}^\infty}
\def\Cinfo{\@ifnextchar\bgroup\Cinfo@arg\Cinfo@noarg}
\newcommand{\Cinfo@arg}[1]{\mathcal{C}^\infty_0( #1 )}
\newcommand{\Cinfo@noarg}{\mathcal{C}^\infty_0}
\def\Cneg{\@ifnextchar\bgroup\Cneg@arg\Cneg@noarg}
\newcommand{\Cneg@arg}[1]{\mathcal{C}^{-1}( #1 )}
\newcommand{\Cneg@noarg}{\mathcal{C}^{-1}}
\def\CI{\@ifnextchar\bgroup\CI@arg\CI@noarg}
\newcommand{\CI@arg}[1]{\mathcal{C}^{1}( #1 )}
\newcommand{\CI@noarg}{\mathcal{C}^{1}}
\def\CII{\@ifnextchar\bgroup\CII@arg\CII@noarg}
\newcommand{\CII@arg}[1]{\mathcal{C}^{2}( #1 )}
\newcommand{\CII@noarg}{\mathcal{C}^{2}}
\makeatother

\makeatletter
\def\Proj{\@ifnextchar\bgroup\Proj@arg\Proj@noarg}
\newcommand{\Proj@arg}[1]{\mathscr{P}_{#1}}
\newcommand{\Proj@noarg}{\mathscr{P}}
\makeatother

\newcommand{\norm}[2][]{|\!| #2 |\!|_{#1}\,\!} 
\newcommand{\jump}[1]{\mbox{$[\![ #1 ]\!]$}}

\renewcommand{\d}[1]{\,\mathrm{d}#1}
\newcommand{\dDom}{\d{\Dom}}

\newcommand{\dbdy}{\d{S}}

\newcommand{\ddd}[2][]{\frac{\mathrm{d}^2 #1}{\mathrm{d} #2 ^2}}
\newcommand{\pdd}[2][]{\frac{\partial #1 }{\partial #2}}
\newcommand{\pddd}[2][]{\frac{\partial^2 #1}{\partial #2 ^2}}

\renewcommand{\Re}{\mathbb{R}}

\newcommand{\eigval}{\lambda}
\newcommand{\eigvalm}{\lambda_{\text{max}}}
\newcommand{\eigvec}{\numvec{\hat{\xi}}}
\newcommand{\eigvecm}{\numvec{\hat{\xi}}_{\text{max}}}
\newcommand{\eigfunc}{\xi}
\newcommand{\hfrac}{\chi}
\newcommand{\bbeta}{\bar{\beta}}
\newcommand{\gammaK}{ \gamma_K }
\newcommand{\bgammaK}{\bar{\gamma}_K}
\newcommand{\gammaM}{ \gamma_M }
\newcommand{\bgammaM}{\bar{\gamma}_M}

\newdefinition{rmk}{Remark}
\newcommand{\redcell}{\cellcolor[HTML]{E94C1F}}
\newcommand{\orangecell}{\cellcolor[HTML]{EE8866}}
\newcommand{\yellowcell}{\cellcolor[HTML]{EEDD88}}
\newcommand{\greencell}{\cellcolor[HTML]{5AAE61}}

\begin{document}

\begin{frontmatter}

\title{\large Critical time-step size analysis and mass scaling by ghost-penalty \\ for immersogeometric~explicit~dynamics}

\author[address1]{Stein K.F. Stoter\corref{cor1}}
\ead{K.F.S.Stoter@tue.nl}
\author[address1]{Sai C. Divi}
\ead{S.C.Divi@tue.nl}
\author[address1]{E.~Harald~van~Brummelen}
\ead{E.H.v.Brummelen@tue.nl}
\author[address2]{Mats~G.~Larson}
\ead{Mats.Larson@umu.se}
\author[address3]{Frits~de~Prenter}
\ead{F.d.Prenter@tudelft.nl}
\author[address1]{Clemens~V.~Verhoosel}
\ead{C.V.Verhoosel@tue.nl}

\cortext[cor1]{Corresponding author}

\address[address1]{Department of Mechanical Engineering, Eindhoven University of Technology, The Netherlands}
\address[address2]{Department of Mathematics and Mathematical Statistics, Ume\aa\ University, 901 87 Ume\aa, Sweden}
\address[address3]{Department of Aerospace Engineering, Delft University of Technology, The Netherlands}

\begin{abstract}
In this article, we study the effect of small-cut elements on the critical time-step size in an immersogeometric context. We analyze different formulations for second-order (membrane) and fourth-order (shell-type) equations, and derive scaling relations between the critical time-step size and the cut-element size for various types of cuts. In particular, we focus on different approaches for the weak imposition of Dirichlet conditions: by penalty enforcement and with Nitsche's method. The stability requirement for Nitsche's method necessitates either a cut-size dependent penalty parameter, or an additional ghost-penalty stabilization term is necessary. Our findings show that both techniques suffer from cut-size dependent critical time-step sizes, but the addition of a ghost-penalty term to the mass matrix serves to mitigate this issue. We confirm that this form of `mass-scaling' does not adversely affect error and convergence characteristics for a transient membrane example, and has the potential to increase the critical time-step size by orders of magnitude. Finally, for a prototypical simulation of a Kirchhoff-Love shell, our stabilized Nitsche formulation reduces the solution error by well over an order of magnitude compared to a penalty formulation at equal time-step size.
\end{abstract}


\end{frontmatter}

\newpage

\tableofcontents

\newpage

\section{Introduction}
\label{sec:Intro}


Explicit analysis forms the backbone of impact and crash-test simulation software. In these highly non-linear, short time-span simulations, the structure that is impacted is often comprised of shell-type components~\cite{LS-DYNA-theory,Abaqus66,AnsysExpDyn}. \textit{Isogeometric analysis} has emerged as a framework that streamlines the design-to-analysis pipeline for these types of simulations \cite{Cottrell2009,Leidinger2019,LS-DYNA-keywords}. In isogeometric analysis, the spline-based geometry representation from CAD drawings is used directly in the analysis software \cite{Cottrell2009}. The higher-order continuity of the spline basis functions makes them particularly well-suitable for handling the higher-order partial differential equations that arise in shell-type formulations \cite{LS-DYNA-keywords,Kiendl2009,Kiendl2015}. Additionally, the increased order of continuity of the basis functions leads to a higher permissible time-step magnitude compared to standard ($\mathcal{C}^0$-continuous) finite element basis functions for a fixed number of degrees of freedom \cite{Hughes2014,Puzyrev2018}. 

\begin{figure}[!b]
\vspace{-0.5cm}
    \centering
    \subfloat[\centering{Domain, basis functions and cut element. Green represents the physical domain interior, and red the exterior.}]{\includegraphics[width=0.57\textwidth]{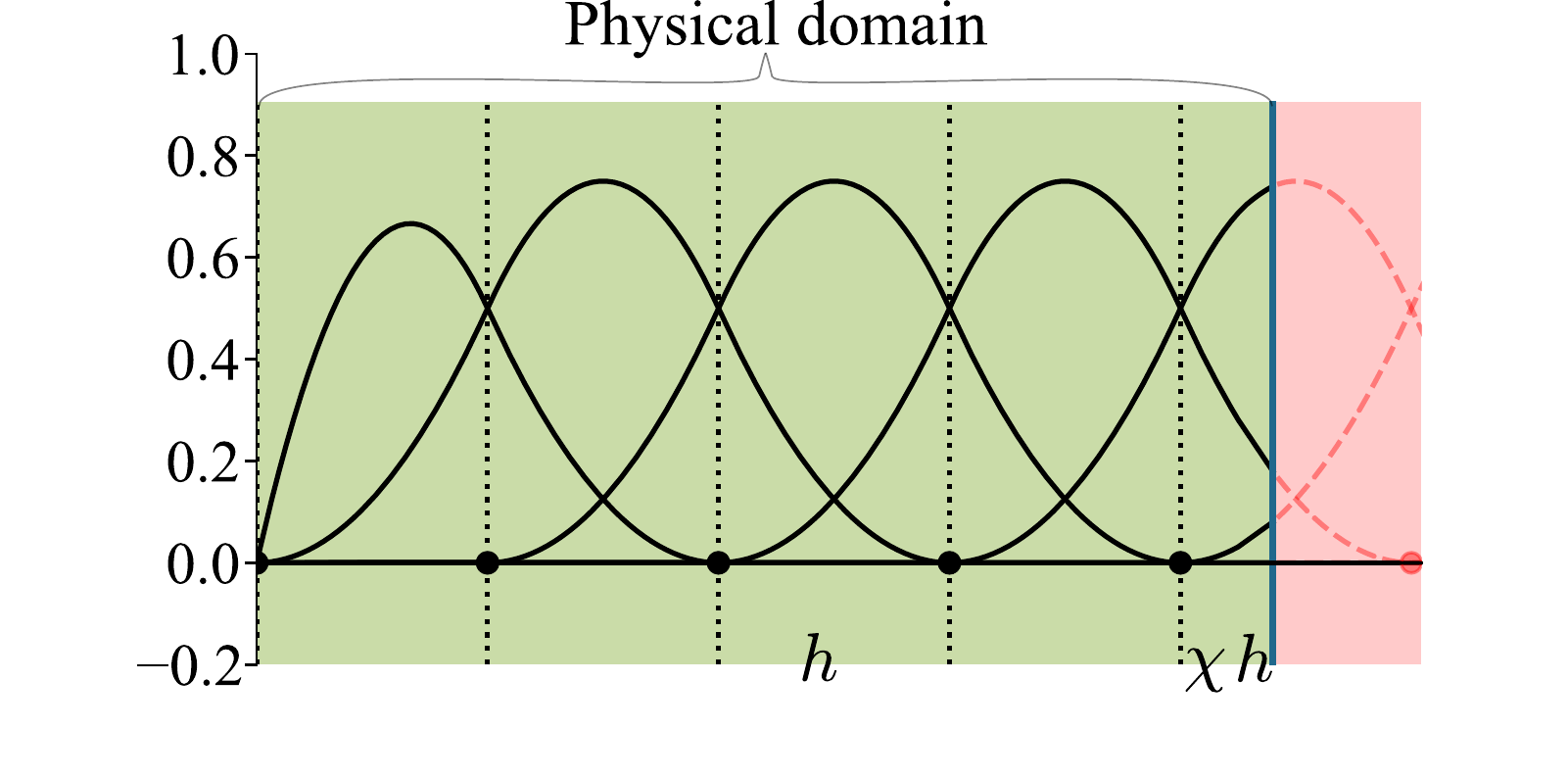} \label{fig:intro_a}}\hfill
    \subfloat[\centering{Spectra for various element size fractions, with and without mass lumping.}]{\includegraphics[width=0.4\textwidth]{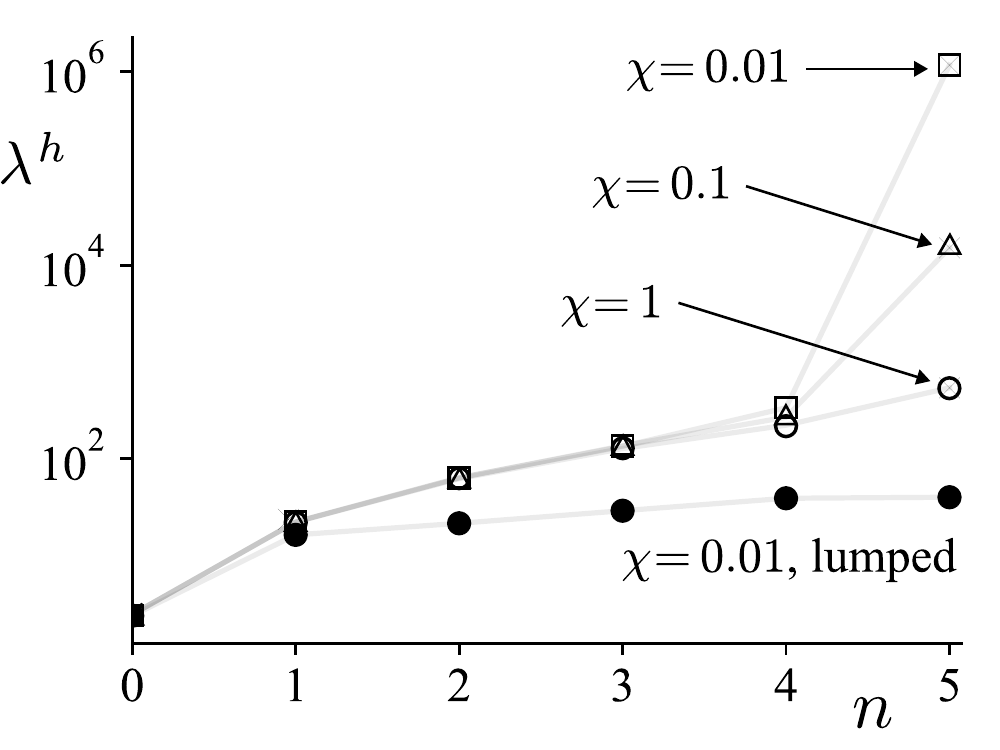} \label{fig:intro_b}}
    \caption{Effect of the small-cut elements (characterized by the element size fraction $\hfrac$ that is illustrated in \cref{fig:intro_a}) on the eigenvalues of a rod discretized with the six basis functions depicted in \cref{fig:intro_a}, for consistent and row-sum lumped mass matrices.}
    \label{fig:freqspec_cut}
\end{figure}

The spline geometry representation from CAD drawings often involves trimmed edges, where the physical domain boundary cuts through the computational background mesh. These ``immersed boundary'' methods, referred to as immersogeometric methods in the context of isogeometric analysis, have been extensively studied in recent years \cite{Hansbo2002a,Parvizian2007,Schillinger2011,Kamensky2015,Schillinger2015,Jonsson2017}. An example of a one-dimensional immersogeometric discretization is shown in \cref{fig:intro_a}, where the B-spline basis functions are not confined to the physical domain. The size of the cut is characterized by the elements size fraction $\hfrac$. The presence of basis functions with small support within the physical domain, i.e., with small $\hfrac$, negatively impacts the stability characteristics of the numerical scheme. In implicit analysis, these issues manifest themselves as poor conditioning of the resulting system of equations~\cite{Prenter2017}. For explicit analysis, a primary concern is the effect of small cut elements on the maximum eigenvalue, $\lambda^h_{\text{max}}$, of the mass-to-stiffness generalized eigenvalue problem. Indeed, for a one-dimensional rod example with the basis functions from \cref{fig:intro_a}, the cut-size heavily affects the maximum eigenfrequency, as can be observed from the diverging open markers in \cref{fig:intro_b}. For explicit time-stepping methods, an increase in the maximum eigenfrequency leads to a tighter bound on the magnitude of the permissible time step. The specific relation between $\lambda^h_{\text{max}}$ and this so-called ``critical time-step size'', $\Delta t_{\text{crit}}$, depends on the explicit time-integration technique that is employed. A typical example for undamped structural dynamics is
\begin{align}
    \Delta t_{\text{crit}}  = \frac{2}{\sqrt{\lambda^h_{\text{max}}}} \,, 
    \label{Dtcrit}
\end{align}
which holds for the Newmark-type central difference method~\cite[p.~493]{Hughes2000Book}. 

The diverging results of the open markers in \Cref{fig:intro_b} are obtained for a ``full'' or ``consistent'' mass matrix (in the sense of \cite{Archer1963}) in the formulation of the eigenvalue problem. In explicit analysis, the practice of mass diagonalization through some form of mass lumping is the gold standard \cite{Hughes2000Book}, as this is essential to attaining a fast solver. As can be observed from the filled marker in \cref{fig:intro_b}, the use of such a (row-sum) lumped mass matrix completely mitigates the detrimental impact of small cuts. This was first observed by Leidinger et al. \cite{LeidingerPhD}, where this behavior is explored analytically in one dimension, and numerically in two dimensions. 

The non-boundary fitted nature of immersed methods also requires specialized techniques for enforcing essential boundary conditions. Often, penalty-based methods are employed. In the context of explicit dynamics, not only the relation between the penalty parameter and the solution quality needs to be understood, but also its impact on the maximum eigenvalue. This is explored in the immersogeometric framework in \cite{Leidinger2019, LeidingerPhD}, where the authors report that sufficient solution accuracy can be attained with a penalty parameter that is small enough to avoid affecting the largest eigenvalue. When the penalty parameter exceeds a certain limit, the maximum eigenvalue scales linearly with the penalty parameter. The same linear scaling is observed for Nitsche's method in \cite{Harari2018}.

The current work aims to study the impact of small cuts on the critical time-step size in a wider range of scenarios, by investigating second- and fourth-order partial differential equations, regardless of the spatial dimension and polynomial order of the discretization, and considering various penalty- and Nitsche-based formulations for enforcing Dirichlet conditions. We also explore the use of ghost-penalty-based stabilization of the stiffness~\cite{Burman:15.2,Burman2010} and/or mass matrix \cite{Stoter2022c,Nguyen2022a}. Our goal is to identify the appropriate formulations in the explicit immersogeometric setting, and to evaluate the potential benefits of the alternative formulations in typical explicit dynamics computations.


The remainder of this article is structured as follows. In \cref{sec:Immersed}, we present the theory and nomenclature relevant to the topic of immersogeometric methods. 
We then focus on a linear wave equation in \cref{sec:SecondOrder}, where we systematically collect the contributions to the mass and stiffness matrices for various immersed boundary formulations. Once the general formulation is derived, we analytically derive minimal scaling-orders of the maximum eigenvalue with respect to the cut-element size. This analysis is repeated for a linear fourth-order equation in \cref{sec:FourthOrder}. 
The results obtained are experimentally verified in \cref{sec:NumExp}, where we also perform convergence studies for the linear wave equation, and transient simulations of a linear Kirchhoff-Love shell model. Finally, in \cref{sec:Conclusion}, we present our concluding remarks.

\section{Immersogeometric methods}
\label{sec:Immersed}

\begin{figure}[!b]
    \centering
    \subfloat[Ambient and physical domain.]{\includegraphics[width=0.33\linewidth]{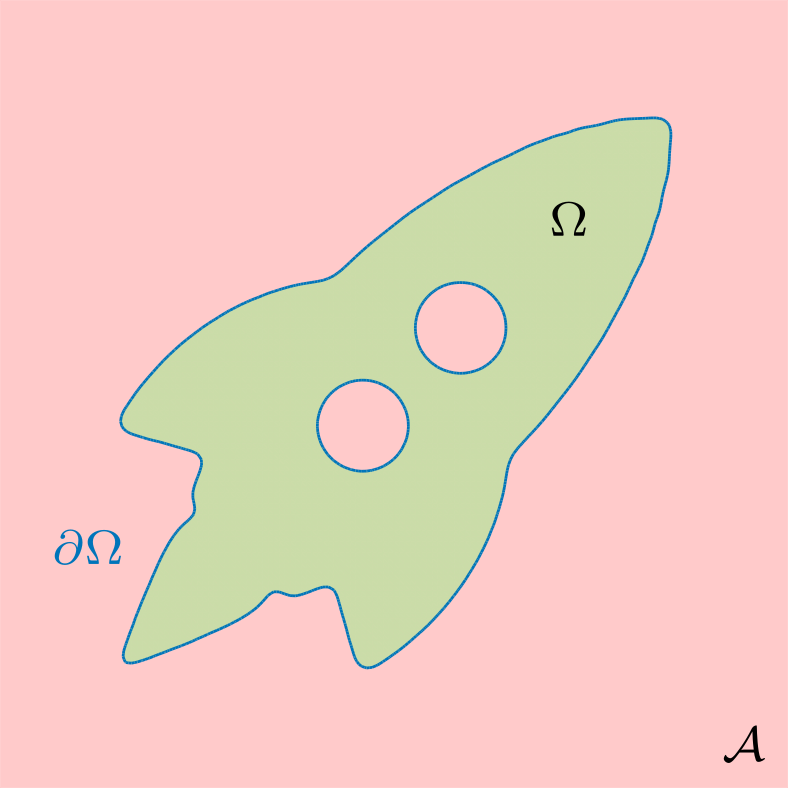}\label{fig:fcmdomaindef}}\hfill
    \subfloat[Ambient and background mesh.]{\includegraphics[width=0.33\linewidth]{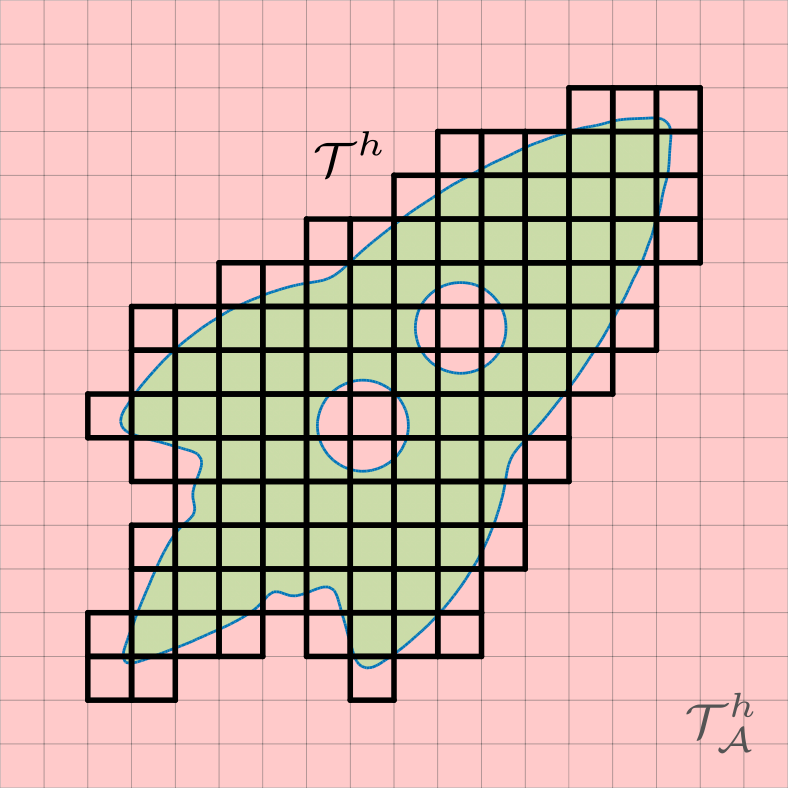}\label{fig:fcmmeshdef}}\hfill
    \subfloat[Ghost faces.]{\includegraphics[width=0.33\linewidth]{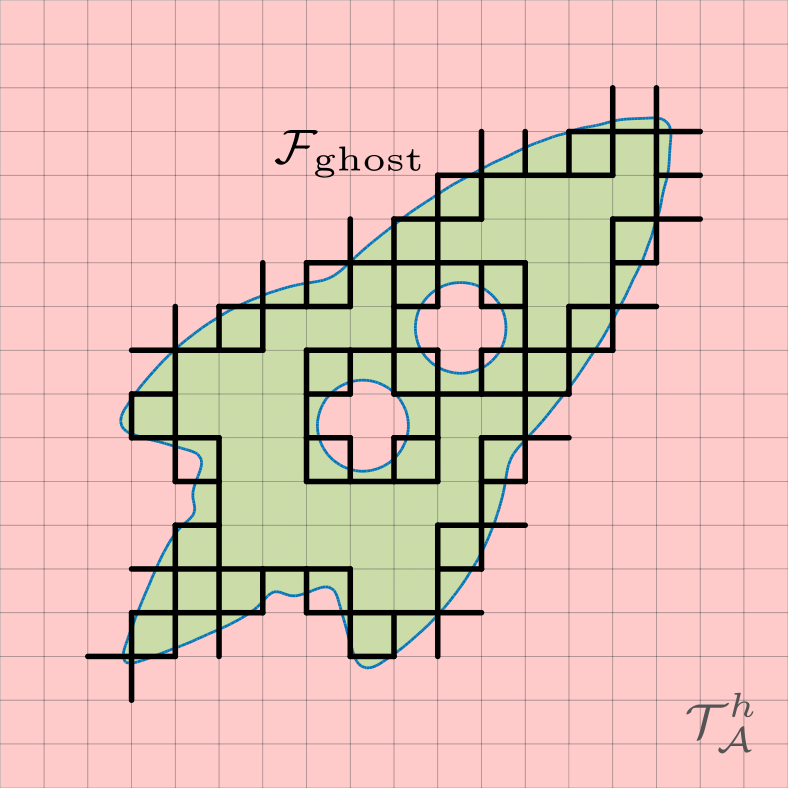}\label{fig:fcmghostfacets}}
    \caption{Definitions of the different meshes and domains in the immersed setting.}
	\label{fig:fcmdomain}
\end{figure}

Consider a physical domain $\Omega \subset \mathbb{R}^{d}$ ($d \in \{2,3\}$) and its boundary $\partial \Omega$, immersed in an ambient domain $\mathcal{A} \supset \Omega$, as shown in \cref{fig:fcmdomaindef}. The ambient mesh $\mathcal{T}_{\mathcal{A}}^{h}$ is a structured mesh covering the ambient domain $\mathcal{A}$. The ``background mesh'' is then defined as the collection of elements from the ambient mesh that intersect the physical domain:
\begin{equation}
\mathcal{T}^{h} := \big\{ K \in  \mathcal{T}^{h}_{\mathcal{A}} : K \cap \Omega \neq \emptyset\big\} \,.
\end{equation}
Both the ambient mesh and the background mesh are illustrated in \cref{fig:fcmmeshdef}. 

We also define a set of ``ghost faces'' as the collection of faces of elements from the background mesh that are intersected by the domain boundary:
\begin{equation} \label{eq:ghostfacets}
\mathcal{F}_{\rm ghost} = \big\{ \partial K \cap \partial K' \, : \, K \in \mathcal{G}, K'\in \mathcal{T}^{h},  K \neq K' \big\} \,,
\end{equation}
where $\mathcal{G} := \{ K \in \mathcal{T}^{h} \mid K \cap \partial \Omega \neq \emptyset \}$. The ghost faces are illustrated in Fig.~\ref{fig:fcmghostfacets}.

\subsection{Non-boundary fitted B-spline basis functions}
\label{ssec:Bsplines}

\begin{figure}[!b]
    \centering
    \subfloat[$p=1$]{\includegraphics[height=2.2cm,width=0.33\linewidth]{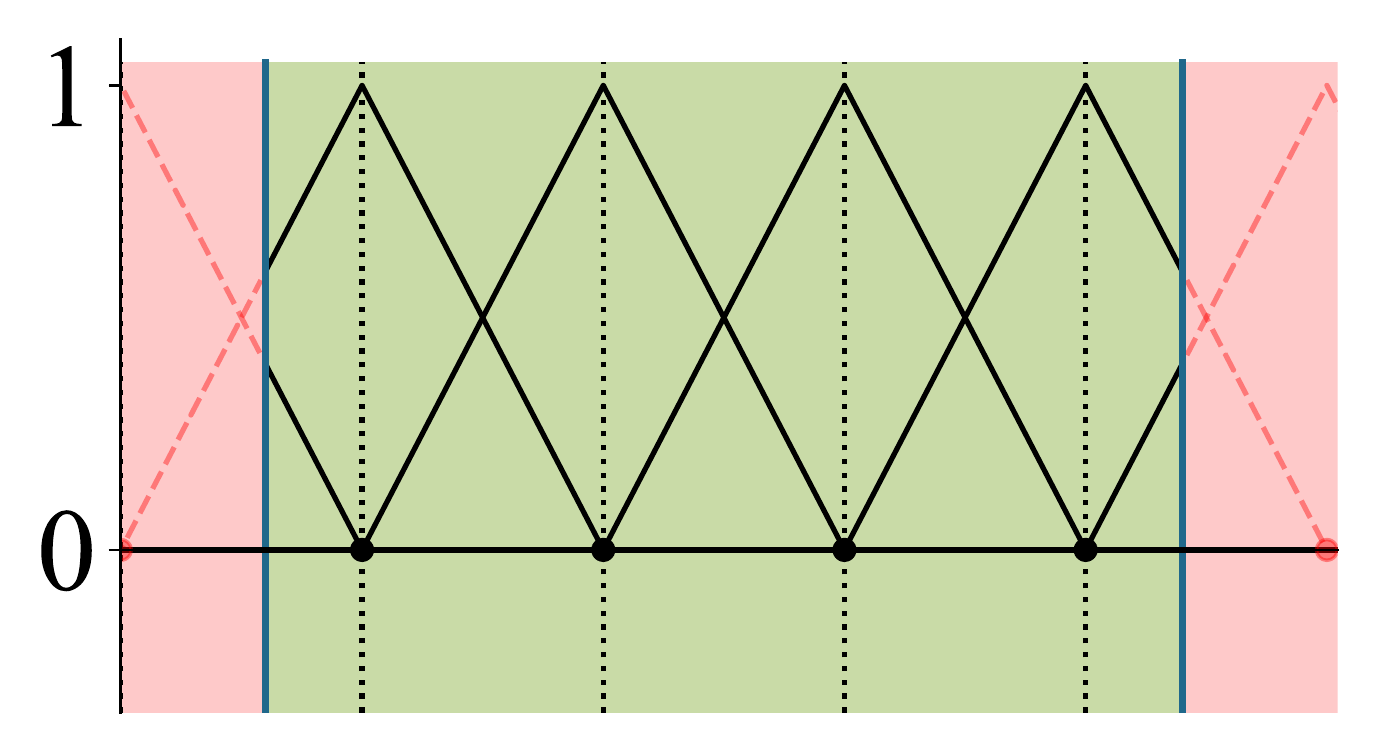}\label{fig:bsplinebasis_p1}}\hfill
    \subfloat[$p=1$]{\includegraphics[height=2.2cm,width=0.33\linewidth]{Figures/basis_spline_P1_k0.pdf}\label{fig:bsplinebasis_c0p1}}\hfill
    \subfloat[$p=1$]{\includegraphics[height=2.2cm,width=0.33\linewidth]{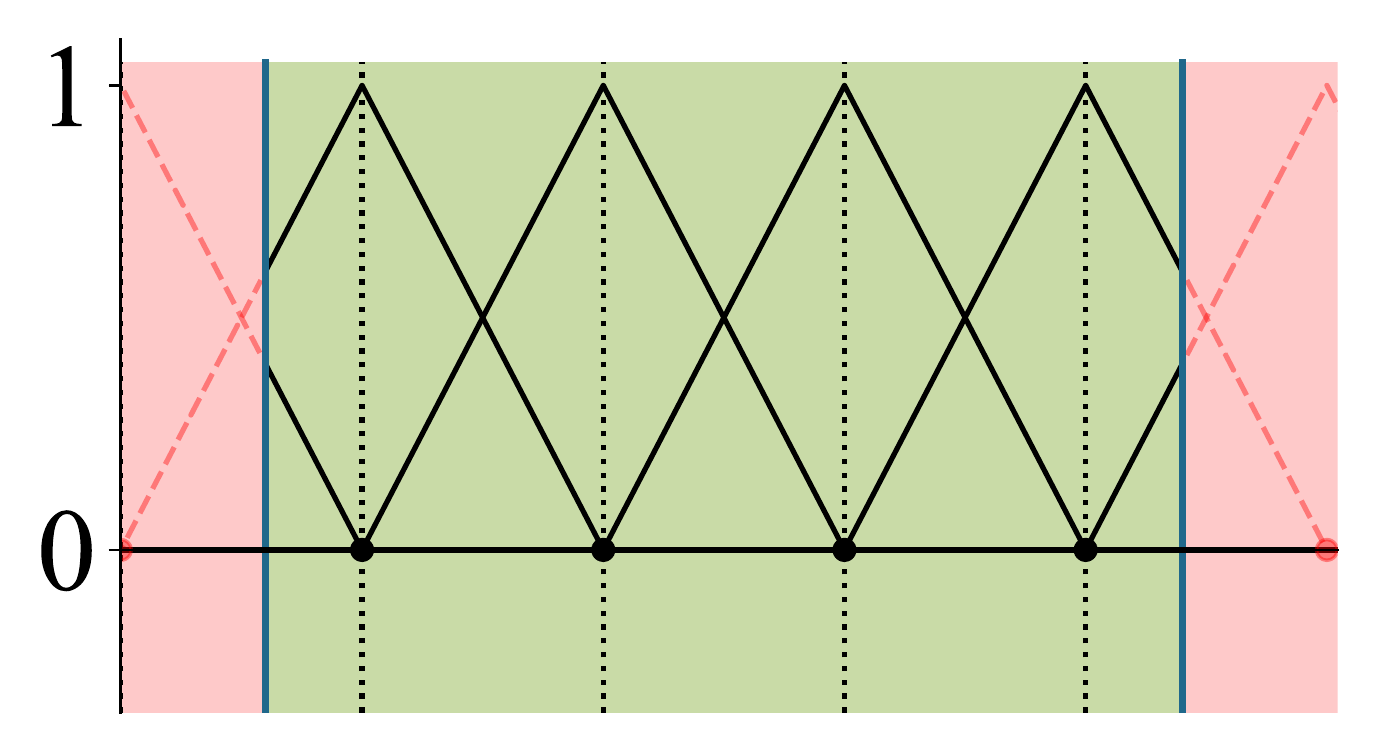}\label{fig:lagrangebasis_p1}}\\
    \subfloat[$p=2$]{\includegraphics[height=2.2cm,width=0.33\linewidth]{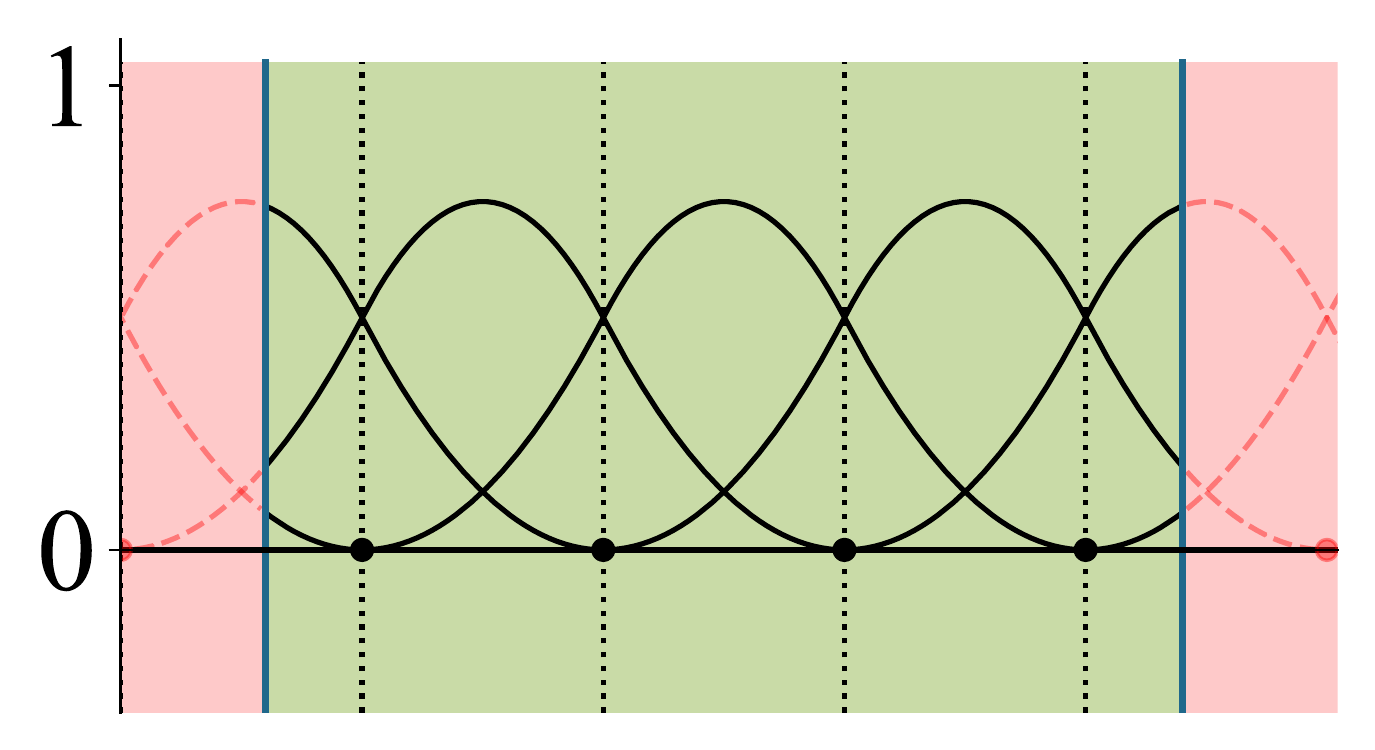}\label{fig:bsplinebasis_p2}}\hfill
    \subfloat[$p=2$]{\includegraphics[height=2.2cm,width=0.33\linewidth]{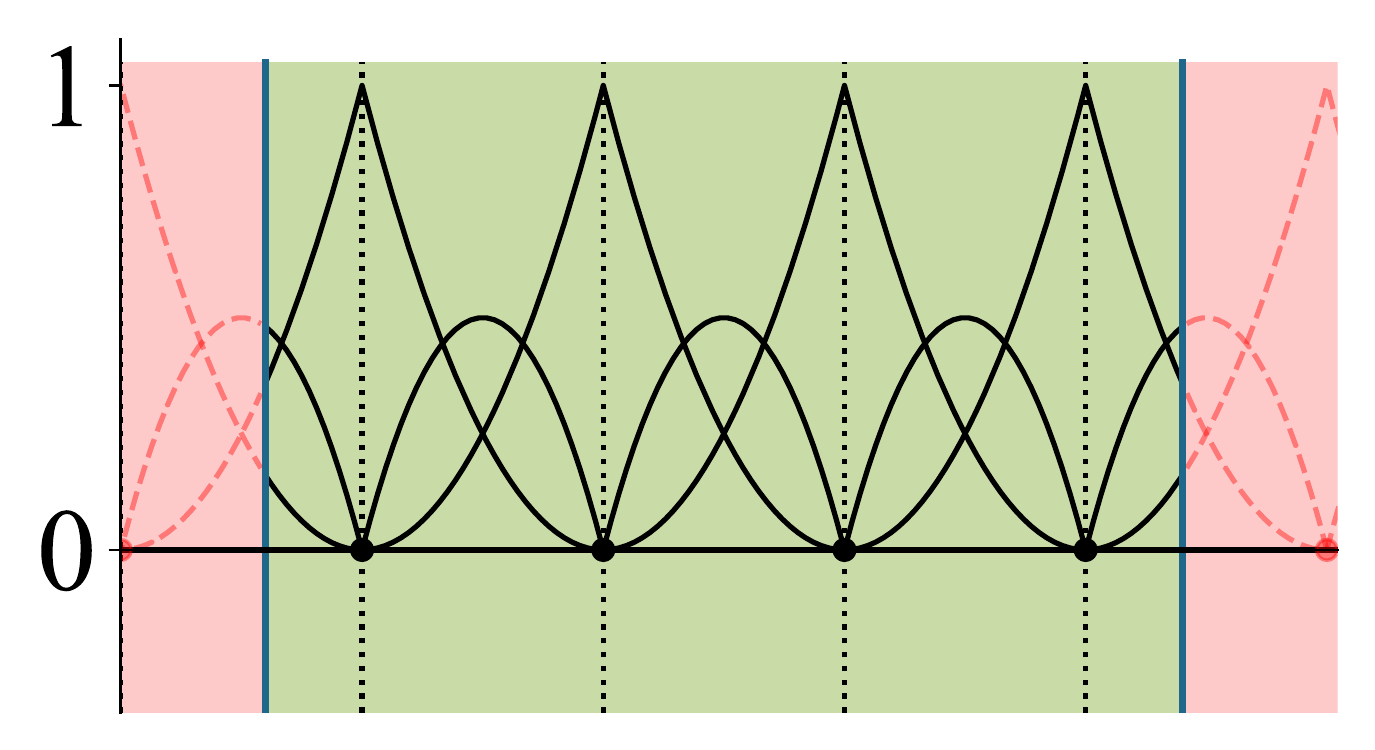}\label{fig:lagrangebasis_p2}}\hfill
    \subfloat[$p=2$]{\includegraphics[height=2.2cm,width=0.33\linewidth]{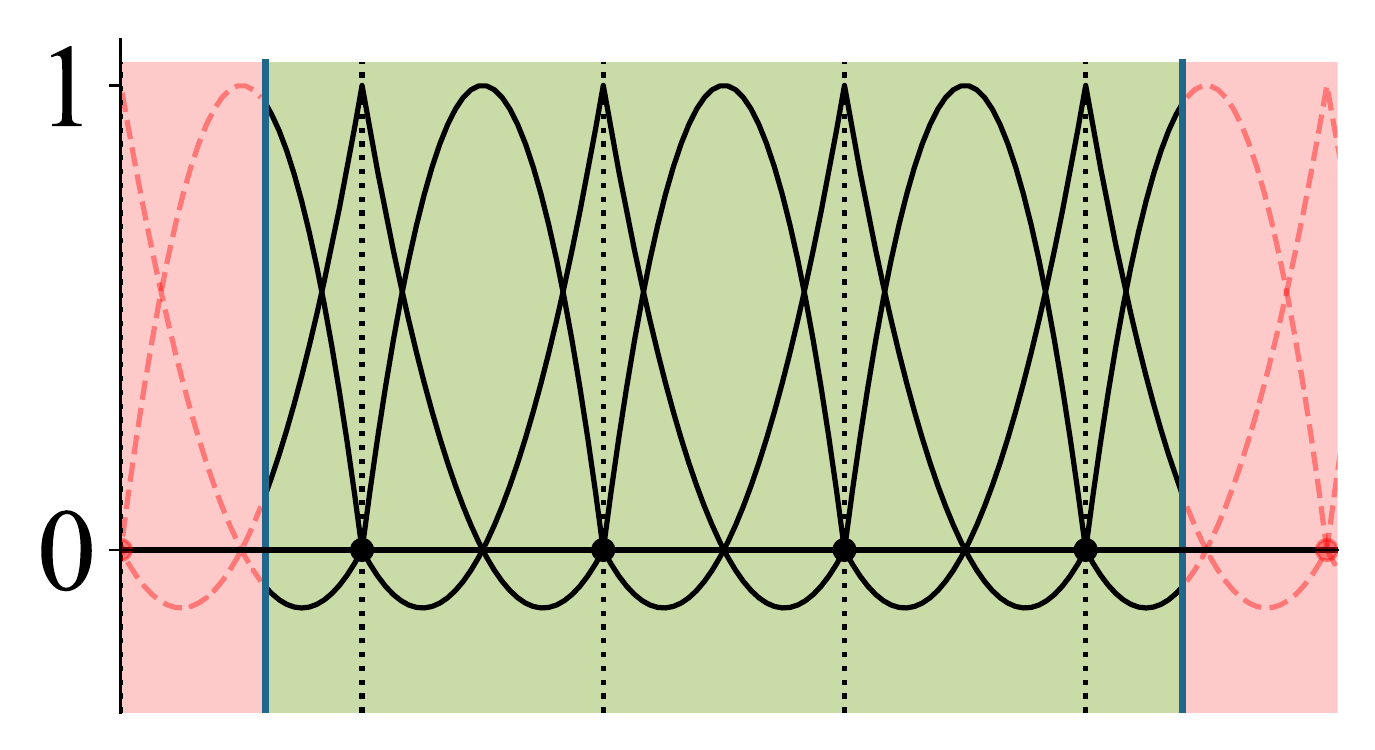}\label{fig:bsplinebasis_c0p2}}\\
    \subfloat[$p=3$]{\includegraphics[height=2.2cm,width=0.33\linewidth]{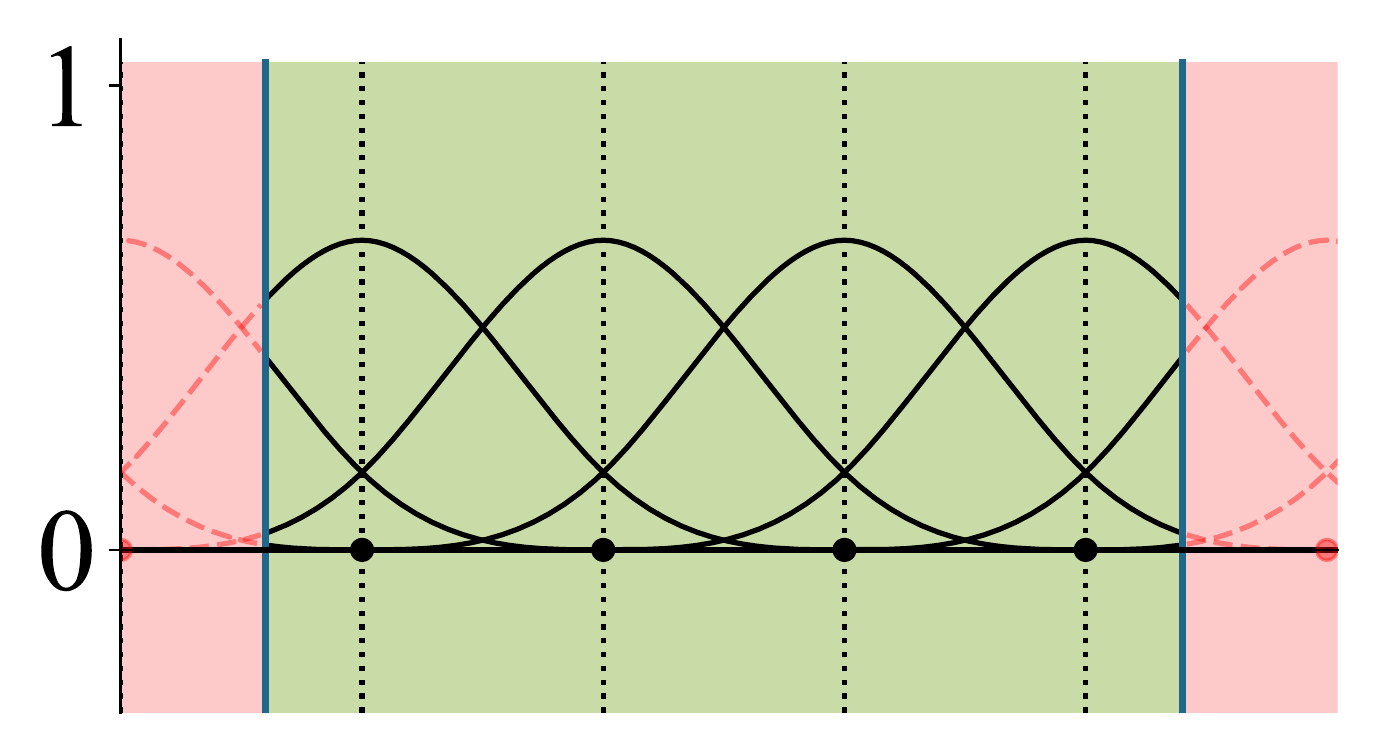}\label{fig:bsplinebasis_p3}}\hfill
    \subfloat[$p=3$]{\includegraphics[height=2.2cm,width=0.33\linewidth]{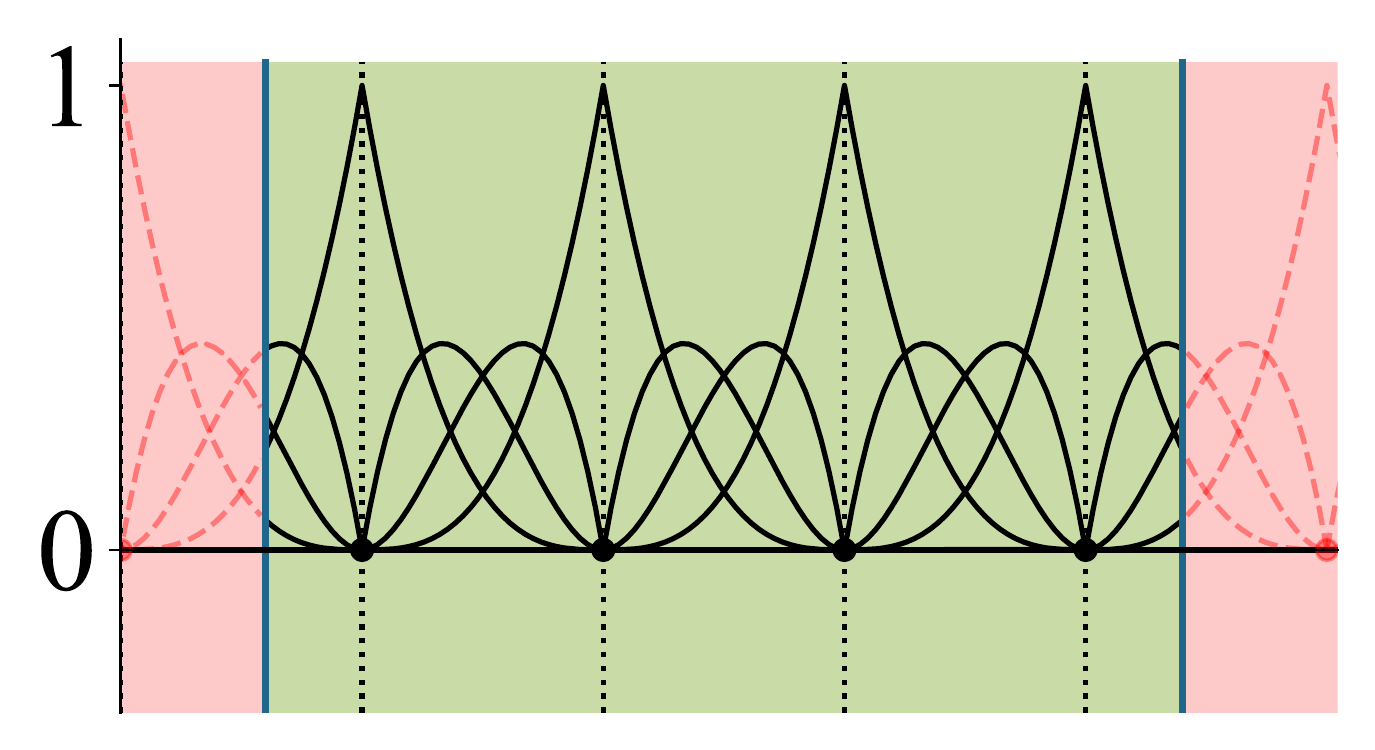}\label{fig:bsplinebasis_c0p3}}\hfill
    \subfloat[$p=3$]{\includegraphics[height=2.2cm,width=0.33\linewidth]{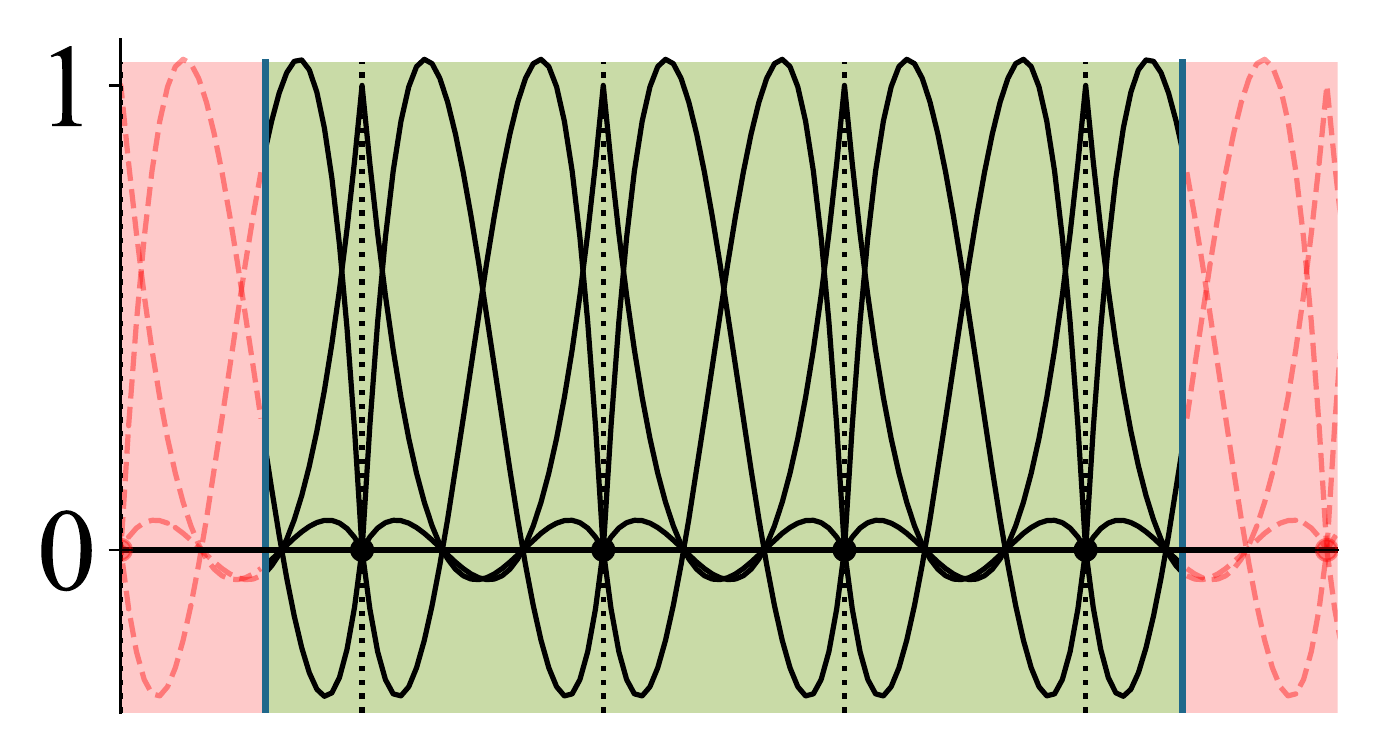}\label{fig:lagrangebasis_p3}}
    \caption{Immersed $\mathcal{C}^{p-1}$-continuous B-splines (left), $\mathcal{C}^{0}$-continuous B-splines (middle) and standard $\mathcal{C}^{0}$-continuous Lagrange basis functions (right).}
	\label{fig:basisfunc}
\end{figure}

The structured nature of the ambient mesh allows for the straightforward construction of multi-variate B-spline basis functions as tensor products between uni-variate B-spline functions. The left column in \cref{fig:basisfunc} depicts various orders of uni-variate B-spline functions of maximum regularity, while the middle column shows $C^0$-continuous uni-variate B-spline functions. An important property of these B-spline basis functions is their non-negativity, meaning that the integral value of each of these functions on any arbitrary subdomain is naturally non-negative as well. The same is not true for the conventional $\mathcal{C}^0$-continuous Lagrange basis functions, as shown in the right column of \cref{fig:basisfunc}.  This non-negativity property is particularly important in the context of explicit dynamics, as it guarantees positivity of the diagonal components of the row-sum lumped mass matrix, which will be addressed in \cref{ssec:rowsummass}.

In this article, we focus on maximum regularity B-splines and address lower-order regularity cases only in remarks. For maximum regularity B-splines, each multi-variate basis function $B_i(\vec{x})$ ($i=1,\cdots,J$) on the background mesh is a member of $\mathcal{C}^{p-1}(\mathcal{A})$, where $p$ is the polynomial order of the B-spline. The precise construction of the tensor product B-spline functions can be found, for example, in \cite{Cottrell2009}. To obtain the approximation space only on the physical domain, we first identify a subset of the maximum regularity B-splines that have non-zero support in the interior of the domain:
\begin{align}
    S = \big\{ N \in \{B_i(\vec{x})\}_{i=1}^J  : \text{supp}_\Omega(N) \neq 0\big\} \,.
\end{align}
We then number the remaining functions in this set from $i=1$ to $\Ndofs$, and define the $\Ndofs$-dimensional approximation space $V^h$ as their span:
\begin{equation}
\scalarhspace = \text{span}\left( S \right) = \text{span} \big\{ N_i(\vec{x}) \big\}_i^\Ndofs . \label{Vh}
\end{equation}
Under the assumption that the approximation functions evolve smoothly from time $t=0$ to time $t=T$, the time-dependent functions belong to the following semi-discrete function space:
\begin{equation}
\scalarhspace_T = \scalarhspace \otimes \mathcal{C}^\infty(0,T) \,. \label{VTh}
\end{equation}

\subsection{Integration of cut elements}
\label{ssec:integration}

For immersed finite element methods, integration procedures are crucial for capturing the physical geometry in the discrete formulation. The octree subdivision integration strategy, described in \cite{duster2008}, is a widely used approach due to its simplicity and robustness. However, the resulting large number of integration points may cause a significant computational cost increase during operation. A myriad of techniques to enhance octree-subdivision has been developed \cite{Verhoosel2022}. A few prominent techniques include error-estimate-based adaptive octree-subdivision \cite{Divi2020}, moment-fitting \cite{joulaian2016}, equivalent polynomial methods \cite{abedian2019}, 
 and the merged sub-cell technique \cite{peto2022}. 

\begin{figure}[!b]
    \centering
    \includegraphics[width=0.63\textwidth]{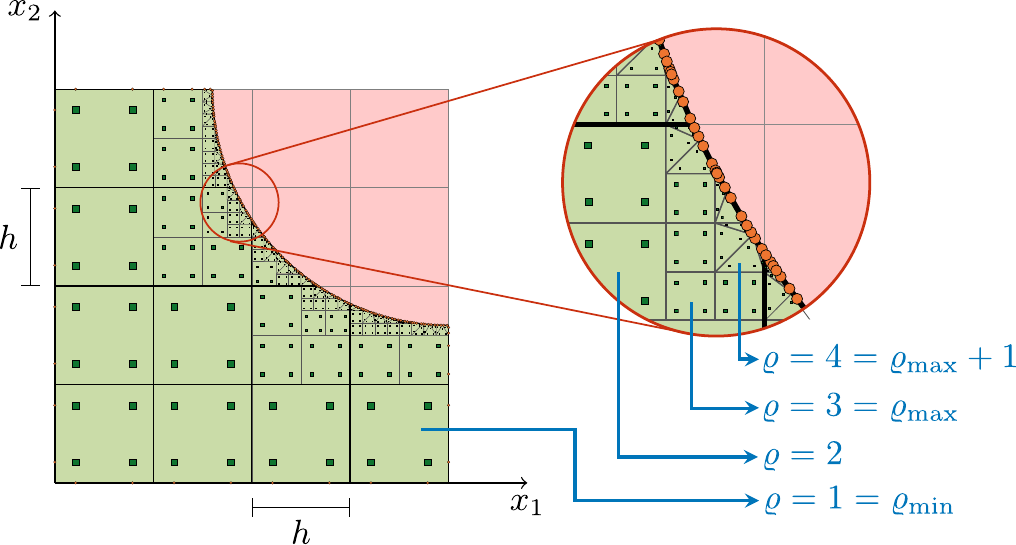}
    \caption{Volumetric (squares) and surface (circles) quadrature rules obtained by the octree integration procedure with tessellation at the lowest bisectioning level.}
    \label{fig:octreeintegration}
\end{figure}

In this article, we make use of the octree subdivision algorithm augmented with a tessellation step \cite{Verhoosel2015}. This approach consists of the following steps, as illustrated in \cref{fig:octreeintegration}: elements in the background mesh that intersect the boundary of the computational domain are bisected into $2^d$ sub-cells. If a sub-cell lies entirely within the domain, it is retained in the partitioning of the cut-cell, whereas it is discarded if it lies entirely outside the domain. This bisectioning procedure is recursively applied to the sub-cells that intersect the boundary, until $\varrho_{\rm max}$-times bisected sub-cells are obtained. At the lowest bisectioning level, a boundary tessellation procedure is applied to construct a $\mathcal{O}( h^2 / 2^{2 \varrho_{\rm max}} )$ accurate parametrization of the interior volume \cite{Verhoosel2015}. This tessellation procedure also provides a parametrization for the trimmed surface. The employed tessellation procedure is further detailed in \cite{Verhoosel2022}.

The octree procedure results in integration rules on each of the sub-cells and each of the surface triangles. Cut-element volumetric integration is then performed by agglomerating all sub-cell quadrature points, and cut-element surface integration by collecting all triangulated surface quadrature points. The accuracy of the cut-element integration scheme can be controlled through the selection of the octree depth $\varrho_{\rm max}$ and the quadrature rules on the sub-cells. 
 A two-dimensional illustration of the cut-element integration scheme is shown in \cref{fig:octreeintegration}, with the volumetric integration points represented by dark-green squares and the surface integration points depicted as orange circles.

\section{Critical time-step analysis for a second-order problem}
\label{sec:SecondOrder}

We center the exposition of our critical-timestep study for second-order initial/boundary-value problems around the linear wave equation:
\begin{subequations}%
    \begin{alignat}{2}
    \rho\pddd{t}\phi-\nabla\cdot(\kappa\nabla \scalar) &= f \qquad && \text{in }\Dom\times\TDom \,,\\
    -\kappa \nabla \scalar\cdot \normal &= g \qquad && \text{on }\bdyNeum\times\TDom\,, \\[0.1cm]
    \scalar &= \scalar_D \qquad && \text{on }\bdyDirch\times\TDom\,, \label{ibvp_Dirch}\\[0.1cm]
    \scalar &= \scalar_0 \qquad && \text{on }\Dom\times\{ 0 \}\,,\\[-0.1cm]
    \pdd{t}\scalar &= \dot{\scalar}_{0} \qquad && \text{on }\Dom\times\{ 0 \}\,,
    \end{alignat}\label{ibvp}%
\end{subequations}
where $\phi$ is the unknown field, $\rho$ and $\kappa$ are parameters of the propagating medium, $f$, $g$, $\phi_D$, $\scalar_{0}$ and $\dot{\scalar}_0$ are the prescribed body force, Neumann data, Dirichlet data and initial state, respectively, $\Dom\subset \Re^d$ and $T>0$ are the spatial domain and final time, and $\bdyNeum$ and $\bdyDirch=\bdyDom\setminus\bdyNeum$ are the Neumann and Dirichlet segments of the domain boundary.

This equation describes, for example, the out-of-plane displacement of a pre-stressed vibrating string or membrane in one, respectively two, dimensions, and the propagation of pressure waves in one, two or three dimensions. A weak formulation of this initial-boundary value problem reads \cite{Lions1973}:
\begin{align}
    &\text{For a.e. }t\in\TDom\text{, find } \scalar \in \H_{\scalar_D}(\Dom) \text{ and } \pddd{t}\scalar \in H^{-1}(\Omega) \text{ s.t. } \forall\,  \scalartest \in \Ho{\Dom}: \nonumber\\
    &\quad \begin{cases} 
    \displaystyle \big\langle \pddd{t}\scalar ,  \rho\,\scalartest \big\rangle + \int\limits_\Dom  \kappa \nabla \scalar \cdot \nabla \scalartest \dDom = \int\limits_\Dom f\,\scalartest \dDom - \int\limits_{\bdyNeum} \! g \, \scalartest \dbdy \,, \\
    \displaystyle \phi\big|_{t=0} = \phi_0 \, , \\[0.2cm]
    \displaystyle \pdd{t}\phi\big|_{t=0} = \dot{\scalar}_{0} \, ,
    \end{cases}\label{weakform}
\end{align}
where $\H_{\scalar_D}(\Dom) = \{ \phi \in H^1(\Omega): \phi\big|_{\bdyDirch} = \phi_D  \}$, $H^{-1}(\Omega)$ is the corresponding dual space, and $\big\langle\cdot,\cdot\big\rangle$ denotes the pairing between them. If the second time derivative of $\phi$ is a member of $L^2(\Omega)$, then $\big\langle \pddd{t}\scalar , \rho\, \scalartest \big\rangle = \int\limits_\Dom  \rho\pddd{t}\scalar \,  \scalartest \dDom$.

\subsection{Semi-discrete formulation}

With the aim of constructing a stable explicit finite element approximation for the weak formulation of \cref{weakform} on a non-boundary-fitted isogeometric mesh, we consider the \hbox{(semi-)discrete} spaces $\scalarhspace$ and $\scalarhspace_T$, from \cref{Vh,VTh}, as our test and trial spaces, respectively. $\scalarhspace$ is defined as the span of basis functions, such that any test function in $\scalarhspace$, for example $\scalartest^h \in \scalarhspace$, can be represented as a linear combination of the basis functions:
\begin{align}
    \scalartest^h(\vec{x}) = \sum\limits_{i=1}^{\Ndofs} \hat{\text{\scalartest}}_i N_i(\vec{x}) = \numvecT{\hat{\scalartest}} \numvec{N}(\vec{x}) \label{vhrep}\,.
\end{align}
Similarly, any trial function in $\scalarhspace_T$, for example $\scalar^h \in \scalarhspace_T$, can be written as a linear combination of the basis functions with time-dependent weighting coefficients:
\begin{align}
    \scalar^h (t,\vec{x}) = \sum\limits_{i=1}^{\Ndofs} \hat{\scalar}_i(t) \, N_i(\vec{x}) = \numvecT{\hat{\scalar}}\!(t) \, \numvec{N}(\vec{x}) \quad \text{with }\hat{\scalar}_i(t) \in\mathcal{C}^\infty (0,T)\,.\label{phihrep}
\end{align}

In the following subsections, we provide a comprehensive overview of the individual components that comprise the mass and stiffness matrices of the semi-discrete formulation.

\subsubsection{Neumann formulation}
First, we consider a pure Neumann problem, i.e., $\bdyDom = \bdyNeum$ and $\bdyDirch=\emptyset$. Then, for any $t\in\TDom$, a finite element approximation produces the following semi-discrete formulation:
\begin{align}
&\text{Find } \scalar^h \in \scalarhspace_T  \text{ s.t. }\forall\,  \scalartest^h \in \scalarhspace : \nonumber\\
&\quad \int\limits_\Dom \rho\ddd{t}\scalar^h \scalartest^h \dDom + \int\limits_\Dom  \kappa \nabla \scalar^h \cdot \nabla \scalartest^h \dDom = \int\limits_\Dom f\,\scalartest^h \dDom - \int\limits_{\bdyNeum} \! g \, \scalartest^h \dbdy \,. \label{neumannprob}
\end{align}

Substitution of the representation of $\scalar^h$ and $\scalartest^h$ from \cref{phihrep,vhrep} into \cref{neumannprob} results in the following system of ordinary differential equations
\begin{align}
& \mat{M} \, \ddd{t}\numvec{\hat{\scalar}} + \mat{K} \, \numvec{\hat{\scalar}}  =  \numvec{F} \label{semidiscrete} \,,
\end{align}
with:
\begin{subequations}\label{MKF}
\begin{alignat}{3}
    \big[\, \mat{M} \,\big]_{ij} & = M\big(N_i(\vec{x}), N_j(\vec{x})\big) = \int\limits_{\Dom} \rho N_i (\vec{x})  N_j (\vec{x}) \d{\Dom} \,, \label{consistentmass}\\
    \big[\, \mat{K} \,\big]_{ij} & = K\big(N_i(\vec{x}), N_j(\vec{x})\big) =  \int\limits_{\Dom} \kappa \nabla N_i (\vec{x}) \cdot \nabla N_j (\vec{x}) \d{\Dom} \,,\\
    \big[\, \mat{F} \,\big]_{i} & =  F\big(N_i(\vec{x})\big) = \int_{\Dom} f(\vec{x}) \, N_i (\vec{x}) \d{\Dom} - \int\limits_{\bdyNeum} g(\vec{x}) \,  N_i (\vec{x}) \dbdy \,.
\end{alignat}
\end{subequations}
For this Neumann problem, the immersogeometric framework primarily impacts the selection of the discrete approximation space $\scalarhspace$ as the span of non-boundary fitted B-spline basis functions, and the procedure for carrying out the integrals in \cref{MKF}.

\subsubsection{Row-sum mass lumping}
\label{ssec:rowsummass}
A fully-discrete formulation of \cref{semidiscrete} typically follows from a finite-difference type approximation to the time-derivative. When an explicit time-stepping scheme is adopted, such as a second-order Newmark-type central difference method or a higher-order explicit Runge-Kutta scheme, then the time-marching does not necessitate an inverse of the stiffness matrix. It does, however, require an inverse of the mass matrix. To avoid significant computational expense (both in terms of storage and operation count), the mass matrix is often manipulated to attain a diagonal matrix, which can be inverted trivially. This manipulation process is referred to as ``mass lumping'' and various methods exist, such as diagonal scaling \cite{Hinton1976}, manifold-based methods \cite{Yang2017} and lumping by nodal quadrature \cite{Strang1973,Patera1984}. In this article, we focus on row-sum lumping, where the diagonal value is set to the sum-total of the row. This sum-total of the row corresponds to the multiplication of the row by a vector of ones. For a partition of unity basis, such as the employed B-spline basis, the vector of ones represents a field of unit value. Consequently, the lumped mass matrix can be written as:
\begin{align}
    \big[\, \mat{M}_D \big]_{ij} = \begin{cases} 0 \qquad & \text{if } i\neq j \\
     M\big(1,N_i(\vec{x}) \big) = \int\limits_\Dom \rho N_i \dDom  \qquad & \text{if } i= j
    \end{cases} \label{lumpedmass}\,.
\end{align}

Ensuring that the diagonal entries of $\mat{M}_D$ are positive is of crucial importance: negative components would cause negative eigenvalues, inducing exponential growth of the corresponding eigenmode. As pointed out in \cref{ssec:Bsplines}, the non-negativity property of the B-spline basis-functions guarantees positivity of $\int_\Dom N_i(x) \dDom$, and thus of the diagonal entries of $\mat{M}_D$, irrespective of how elements are cut. The same cannot be guaranteed for, for example, more classical $\mathcal{C}^0$-continuous Lagrange basis functions.

\begin{rmk}\label{rmk:masslumpprob}
In the analysis in \cref{ssec:analysis2nd,ssec:analysis4th} we require evaluations of the vector-matrix-vector product $\numvecT{\hat{\scalartest}} \mat{M}_D \numvec{\hat{\scalartest}}$ for various $\numvec{\hat{\scalartest}}$. Unlike the consistent mass matrix (for which $ \numvecT{\hat{\scalartest}} \mat{M} \numvec{\hat{\scalartest}} = \int_\Omega \rho\, v^hv^h \dDom$), the ad-hoc nature of the mass-lumping procedure implies that the operation $\numvecT{\hat{\scalartest}} \mat{M}_D \numvec{\hat{\scalartest}}$ does not permit an integral-based bilinear form. However, for those particular functions whose vector of coefficients ($\numvec{\hat{\scalartest}}$ in \cref{vhrep}) exclusively consists of ones and zeros ($[\numvec{\hat{\scalartest}}]_i \in \{ 0,1\}$), the following integral evaluation of the vector-matrix-vector product is valid:
\begin{align}
    \numvecT{\hat{\scalartest}} \mat{M}_D \numvec{\hat{\scalartest}} =  M\big(1,\scalartest^h(\vec{x}) \big) = \int\limits_\Dom \rho v^h(\vec{x}) \dDom =: M_D\big(\scalartest^h(\vec{x})\big) \,.
\end{align}
\end{rmk}

\subsubsection{Dirichlet condition enforcement by penalization}

The remaining challenge in solving \cref{weakform} in an immersed setting, is the imposition of Dirichlet boundary conditions. In \cref{weakform}, these are essential conditions, and they are imposed on the function space itself. For boundary-fitted methods, this can be mimicked by strongly prescribing nodal values. For immersed methods, however, the nodes are no longer placed on the domain boundary and the conventional procedure is not feasible. However, we can still integrate along this immersed boundary, as described in \cref{ssec:integration}. Hence, boundary conditions can be incorporated weakly; through the addition of integral terms in the weak formulation targeted at enforcing the constraints. 

The most common approach for weak imposition of Dirichlet conditions in explicit analysis is by penalty enforcement \cite{Leidinger2019}:
\begin{align}
&\text{Find } \scalar^h \in \scalarhspace_T \text{ s.t. } \forall\,  \scalartest^h \in \scalarhspace: \nonumber\\
&\quad M(\ddd{t}\scalar^h,\scalartest^h) + K(\scalar^h,\scalartest^h) +  \underbrace{\int_\bdyDirch \kappa \beta\, \scalar^h \, \scalartest^h \dbdy}_{\displaystyle K_\beta \big(\phi^h, v^h\big)}  = F(\scalartest^h) + \underbrace{\int_\bdyDirch \kappa \beta\, \scalar_D\, \scalartest^h \dbdy}_{\displaystyle F_\beta \big( v^h \big)}\,, \label{penaltyform}
\end{align}
where $\beta>0$. The additional contributions to the stiffness matrix and force vector read:
\begin{subequations}\label{Kbeta}
\begin{alignat}{3}
    \big[\, \mat{K}_\beta \,\big]_{ij} & = K_\beta \big(N_i(\vec{x}), N_j(\vec{x})\big) = \int\limits_\bdyDirch \kappa \beta\, N_i(\vec{x}) N_j(\vec{x}) \dbdy\\
    \big[\, \mat{F} _\beta\,\big]_{i} & = F_\beta \big(N_i(\vec{x})\big) = \int\limits_\bdyDirch \kappa \beta\, \scalar_D(\vec{x}) N_i(\vec{x}) \dbdy \,.
\end{alignat}
\end{subequations}
To be dimensionally consistent, $\beta^{-1}$ needs to have unit length. Suitable scaling of the eigenvalues associated to the penalty term is achieved when the penalty is chosen to scale inversely with the size of the elements of the background mesh $h_\DomEl$ \cite{Leidinger2019,Apostolatos2015,Herrema2019,Pasch2021}:
\begin{align}
    \beta \big|_\DomEl = \bar{\beta} \, h_\DomEl^{-1} \,,
\end{align}
where $\bar{\beta}$ is a dimensionless global constant. 

The penalty method offers a number of advantages, such as the absence of stringent restrictions on the penalty parameter to ensure positivity of eigenvalues, owing to the positive semi-definiteness of the contribution to the stiffness matrix, and its ease of implementation. A significant drawback is its variationally inconsistent nature (in the sense of \cite{Strang1973}), which leads to a loss of optimal convergence rate that may result in error increases of orders of magnitude. The optimal convergence rate may be retrieved by choosing a penalty scaling stronger than $h^{-1}$ \cite{Barrett1987}. This is, however, not a viable option in the context of explicit analysis, as the larger penalty value would soon increase the largest eigenvalue, and hence negatively affects the critical time-step size. Moreover, the solution quality is sensitive to the choice of penalty parameter, for which rigorous estimates are not available. Also in this regard, one must exercise caution not to choose the penalty values too large to avoid impacting the critical time-step size \cite{Leidinger2019}.

\subsubsection{Dirichlet condition enforcement by Nitsche's method}
\label{ssec:Nitsche}

The penalty method can be augmented by terms to make it variationally consistent, mitigating many of the just mentioned drawbacks. When this is done in a symmetric manner, this is called Nitsche's method \cite{Nitsche1971}:
\begin{align}
&\text{Find } \scalar^h \in \scalarhspace_T \text{ s.t. } \forall\,  \scalartest^h \in \scalarhspace: \nonumber\\
\begin{split}
&\quad M(\ddd{t}\scalar^h,\scalartest^h) + K(\scalar^h,\scalartest^h)  +  K_\beta ( \scalar^h , \scalartest^h ) \underbrace{ - \int_\bdyDirch \!\! \kappa \nabla \scalar^h\cdot\normal\, \scalartest^h \dbdy - \int_\bdyDirch \!\! \kappa \nabla \scalartest^h\cdot\normal \,  \scalar^h \dbdy }_{\displaystyle K_{\text{cs}} \big( \phi^h , v^h \big)}   \\
&\quad  \hspace{4cm} = F(\scalartest^h) + F_\beta (\scalartest^h) \underbrace{ - \int_\bdyDirch \!\! \kappa \scalar_D \nabla  \,  \scalartest^h\cdot\normal\dbdy }_{\displaystyle F_{\text{s}} \big( v^h \big)} \,. 
\end{split} \raisetag{2cm}\label{nitscheform}
\end{align}
The new matrix and vector contributions originating from the consistency and symmetry terms may be identified as:
\begin{subequations}\label{Ksym}
\begin{alignat}{3}
    \big[\, \mat{K}_{\text{cs}} \,\big]_{ij} & = K_{\text{cs}} \big(N_i(\vec{x}), N_j(\vec{x})\big) = - \int\limits_\bdyDirch  \kappa \nabla N_i\cdot\normal\, N_j \dbdy  - \int\limits_\bdyDirch \kappa \nabla N_j\cdot\normal \,  N_i \dbdy  \,,\\
    \big[\, \mat{F}_{\text{s}} \,\big]_{i} & = F_{\text{s}} \big(N_i(\vec{x})\big) = - \int\limits_\bdyDirch \kappa \scalar_D \,  \nabla N_i\cdot\normal \dbdy  \,.
\end{alignat}
\end{subequations}

Due to the added terms, the stiffness matrix is no-longer unconditionally positive definite. When $\beta$ is chosen too small, the stiffness matrix may include negative eigenvalues. Negativity of eigenvalues carries over to the mass-to-stiffness generalized eigenvalue problem, again leading to detrimental exponential growth of the corresponding eigenmodes in time. The restriction on $\beta$ to ensure positive definiteness has been studied extensively \cite{Barbosa1991,Stenberg1995,Dolbow2009}, and follows in each element $\DomEl$ from a local inverse estimate. We make use of the following choice:
\begin{align}
    \beta \big|_\DomEl = 2 \sup\limits_{\scalartest^h\in\scalarhspace} \frac{ \norm{\nabla \scalartest^h \cdot \normal }^2_{\bdyDirch\cap \DomEl}  }{ \norm{ \nabla \scalartest^h }^2_{\Dom\cap \DomEl}  } \propto \frac{1}{h_c}\, \Big|_\DomEl \,, \label{betac}
\end{align}
where $h_c$ is a length-scale associated to the cut element. 

The element size for cut elements is not unambiguously defined. The inverse estimate in \cref{betac} can be bound by a size-independent factor multiplied by the ratio of the area of the Dirichlet surface to the volume of the cut-element \cite{Warburton2003,Prenter2018}. This length parameter provides a suitable size measure for small-cut elements, however it does not limit to $h_K$ for (nearly) uncut cases. In what follows, we make use of the following definition of $h_c$, which is consistent for both small and large cuts:
\begin{align}\label{hc}
    h_c \big|_\DomEl = \min\Big\{  \,\,\,  \frac{ \int_{\Dom \cap \DomEl} \d{\Dom} }{ \int_{ \bdyDirch\cap \DomEl } \dbdy } \,\,\,\,,\,\,\,\, \left[ \int_{\Dom \cap \DomEl } \d{\Dom} \right]^{\frac{1}{d}}  \,\,\,  \Big\} =: \hfrac  h_\DomEl \,.
\end{align}
The element size fraction $\hfrac$ that implicitly follows from this definition as $\hfrac = h_c / h_\DomEl$ represents the central quantity in our study on the cut-sensitivity of the largest eigenvalues later on. For sliver cuts, $\hfrac$ equates to the cut element volume fraction $\eta=\int_{\Dom \cap \DomEl }\! \d{\Dom} / \int_{\DomEl }\! \d{\Dom} $ (used in earlier analysis \cite{Prenter2017}), and for shape-regular cut elements they relate as $\hfrac\propto\eta^{\frac{1}{d}}$.

We define the constant of proportionality in \cref{betac} as $\bbeta\big|_\DomEl$, such that, by definition, we can make use of the following expression for $\beta\big|_\DomEl$:
\begin{align}
    \beta\big|_\DomEl = \bar{\beta}\big|_\DomEl \, (\hfrac \, h_\DomEl)^{-1} \,.
\end{align}

\subsubsection{Ghost-penalty stabilization}

For elements with vanishing support in the physical domain, the coefficient $\beta$ required to ensure positive definiteness of the complete stiffness matrix can become arbitrarily large, as is reflected by the requirement in \cref{betac}. For implicit and steady analysis, such a large penalty factor has negative consequences on the stability and conditioning of the ensuing system \cite{Prenter2018}. To mitigate this issue, the support of small-cut basis functions must be extended into the domain interior. Doing so strongly, for example by adopting weighted extended basis B-splines (WEB-splines) \cite{Hollig2003}, or by performing cell aggregation \cite{Badia2018}, requires manipulation of the basis functions. If this is undesirable, then a potential alternative is the addition of ``ghost-penalty'' stabilization \cite{Burman:15.2,Burman2010}:
\begin{align}
&\text{Find } \scalar^h \in \scalarhspace_T \text{ s.t. } \forall\,  \scalartest^h \in \scalarhspace: \nonumber\\
\begin{split}
&\,\, M(\ddd{t}\scalar^h,\scalartest^h) + K(\scalar^h,\scalartest^h)  +  K_\beta ( \scalar , \scalartest^h ) +  K_{\text{cs}} ( \scalar , \scalartest^h )  +  \underbrace{   \int_{\IF_g} \kappa \gammaK \, \jump{\partial^{k+1}_n \scalar^h } \jump{\partial^{k+1}_n \scalar^h }  \dbdy }_{\displaystyle K_{\gamma} \big( \phi^h , v^h \big)}   \\[-0.7cm]
&\quad  \hspace{3cm} = F(\scalartest^h) + F_\beta (\scalartest^h)  + F_{\text{cs}} (\scalartest^h)  \,,
\end{split} \label{ghostform}
\end{align}
where $\IF_g = \bigcup \mathcal{F}_{\rm ghost} $ is the union of the ghost faces from \cref{eq:ghostfacets} and \cref{fig:fcmghostfacets}, $\jump{\cdot}$ is the jump operator, and $\partial^{k+1}_n$ is the normal gradient of order $k+1$, with $k$ the order of continuity of the B-spline basis functions. In a more general sense, $K_{\gamma} \big( \cdot , \cdot \big)$ should include jump-terms of the normal derivatives of order 1 until the polynomial order $p$. Since we only consider maximum order of continuity splines, for which $k=p-1$, all but the highest normal derivatives vanish on element interfaces. The contribution to the stiffness matrix is:
\begin{align} \label{Kgamma}
    \big[\, \mat{K}_{\gamma} \,\big]_{ij} & = K_{\gamma} \big( N_i(\vec{x}), N_j(\vec{x}) \big) = \int\limits_{\IF_g} \kappa \gammaK \, \jump{\partial^{k+1}_n N_i } \jump{\partial^{k+1}_n N_j }  \dbdy  \,.
\end{align}

The new penalty between the domain interior elements and the cut elements effectively adds a stiffness to the deflection of weakly supported degrees of freedom. More technically, the ghost-penalty term extends coercivity from the domain interior to the background mesh \cite{Divi2022}. When the parameter $\gammaK$ is large enough and scales with $h_{\DomEl}^{2p-1}$, $\beta$ is permitted to scale with the background-element size $h_\DomEl$ rather than with the cut-element size $h_c = \hfrac\, h_\DomEl$:
\begin{subequations}
    \begin{alignat}{2}
    \gammaK\big|_\DomEl & = \bgammaK \,  h_{\DomEl}^{2p-1} \,,  \\[0.1cm]
    \beta\big|_\DomEl & = \bbeta \, h_{\DomEl}^{-1} \,.
    \end{alignat}
\end{subequations}
The minimal permitted values of the pre-factors $\bgammaK$ and $\bbeta$ is still an open research question~\cite{Divi2022,Badia2022}. In the current work, we choose a sufficiently large $\bgammaK$ to enable the use of a small $\bbeta$, and experimentally verify that the ensuing stiffness matrix is positive definite.

\subsubsection{Ghost mass}

The final ingredient that we choose to add to our explicit immersogeometric formulation is a ``ghost-mass'' term. We propose this term as a type of consistent mass scaling, with the intent of reducing the maximum eigenvalues that are caused by small-cut elements. The following additional component is introduced to the mass matrix \cite{Deng2021,Stoter2022c,Nguyen2022a}:
\begin{align} \label{inertialghost}
    &M_\gamma ( \ddd{t} \scalar^h , \scalartest^h ) =   \int\limits_{\IF_g} \rho \gammaM \, \jump{\partial^{k+1}_n \ddd{t}\scalar^h } \jump{\partial^{k+1}_n \scalar^h }  \dbdy\,, \\
    &\big[\,\mat{M}_\gamma \big]_{ij} = M_\gamma \big( N_i(\vec{x}), N_j(\vec{x}) \big) \,. \label{ghostmass}
\end{align}
If lower-order regularity B-splines are used, $M_\gamma (\cdot,\cdot)$ should include penalties on the jumps of all lower-order normal derivatives as well. These jumps vanish for the $\mathcal{C}^{p-1}$-continuous B-splines considered herein.
 
The ghost-mass term adds inertia to the acceleration of the deflection of the weakly supported degrees of freedom. As $\mat{M}_\gamma$ represents a positive semi-definite contribution to the mass matrix, it serves to reduce the eigenvalues of modes that excite the $M_\gamma ( \cdot , \cdot )$ term, i.e., those with derivative changes across boundaries of elements with small interior support. At the same time, the term is consistent for sufficiently smooth solutions, such that it does not introduce a modeling error.

The required scaling of the penalty parameter $\gammaM$ is different than that of $\gammaK$ in the stiffness matrix, as already follows from a dimensional consistency argument. We make use of the following scaling:
\begin{align}
    \gammaM = \bar{\gamma}_M \,  h_\DomEl^{2p+1} \,.
\end{align}
where the appropriate choice of $\bar{\gamma}_M$ is investigated in \cref{ssec:analysis2nd}.

\begin{rmk} \label{rmk:ghostmasslump}
The matrix that follows from the ghost-mass term is not amendable to standard row-sum lumping. Recall from \cref{ssec:rowsummass} that the row-sum represents the action of the corresponding bilinear form on a field of unit value, and, as can be observed from \cref{inertialghost}, $M_\gamma ( 1 , \scalartest^h ) = 0 $. When the immersed boundary elements make up a relatively small portion of the mesh, the computational expense involved in inverting the mass matrix may not cause a bottleneck (especially when the LU-factorization is stored and reused). Nevertheless, various strategies could be considered to efficiently approximate an inverse to the scaled mass matrix. One could reduce the number of ghost faces in such a way that $\mat{M}_\gamma$ becomes block-diagonal \cite{Larson2021}, or approximate the inverse in a (block) Jacobi sense, or attempt row-sum lumping of the absolute values of the matrix components, or consider different classes of discrete extension operators altogether \cite{Burman2022}. However, since the primary focus of this work is the analysis of the explicit immersogeometric formulations, we consider the development of such optimized implementation strategies beyond the scope of this article, and we exclusively consider $\mat{M}_\gamma$ in its full form.
\end{rmk}

\subsection{Analysis of the critical time-step size}
\label{ssec:analysis2nd}

The various stiffness terms proposed in the previous section can be collected in the total stiffness bilinear form $\tilde{K}(\cdot,\cdot)$, to produce the total stiffness matrix $\mat{\tilde{K}}$. Similarly, the total inertial bilinear form $\tilde{M}(\cdot,\cdot)$ and mass matrix $\mat{\tilde{M}}$ follow from the chosen combination of a consistent or lumped mass matrix, and potentially a ghost-mass contribution. The general semi-discrete form then reads:
\begin{align}
\begin{split}
&\text{Find } \scalar^h \in \scalarhspace_T \text{ s.t. } \forall\,  \scalartest^h \in \scalarhspace: \\[-0.1cm]
&\quad \tilde{M}(\ddd{t}\scalar^h,\scalartest^h) + \tilde{K}(\scalar^h,\scalartest^h) = \tilde{F}(\scalartest^h) \,,
\end{split} \qquad \Leftrightarrow \qquad \mat{\tilde{M}} \, \ddd{t}\numvec{\hat{\scalar}} + \mat{\tilde{K}} \, \numvec{\hat{\scalar}} = \numvec{\tilde{F}} \,. \label{generalizedform}
\end{align}
As addressed in the introduction, the critical time-step size of an explicit time-stepping treatment of \cref{generalizedform} is inversely related to the maximum eigenvalue of the generalized eigenvalue problem:
\begin{align}
    \tilde{K}(\eigfunc^h,\scalartest^h)  = \lambda \tilde{M}(\eigfunc^h,\scalartest^h) \quad\forall\,\scalartest^h\in\scalarhspace \qquad \Leftrightarrow \qquad  \mat{\tilde{K}} \, \eigvec = \eigval \mat{\tilde{M}}\, \eigvec \,.
\end{align}

For symmetric matrices, the maximum eigenvalue is the largest value of the generalized Rayleigh quotient:
\begin{align}
   \eigvalm \geq \mathcal{R}(\mat{\tilde{K}},\mat{\tilde{M}},\eigvec) = \frac{\eigvec^\text{T}\,\mat{\tilde{K}} \, \eigvec  }{\eigvec^\text{T}\,\mat{\tilde{M}} \, \eigvec  } = \frac{ \tilde{K}(\eigfunc^h,\eigfunc^h) }{\tilde{M}(\eigfunc^h,\eigfunc^h)  } \qquad \forall\, \eigvec\in\Re^N \,,\,\, \eigfunc^h =\sum\limits_{n=0}^N \eigvec_n N_n(\vec{x})\,, \label{eigvalRayleigh}
\end{align}
where the equality holds for $\eigvec = \eigvecm$. Due to the bilinearity of the forms, the various components separate into individual contributions:
\begin{align}
 \mathcal{R}(\mat{\tilde{K}},\mat{\tilde{M}},\eigvec) = \frac{ K(\eigfunc,\eigfunc) +  K_\text{cs}(\eigfunc,\eigfunc) +  K_\beta(\eigfunc,\eigfunc) + K_\gamma (\eigfunc,\eigfunc) }{M(\eigfunc,\eigfunc) + M_\gamma (\eigfunc,\eigfunc)  }\,.
\end{align}

We now wish to analyse whether $\eigvalm$, and by induction, the critical time-step size, is sensitive to the size and shape of the cut elements, as characterized by the parameter $\hfrac$ in \cref{hc}. For $\hfrac\rightarrow 0$, the generalized Rayleigh quotient approaches the ratio of the lowest-order scaling of each of the components in the numerator divided by the lowest-order scaling of each of the components in the denominator. For some $\eigvec$, we can then characterize the scaling of the generalized Rayleigh quotient as:
\begin{align}\label{scalingsconcl}
    \mathcal{R}(\mat{K},\mat{M},\eigvec) = \mathcal{O}\big(\hfrac^q) \,.
\end{align}
If there exists any $\eigvec$ such that $q<0$, then, by \cref{eigvalRayleigh}, the maximum eigenvalue can become arbitrarily large for arbitrarily small cuts, resulting in an unfeasibly small critical time-step size. Those $\eigvec$ for which $q=0$ cause cut-\textit{shape} dependent eigenvalues that may turn out to be dominant. Those corresponding to $q>0$ are suppressed as the cut becomes small.

\begin{figure}[!b]\vspace{-0.5cm}
    \centering
    \subfloat[First corner-cut function, \\\centering$\eigvec= \begin{bmatrix}\cdots 0,1,0,\cdots \end{bmatrix}^\mathrm{T}$.] {\includegraphics[trim=0 70 0 0, clip, width=0.42\linewidth]{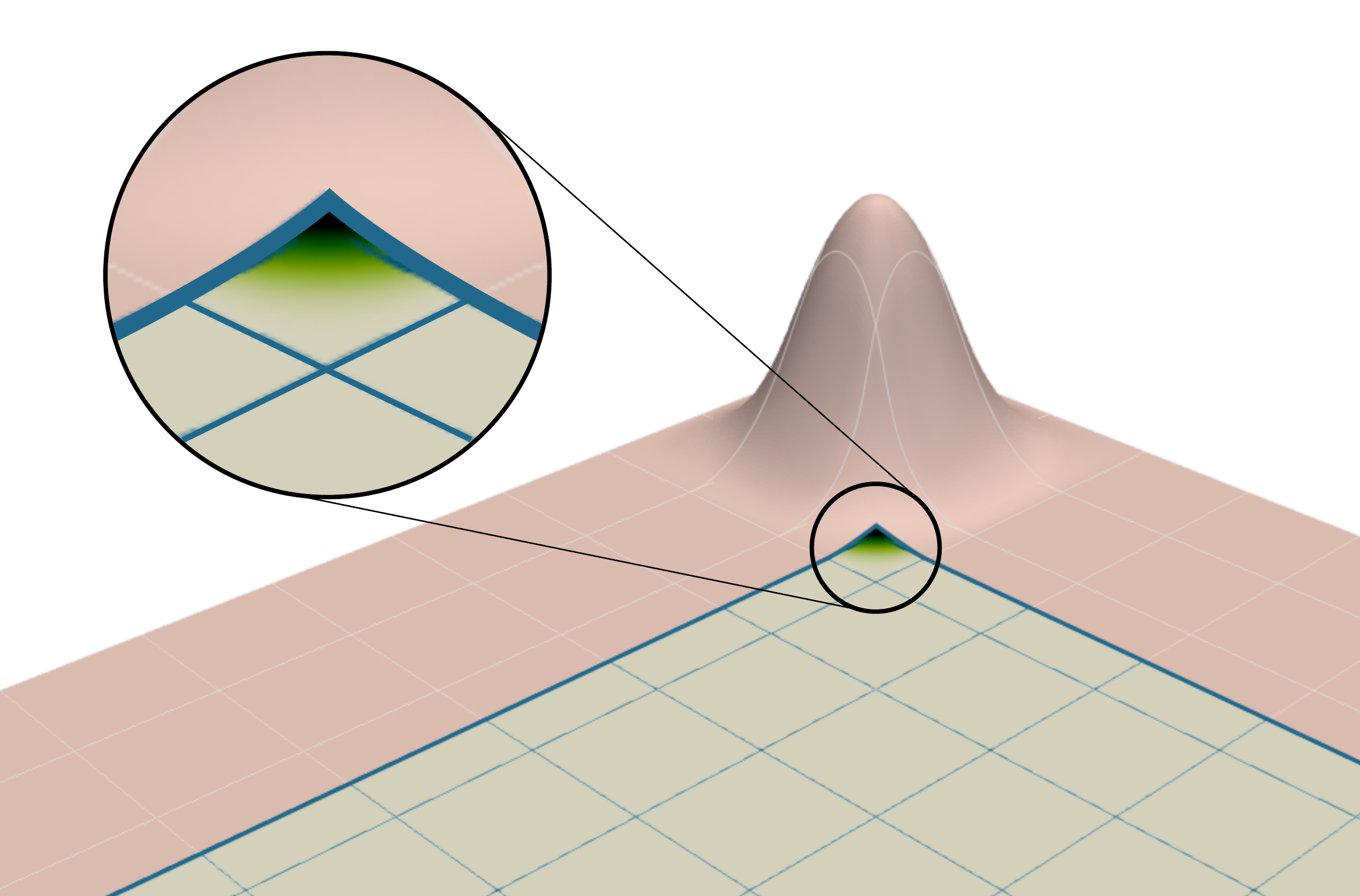}\label{fig:cut_funcs1}}\hspace{1.3cm}
    \subfloat[Second corner-cut function, \\\centering$\eigvec= \begin{bmatrix}\cdots0,1,0,\cdots \end{bmatrix}^\mathrm{T}$.]{\includegraphics[trim=0 70 0 0, clip, width=0.42\linewidth]{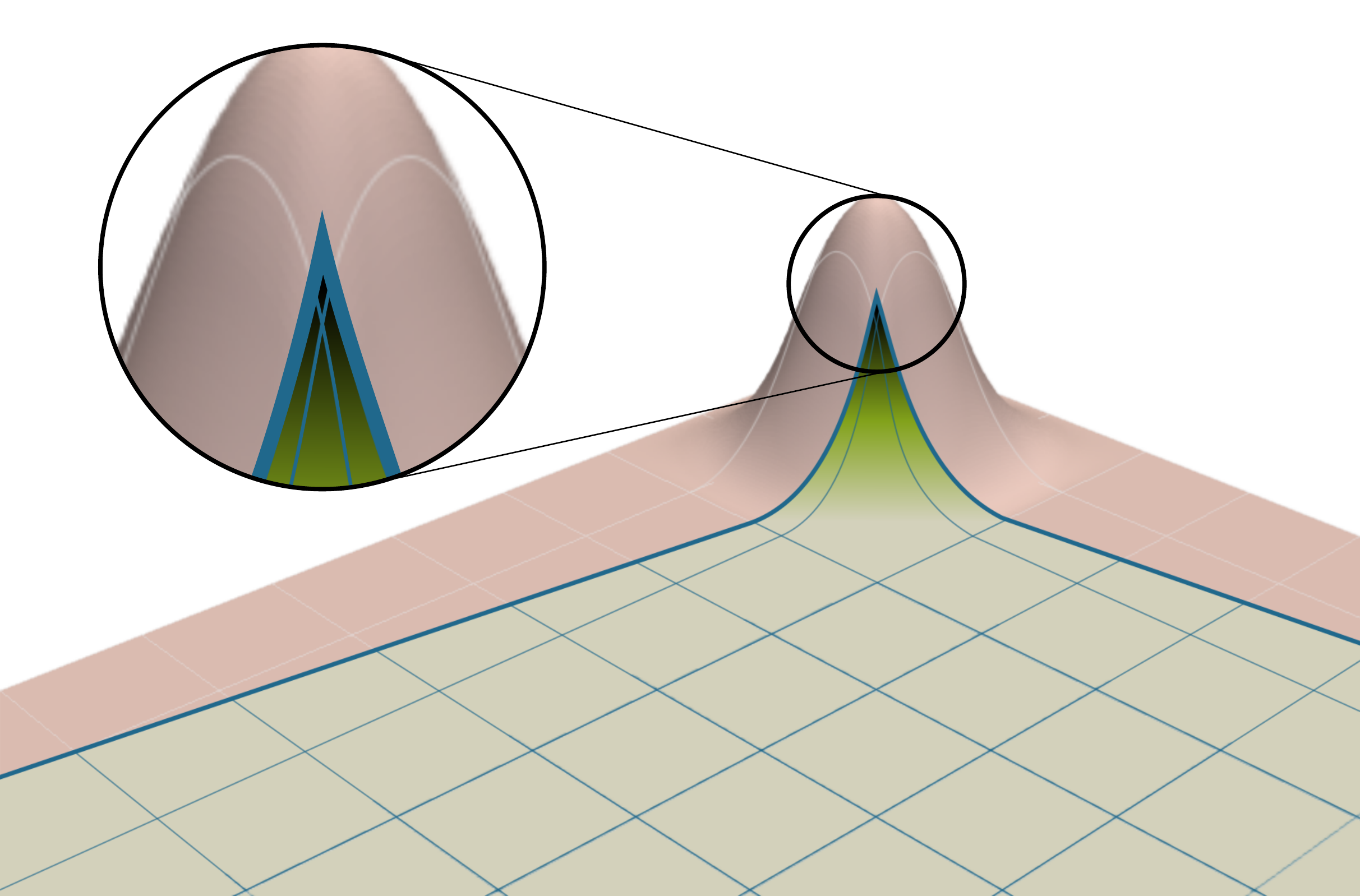}\label{fig:cut_funcs1b}}\hfill\\
    \subfloat[First sliver-cut function, \\\centering$\eigvec= \begin{bmatrix}\cdots0,1,1,\cdots1,1,0,\cdots \end{bmatrix}^\mathrm{T}$.]{\includegraphics[trim=0 70 0 0, clip, width=0.42\linewidth]{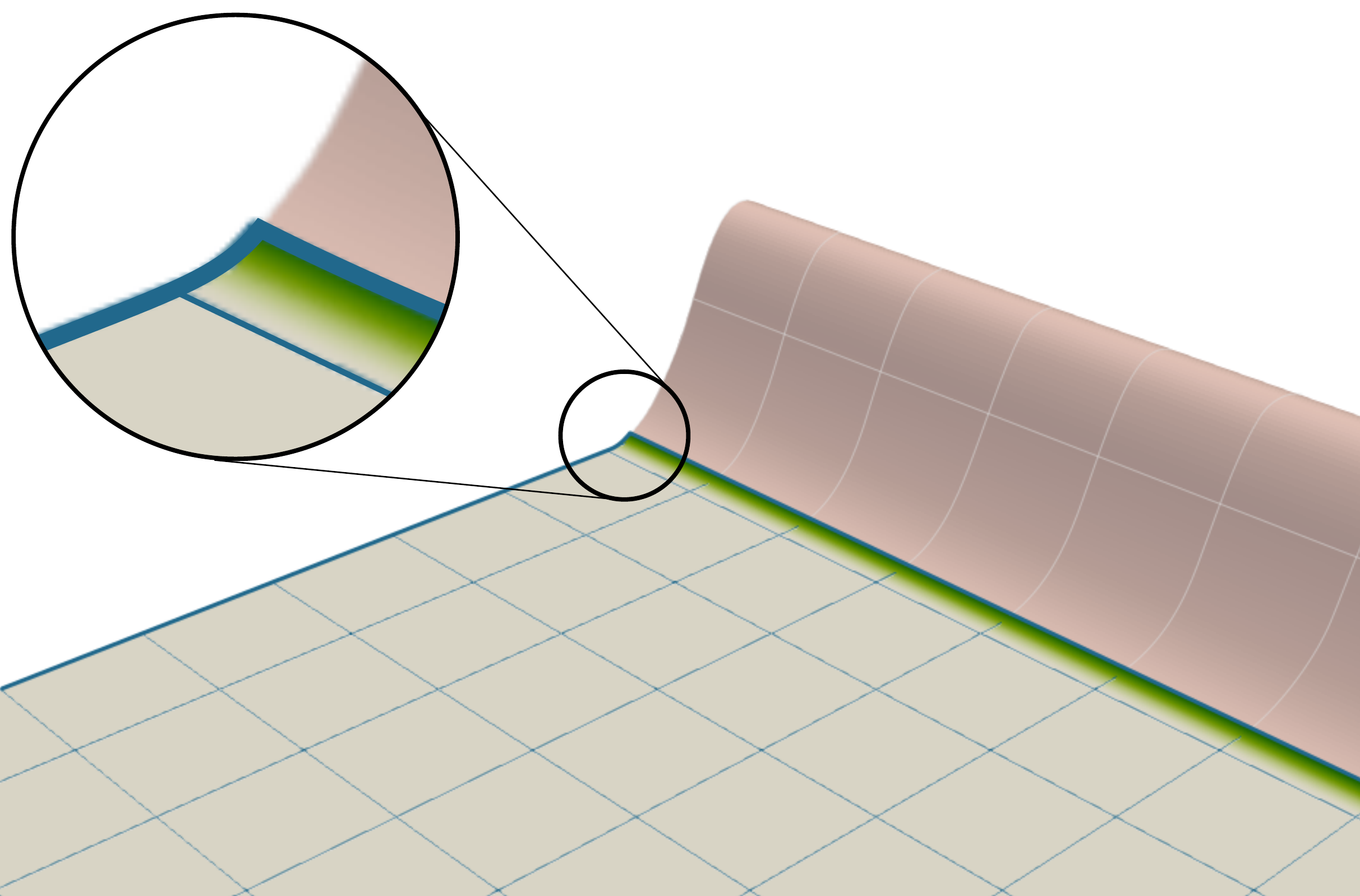}\label{fig:cut_funcs2}}\hspace{1.3cm}
    \subfloat[Second sliver-cut function, \\\centering$\eigvec= \begin{bmatrix}\cdots0,1,1,\cdots1,1,0,\cdots \end{bmatrix}^\mathrm{T}$.]{\includegraphics[trim=0 70 0 0, clip, width=0.42\linewidth]{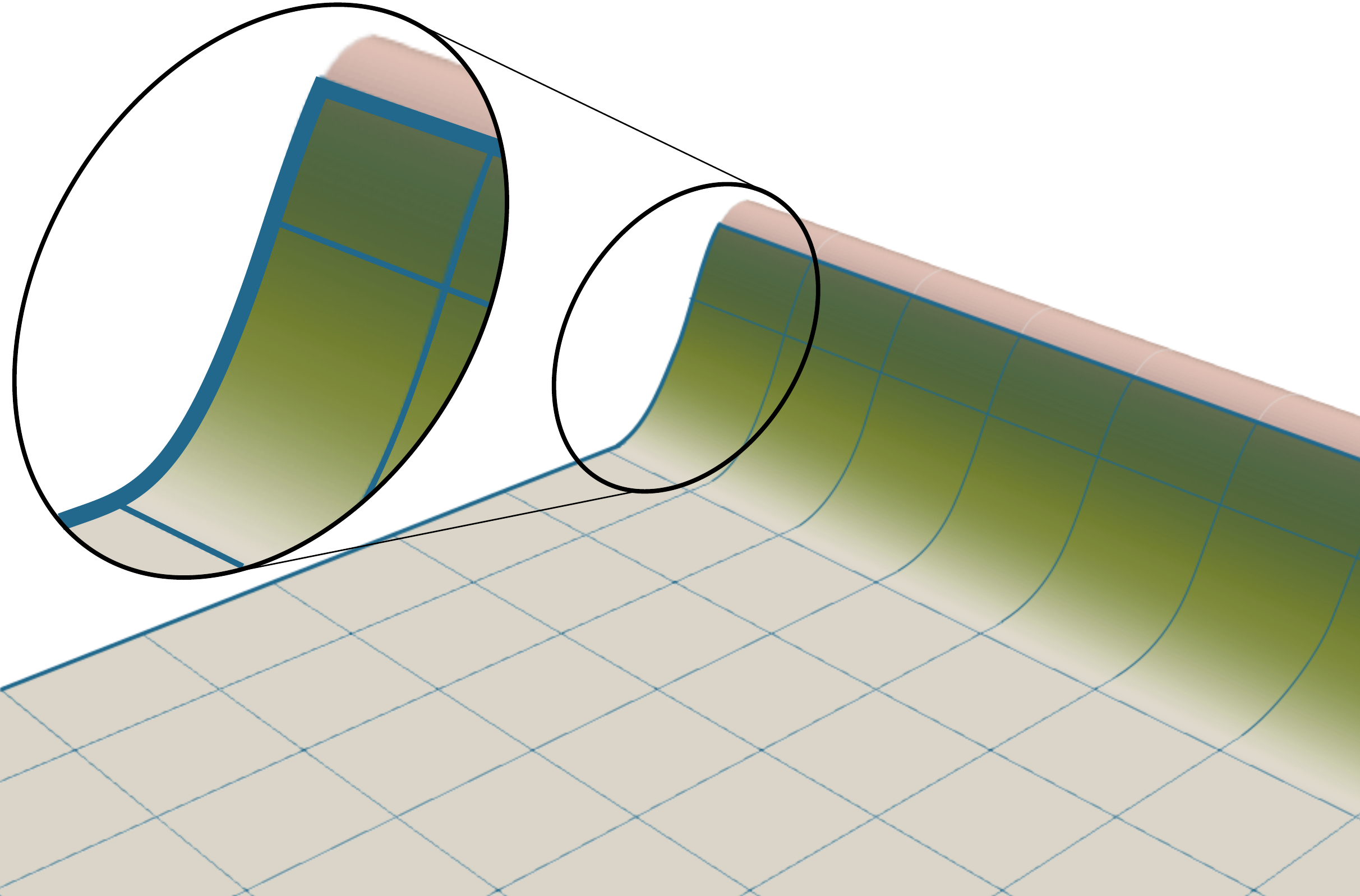}\label{fig:cut_funcs3}}\hfill \\
    \subfloat[First cut function in one dimension, \\\centering$\eigvec= \begin{bmatrix}\cdots 0,1,0,\cdots \end{bmatrix}^\mathrm{T}$.] {\includegraphics[width=0.49\linewidth]{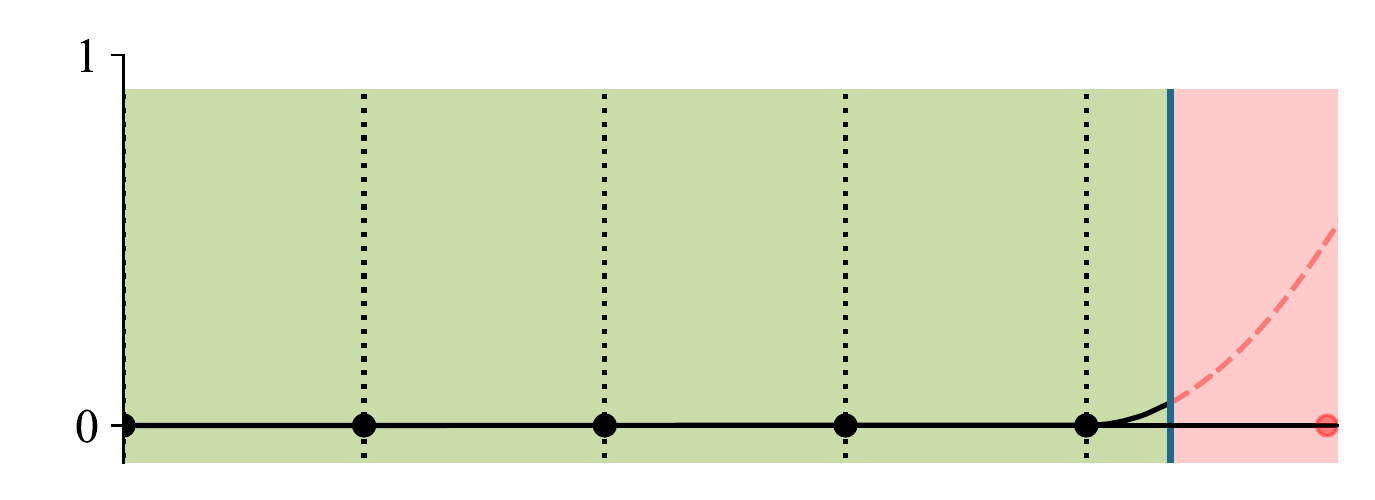}\label{fig:cut_funcs1D_1}}\hfill
    \subfloat[Second cut function in one dimension, \\\centering$\eigvec= \begin{bmatrix}\cdots 0,1,0,\cdots \end{bmatrix}^\mathrm{T}$.] {\includegraphics[width=0.49\linewidth]{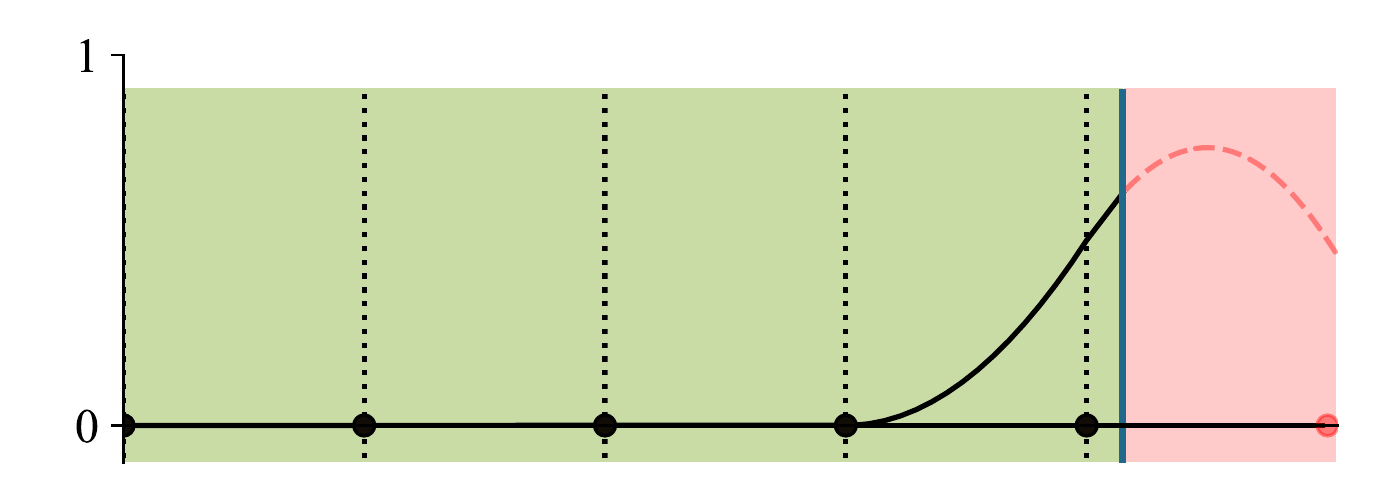}\label{fig:cut_funcs1D_2}} 
    \caption{Functions under consideration for cut-size scaling. Green represents the domain interior and red the domain exterior.}
    \label{fig:cut_funcs}
\end{figure}

Our study on the impact of the cut elements on the time-step estimate thus revolves around the examination of the generalized Rayleigh quotient for selected examples of problematic cut cases. In this regard, we consider two distinct types of cuts: corner cuts and sliver cuts. These may arise even when the background mesh is sufficiently refined to resolve small-scale geometric features. On these cut elements, we consider two different functions: one predominantly supported in the domain interior and one predominantly supported in the domain exterior. All four functions are depicted for the example of a two-dimensional domain in \cref{fig:cut_funcs1,fig:cut_funcs2,fig:cut_funcs3}. In one dimension, the four cases reduce to the two functions plotted in \cref{fig:cut_funcs1D_1,fig:cut_funcs1D_2}. For all functions considered, the coefficient vector consists exclusively of ones and zeros, whereby the condition of Remark~\ref{rmk:masslumpprob} is satisfied, and $\eigvec^\text{T}\,\mat{M}_D \, \eigvec $ may be evaluated as $M_D(\eigfunc)$.

\enlargethispage{0.5cm}

We determine the generalized Rayleigh quotient for different collections of components to the mass and stiffness matrix, i.e., for different finite element formulations. For the stiffness matrix, we consider a pure Neumann boundary, a pure penalty method with $\beta= \bbeta h_\DomEl^{-1}$, Nitsche's method with $\beta = \bbeta h_c^{-1} = \bbeta (\hfrac\, h_\DomEl)^{-1}$, and Nitsche's method combined with a ghost-penalty term and the choices $\beta = \bbeta h_\DomEl^{-1}$ and $\gammaK = \bgammaK h_{\DomEl}^{2p-1} $. For the mass matrix, we consider the consistent mass matrix of \cref{consistentmass}, a consistent mass matrix with ghost mass with $\gammaM = \bgammaM h_\DomEl^{2p+1} $, a lumped mass matrix per \cref{lumpedmass}, and a lumped mass matrix with (non-lumped) ghost mass. 

\enlargethispage{0.20cm}

Carrying out all the generalized Rayleigh quotients computations for the corner-cut functions of \cref{fig:cut_funcs1,fig:cut_funcs1b} results in the scalings detailed in \cref{tab:R1,tab:R1b}, respectively, and for the sliver functions depicted in \cref{fig:cut_funcs2,fig:cut_funcs3} results in \cref{tab:R2,tab:R3}, respectively. Cells in the tables marked in red indicate formulations and cut-cases for which the generalized Rayleigh quotient increases for decreasing cut-size, i.e., for $\hfrac\rightarrow +0$. Green cells signify formulations and cut-cases for which the generalized Rayleigh quotient decreases for decreasing cut-size, representing stable cases. The beige and yellow cells denote conditionally stable cases: there are restrictions on the polynomial order, spatial dimension and/or the penalty parameter for which the generalized Rayleigh quotient may or may not dominate and/or increase with decreasing cut-size.

\begin{table}[H]
\caption{Resulting scaling with cut size for the corner-cut function for the second-order problem.}\label{tab:R1}
\scriptsize
\begin{tabular}{|l|c|ccc|}\hline
                                                                              &                                                 			   &  \begin{tabular}[c]{@{}c@{}}Consistent mass  \\ (\cref{consistentmass}) \end{tabular}               & \begin{tabular}[c]{@{}c@{}} Lumped mass  \\ (\cref{lumpedmass}) \end{tabular}                      & \begin{tabular}[c]{@{}c@{}}Consistent/lumped  mass\\and ghost mass (\cref{ghostmass})\\ with $\gammaM = \bgammaM h^{2p+1}$ \end{tabular} \\\hline
                                                                              &    $\qquad\,\, \diagdown\mathcal{O}(\tilde{M})$ 			   & \multirow{2}{*}{$\hfrac^{2pd+d}$}  & \multirow{2}{*}{$\hfrac^{pd+d}$}     & \multirow{2}{*}{$\bgammaM\hfrac^{0}$}                                                              \\
                                                                              &               $\mathcal{O}(\tilde{K})\diagdown$ 			   &                                  &                                    &                                                                                           		  \\\hline
                                                                              &                                                				   & \redcell                         & \greencell                        & \greencell                                                                              		  \\[-0.1cm]
\begin{tabular}[c]{@{}l@{}}   Neumann (\cref{neumannprob}) \end{tabular}                           &                               $\hfrac^{2pd+d-2}$  			   & \redcell $\hfrac^{-2}$             & \greencell  $^*\,\,\hfrac^{pd-2} $         & \greencell $\bgammaM^{-1} \hfrac^{2pd+d-2}$       												  \\[0.3cm]
\begin{tabular}[c]{@{}l@{}}Penalty (\cref{penaltyform}) \\ $\quad \text{with } \beta = \bbeta h_K^{-1}$ \end{tabular}  &                   $\hfrac^{2pd+d-2}$ 			           & \redcell $\hfrac^{-2}$             & \greencell  $^*\,\,\hfrac^{pd-2}$         & \greencell $\bgammaM^{-1} \hfrac^{2pd+d-2}$   													  \\[0.3cm]
\begin{tabular}[c]{@{}l@{}}Nitsche (\cref{nitscheform}) \\ $\quad \text{with } \beta = \bbeta(\hfrac h_K)^{-1}$ \end{tabular} &              $\bbeta\hfrac^{2pd+d-2}$  			   & \redcell $\bbeta\hfrac^{-2}$       & \yellowcell  $^{**}\,\,\bbeta\hfrac^{pd-2}$   & \greencell $\bbeta\bgammaM^{-1} \hfrac^{2pd+d-2}$                                                  \\[0.3cm]
\begin{tabular}[c]{@{}l@{}}Nitsche and ghost (\cref{ghostform})\\  \,\,\,\, $\beta = \bbeta h_K^{-1}$, $\gammaK = \bgammaK h^{2p-1}$ \end{tabular} & $ \bgammaK \hfrac^0 $   & \redcell $\bgammaK\hfrac^{-2pd-d}$ & \redcell $\bgammaK\hfrac^{-pd-d}$    & \greencell $\bgammaK\bgammaM^{-1} \hfrac^0$ 												          \\[0.3cm]\hline
\end{tabular}\\[0.2cm]
\raggedright $\hspace{2cm}^{*\phantom{*}}\,$: Unstable for $d=1$ and $p=1$.\\
\raggedright $\hspace{2cm}^{**}\,$: Unstable for $d=1$ and $p=1$, and may dominate for $d=2$ and $p=1$.
\end{table}

\begin{table}[H]
\caption{Resulting scaling with cut size for the second corner-cut function for the second-order problem.}\label{tab:R1b}
\scriptsize
\begin{tabular}{|l|c|ccc|}\hline
                                                                                                                                 &                                                  &  \begin{tabular}[c]{@{}c@{}}Consistent mass  \\ (\cref{consistentmass}) \end{tabular}               & \begin{tabular}[c]{@{}c@{}} Lumped mass  \\ (\cref{lumpedmass}) \end{tabular}                      & \begin{tabular}[c]{@{}c@{}}Consistent/lumped  mass\\and ghost mass (\cref{ghostmass})\\ with $\gammaM = \bgammaM h^{2p+1}$ \end{tabular} \\\hline
                                                                                                                                 &  $\qquad\,\, \diagdown\mathcal{O}(\tilde{M})$    & \multirow{2}{*}{$\hfrac^{0}$}           & \multirow{2}{*}{$\hfrac^{0}$}           & \multirow{2}{*}{$\hfrac^{0}$}                                                                       \\
                                                                                                                                 & $\mathcal{O}(\tilde{K})\diagdown$                &                                       &                                       &                                                                                                    \\\hline
                                                                                                                                 &                                                  &   \greencell                            & \greencell                           & \greencell                                                                                         \\[-0.1cm]
\begin{tabular}[c]{@{}l@{}}   Neumann (\cref{neumannprob}) \end{tabular}                                                                              & $\hfrac^{0}$                                       &   \greencell $\hfrac^{0}$                 & \greencell  $\hfrac^{0}$               & \greencell $\hfrac^{0}$                                                                             \\[0.3cm]
\begin{tabular}[c]{@{}l@{}}Penalty (\cref{penaltyform}) \\ $\quad \text{with } \beta = \bbeta h_K^{-1}$ \end{tabular}                                         & $\bbeta \hfrac^{0}$                                &   \yellowcell $^*\,\,\bbeta \hfrac^{0}$          & \yellowcell  $^*\,\,\bbeta \hfrac^{0}$        & \yellowcell $^*\,\,\bbeta \hfrac^{0}$                                                                       \\[0.3cm]
\begin{tabular}[c]{@{}l@{}}Nitsche (\cref{nitscheform}) \\ $\quad \text{with } \beta = \bbeta(\hfrac h_K)^{-1}$ \end{tabular}                                   & $\bbeta \hfrac^{d-2}$                                  &   \greencell $^{**}\bbeta \hfrac^{d-2}$            & \greencell  $^{**}\bbeta \hfrac^{d-2}$             & \greencell $^{**}\bbeta  \hfrac^{d-2}$                                                                         \\[0.3cm]
\begin{tabular}[c]{@{}l@{}}Nitsche and ghost (\cref{ghostform})\\  \,\,\,\, $\beta = \bbeta h_K^{-1}$, $\gammaK = \bgammaK h^{2p-1}$ \end{tabular}  & $\bbeta \hfrac^0 $                                & \greencell $^{***}\,\,\bbeta\hfrac^{0}$             & \greencell  $^{***}\,\,\bbeta\hfrac^{0}$         & \greencell $^{***}\,\,\bbeta \hfrac^0$                                                         \\[0.3cm] \hline
\end{tabular}\\[0.2cm]
\raggedright $\hspace{2cm}^{*\phantom{**}}\,$: Adverse scaling with $\bbeta$.\\
\raggedright $\hspace{2cm}^{**\phantom{*}}\,$: Unstable for $d=1$.\\
\raggedright $\hspace{2cm}^{***}\,$: Adverse scaling with $\bbeta$, but only small $\bbeta$ is required for sufficient accuracy.
\end{table}

\begin{table}[H]
\caption{Resulting scaling with cut size for the first sliver-cut function for the second-order problem.}\label{tab:R2}
\scriptsize
\begin{tabular}{|l|c|ccc|}\hline
                                                                              &                                                  &  \begin{tabular}[c]{@{}c@{}}Consistent mass  \\ (\cref{consistentmass}) \end{tabular}               & \begin{tabular}[c]{@{}c@{}} Lumped mass  \\ (\cref{lumpedmass}) \end{tabular}                      & \begin{tabular}[c]{@{}c@{}}Consistent/lumped  mass\\and ghost mass (\cref{ghostmass})\\ with $\gammaM = \bgammaM h^{2p+1}$ \end{tabular} \\\hline
                                                                              &    $\qquad\,\, \diagdown\mathcal{O}(\tilde{M})$  & \multirow{2}{*}{$\hfrac^{2p+1}$}   & \multirow{2}{*}{$\hfrac^{p+1}$} & \multirow{2}{*}{$\bgammaM\hfrac^{0}$}                                                \\
                                                                              &               $\mathcal{O}(\tilde{K})\diagdown$  &                                  &                               &                                                                                        \\\hline
                                                                              &                                                  &   \redcell                       & \orangecell                   & \greencell                                                                              \\[-0.1cm]
\begin{tabular}[c]{@{}l@{}}   Neumann (\cref{neumannprob}) \end{tabular}                           &                               $\hfrac^{2p-1}$      &   \redcell $\hfrac^{-2}$           & \orangecell  $\hfrac^{p-2}$     & \greencell $\bgammaM^{-1}\hfrac^{2p-1}$                                                   \\[0.3cm]
\begin{tabular}[c]{@{}l@{}}Penalty (\cref{penaltyform}) \\ $\quad \text{with } \beta = \bbeta h_K^{-1}$ \end{tabular}  &                        $\hfrac^{2p-1}$      &   \redcell $\hfrac^{-2}$           & \orangecell  $\hfrac^{p-2}$     & \greencell $\bgammaM^{-1}\hfrac^{2p-1}$                                                   \\[0.3cm]
\begin{tabular}[c]{@{}l@{}}Nitsche (\cref{nitscheform}) \\ $\quad \text{with } \beta = \bbeta(\hfrac h_K)^{-1}$ \end{tabular} &             $\bbeta\hfrac^{2p-1}$      &   \redcell $\bbeta\hfrac^{-2}$     & \orangecell  $\bbeta\hfrac^{p-2}$     & \greencell $\bbeta\bgammaM^{-1}\hfrac^{2p-1}$                                             \\[0.3cm]
\begin{tabular}[c]{@{}l@{}}Nitsche and ghost (\cref{ghostform})\\  \,\,\,\, $\beta = \bbeta h_K^{-1}$, $\gammaK = \bgammaK h^{2p-1}$ \end{tabular} & $ \bgammaK\hfrac^0 $ & \redcell $\bgammaK\hfrac^{-2p-1}$ & \redcell $\bgammaK\hfrac^{-p-1}$   & \greencell $\bgammaK\bgammaM^{-1}\hfrac^0$                                               \\[0.3cm]\hline
\end{tabular}
\end{table}

\begin{table}[H]
\caption{Resulting scaling with cut size for the second sliver-cut function for the second-order problem.}\label{tab:R3}
\scriptsize
\begin{tabular}{|l|c|ccc|}\hline
                                                                                                                                 &                                                  &  \begin{tabular}[c]{@{}c@{}}Consistent mass  \\ (\cref{consistentmass}) \end{tabular}               & \begin{tabular}[c]{@{}c@{}} Lumped mass  \\ (\cref{lumpedmass}) \end{tabular}                      & \begin{tabular}[c]{@{}c@{}}Consistent/lumped  mass\\and ghost mass (\cref{ghostmass})\\ with $\gammaM = \bgammaM h^{2p+1}$ \end{tabular} \\\hline
                                                                                                                                 &  $\qquad\,\, \diagdown\mathcal{O}(\tilde{M})$    & \multirow{2}{*}{$\hfrac^{0}$}           & \multirow{2}{*}{$\hfrac^{0}$}           & \multirow{2}{*}{$\hfrac^{0}$}                                                                       \\
                                                                                                                                 & $\mathcal{O}(\tilde{K})\diagdown$                &                                       &                                       &                                                                                                    \\\hline
                                                                                                                                 &                                                  &   \greencell                            & \greencell                           & \greencell                                                                                         \\[-0.1cm]
\begin{tabular}[c]{@{}l@{}}   Neumann (\cref{neumannprob}) \end{tabular}                                                                              & $\hfrac^{0}$                                       &   \greencell $\hfrac^{0}$                 & \greencell  $\hfrac^{0}$               & \greencell $\hfrac^{0}$                                                                             \\[0.3cm]
\begin{tabular}[c]{@{}l@{}}Penalty (\cref{penaltyform}) \\ $\quad \text{with } \beta = \bbeta h_K^{-1}$ \end{tabular}                                         & $\bbeta \hfrac^{0}$                                &   \yellowcell $^*\,\,\bbeta \hfrac^{0}$          & \yellowcell  $^*\,\,\bbeta \hfrac^{0}$        & \yellowcell $^*\,\,\bbeta \hfrac^{0}$                                                                       \\[0.3cm]
\begin{tabular}[c]{@{}l@{}}Nitsche (\cref{nitscheform}) \\ $\quad \text{with } \beta = \bbeta(\hfrac h_K)^{-1}$ \end{tabular}                                   & $\bbeta \hfrac^{-1}$                                  &   \redcell $\bbeta \hfrac^{-1}$            & \redcell  $\bbeta \hfrac^{-1}$             & \redcell $\bbeta  \hfrac^{-1}$                                                                         \\[0.3cm]
\begin{tabular}[c]{@{}l@{}}Nitsche and ghost (\cref{ghostform})\\  \,\,\,\, $\beta = \bbeta h_K^{-1}$, $\gammaK = \bgammaK h^{2p-1}$ \end{tabular}  & $\bbeta \hfrac^0 $                                & \greencell $^{**}\,\,\bbeta\hfrac^{0}$             & \greencell  $^{**}\,\,\bbeta\hfrac^{0}$         & \greencell $^{**}\,\,\bbeta \hfrac^0$                                                         \\[0.3cm] \hline
\end{tabular}\\[0.2cm]
\raggedright $\hspace{2cm}^{*\phantom{*}}\,$: Adverse scaling with $\bbeta$.\\
\raggedright $\hspace{2cm}^{**}\,$: Adverse scaling with $\bbeta$, but only small $\bbeta$ is required for sufficient accuracy.
\end{table}

\begin{table}[H]
\caption{Overview of the stability characteristics for the second-order problem. Worst-case scaling of the cut scenarios considered in \cref{tab:R1,tab:R1b,tab:R2,tab:R3}.}\label{tab:Rconc}
\scriptsize
\begin{tabular}{|l|ccc|}\hline
                                                                              &  \begin{tabular}[c]{@{}c@{}}Consistent mass  \\ (\cref{consistentmass}) \end{tabular}               & \begin{tabular}[c]{@{}c@{}} Lumped mass  \\ (\cref{lumpedmass}) \end{tabular}                      & \begin{tabular}[c]{@{}c@{}}Consistent/lumped  mass\\and ghost mass (\cref{ghostmass})\\ with $\gammaM = \bgammaM h^{2p+1}$ \end{tabular} \\\hline
                                                                              &                            \redcell   & \orangecell                    & \greencell                                                                               \\[-0.1cm]
\begin{tabular}[c]{@{}l@{}}   Neumann (\cref{neumannprob}) \end{tabular}                           &                            \redcell   & \orangecell Unstable for $p=1$ & \greencell  Unnecessary for $p\geq 2$                                                      \\[0.3cm]
\begin{tabular}[c]{@{}l@{}}Penalty (\cref{penaltyform}) \\ $\quad \text{with } \beta = \bbeta h_K^{-1}$ \end{tabular}  &                            \redcell   & \orangecell \begin{tabular}[c]{@{}c@{}}Unstable for $p=1$\\ Scales with $\bbeta$  \end{tabular}  & \yellowcell    Scales with $\bbeta$                                                     \\[0.3cm]
\begin{tabular}[c]{@{}l@{}}Nitsche (\cref{nitscheform}) \\ $\quad \text{with } \beta = \bbeta (\hfrac h_K)^{-1}$ \end{tabular} &                      \redcell   & \redcell                       & \redcell                                                                                  \\[0.3cm]
\begin{tabular}[c]{@{}l@{}}Nitsche and ghost (\cref{ghostform})\\  \,\,\,\, $\beta = \bbeta h_K^{-1}$, $\gammaK = \bgammaK h^{2p-1}$ \end{tabular} & \redcell & \redcell                       & \greencell \begin{tabular}[c]{@{}c@{}}Requires $\bgammaM\sim \bgammaK$\\ and small $\bbeta$\end{tabular}                                                 \\[0.3cm] \hline
\end{tabular}
\end{table}

\newpage
\Cref{tab:Rconc} provides a comprehensive summary of findings that follow from the scaling relations detailed in \cref{tab:R1,tab:R1b,tab:R2,tab:R3}. The following conclusions stand out:
\begin{itemize}
    \item The first column of \cref{tab:Rconc} is fully colored red, expressing that the use of a consistent mass matrix without a form of mass scaling is not applicable for explicit immersed computation. A consistent mass matrix is unconventional for explicit analysis due to the implied cost of inversion, but in an immersed setting the largest eigenvalue may also become arbitrarily large for small cut elements, causing unfeasibly small critical time-step sizes. 
    \item The beige coloring of the first two rows of the second column in \cref{tab:Rconc} indicates that lumping the mass matrix can mitigate the problematic cut-size dependent scaling of the maximum eigenvalue, but only when the polynomial order of the basis function is at least quadratic. This was first observed in \cite{LeidingerPhD}.
    \item To enable the use of linear basis functions, a ghost mass must be added to ensure that the critical time-step sizes remain independent of the cut-size. This formulation corresponds to the last column in \cref{tab:Rconc}, which is the only column with green cells.
    \item When a penalty method is used for the enforcement of Dirichlet constraints, the maximum eigenvalue scales linearly with the non-dimensionalized penalty parameter $\bbeta$, as is indicated by the beige cells in the second row of \cref{tab:Rconc}. This is a known issue \cite{Leidinger2019}, and this undesirable scaling cannot be fixed by adding a ghost-mass term, nor by raising the polynomial order.
    \item A Nitsche formulation with local penalty parameter (i.e., without the ghost-stiffness term) does not yield a cut-size independent critical time-step size. The resulting scaling reminds of the degenerative error bounds due to sliver cuts proven in \cite{Prenter2018}.
    \item To adopt Nitsche's method in a stable manner, one requires both a ghost-mass term and as a ghost-stiffness term. Both penalty terms must scale with the uncut-element size $h_\DomEl$. To keep the maximum eigenvalue small, the ghost-mass penalty $\bgammaM$ must be of the same order of magnitude as the ghost-stiffness penalty $\bgammaK$, and the Nitsche penalty $\bbeta$ must be as small as possible. In contrast to a penalty method, using a small $\bbeta$ does not adversely impact solution quality, and hence we color this bottom-right box green. 
\end{itemize}

\begin{figure}[!t]
\vspace{-0.5cm}
    \centering
    \subfloat[A regular cut $p=2$ basis function.]{\includegraphics[width=0.48\linewidth]{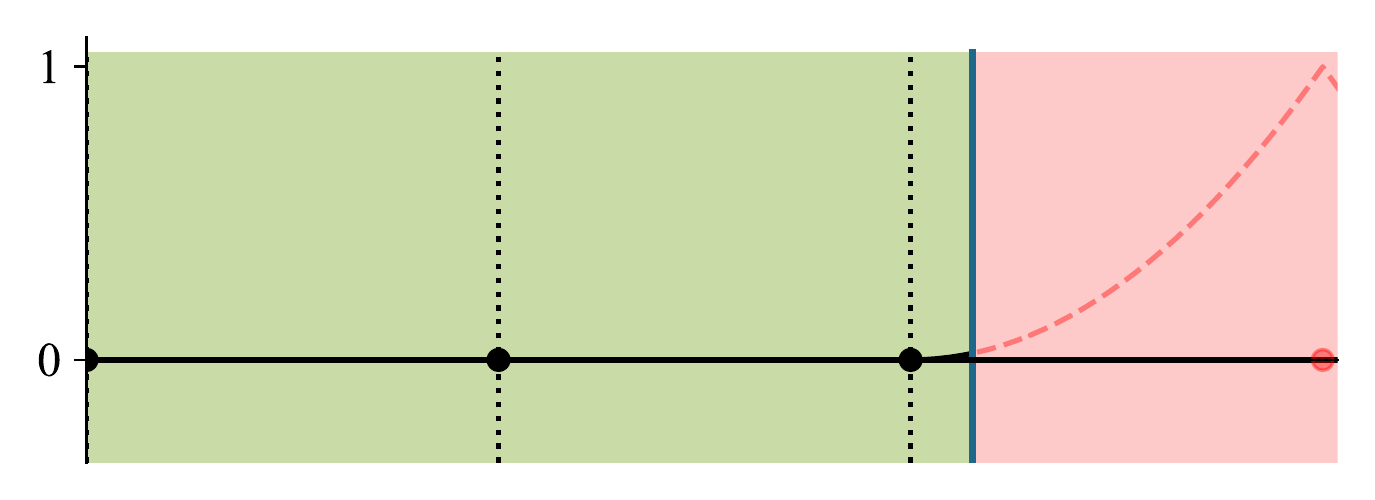}}
    \subfloat[\centering A $p=2$ cut basis function that approaches a cut $p=1$ function as $\hfrac\rightarrow +0$.]{\includegraphics[width=0.48\linewidth]{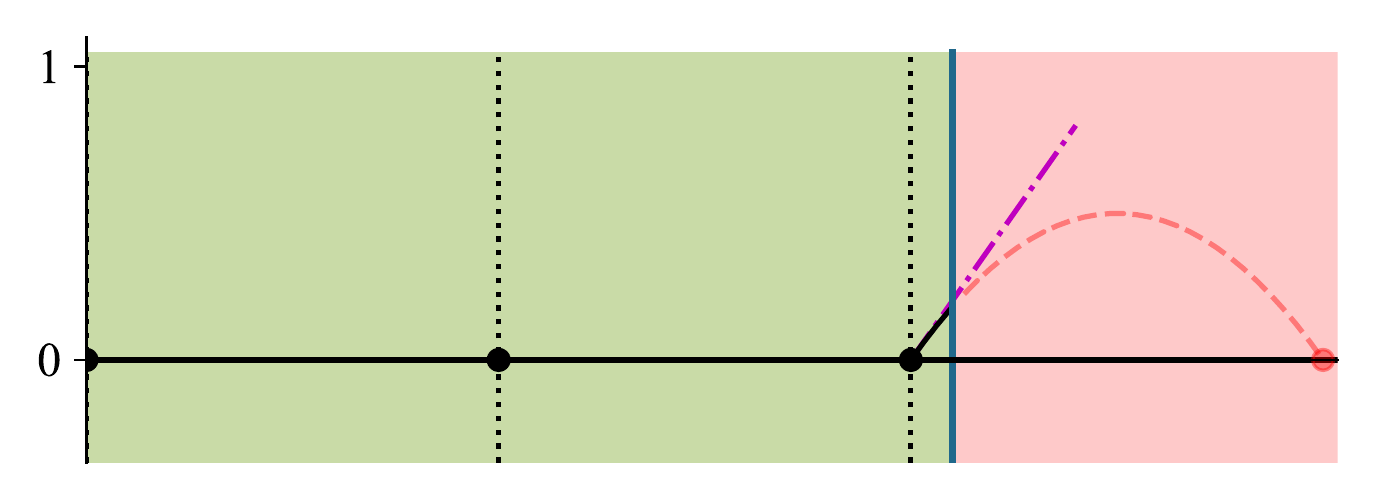}}
    \caption{Cut basis functions for a quadratic $\mathcal{C}^0$-continuous B-spline.}
    \label{fig:polapprox}
\end{figure}

\begin{rmk}\label{rmk:cont}
    The utilization of maximum regularity B-spline basis functions in the analysis reduces the number of cut functions that need to be considered. As demonstrated in the subfigures of the left column of \cref{fig:basisfunc}, all basis functions are identical up to affine transformation. 
    The adoption of lower-order continuous basis functions, such as the $\mathcal{C}^0$-continuous functions depicted in \cref{fig:bsplinebasis_c0p1,fig:bsplinebasis_c0p2,fig:bsplinebasis_c0p3}, necessitates the examination of a more extensive collection of cut functions. In particular, as $\hfrac\rightarrow +0$, the $\mathcal{C}^0$-continuous B-spline basis functions with small support locally approach lower-order polynomial functions, as is illustrated in \cref{fig:polapprox}. The lower-order polynomial functions are the critical cases in the first three rows of the middle column in \cref{tab:R1,tab:R3}. For lower-order continuous B-splines, the polynomial order in those scaling relations should thus be replaced by $k+1$, with $k$ the order of continuity. This corroborates the conclusion from \cite{LeidingerPhD} that a sufficiently high regularity is required to mitigate the cut-size dependency of the critical time-step size.
\end{rmk}

\begin{rmk}
    The scaling relations documented in \cref{tab:R2,tab:R3} are precisely those of \cref{tab:R1,tab:R1b}, with $d=1$. This concurrence is anticipated, as the cut functions depicted in \cref{fig:cut_funcs2,fig:cut_funcs3} are effectively one-dimensional. Thus, for an arbitrary $d$-dimensional domain, one only needs to determine the scaling relations for the $d$-dimensional cut functions (e.g., those in \cref{fig:cut_funcs1,fig:cut_funcs1b}), and the scaling relations for the lower-dimensional cut functions can then be deduced by successively replacing $d$ by the positive lower integer values below~$d$.
\end{rmk}

\section{Critical time-step analysis for a fourth-order problem}
\label{sec:FourthOrder}

While the wave equation studied in the previous section is a prototypical model equation for explicit analysis, its second-order nature limits the relevance of the conclusions from \cref{ssec:analysis2nd} when considering higher-order formulations, such as shell formulations. To understand how the conclusions translate, we repeat the analysis of \cref{sec:SecondOrder} in brief for the following fourth-order problem:
\begin{subequations}\label{ibvp4}
    \begin{alignat}{2}
    \rho\pddd{t}\phi+\Delta ( \kappa \Delta \scalar) &= f \qquad && \text{in }\Dom\times\TDom \subset \Re^d\times[0,T] \,,\\
    \nabla (\kappa \Delta \scalar)\cdot\normal &= q \qquad && \text{on }\bdyDom_q\times\TDom\,, \label{ibvp4_Nat2}\\[0.1cm]
    -\kappa \Delta \scalar &= m \qquad && \text{on }\bdyDom_m\times\TDom\,, \label{ibvp4_Nat1}\\[0.1cm]
    \nabla\scalar\cdot\normal &= g \qquad && \text{on }\bdyDom_N\times\TDom\,, \label{ibvp4_Ess2}\\[0.1cm]
    \scalar &= \scalar_D \qquad && \text{on }\bdyDirch\times\TDom\,, \label{ibvp4_Ess1}\\[0.1cm]
    \scalar &= \scalar_0 \qquad && \text{on }\Dom\times\{ 0 \}\,,\\[-0.1cm]
    \pdd{t}\scalar &= \dot{\scalar}_{0} \qquad && \text{on }\Dom\times\{ 0 \}\,.
    \end{alignat}
\end{subequations}
This equation describes the bending of an Euler-Bernoulli beam in one spatial dimension and is a shell-type analog in two dimensions \cite{Szilard2004,Vinson1990}. \Cref{ibvp4_Ess1,ibvp4_Ess2} are essential boundary conditions, describing displacement and rotation at the edge, respectively. \Cref{ibvp4_Nat1,ibvp4_Nat2} are natural conditions. A more direct link to shell formulations requires different forms of the natural boundary conditions, to represent the applied bending moment and shear force \cite{Vinson1990,Ying2021}, but the expressions of \cref{ibvp4_Nat1,ibvp4_Nat2} are applicable for arbitrary spatial dimension, making them more suitable for a general time-step stability analysis. Note that, due to the fourth-order parabolic nature of the PDE, two boundary conditions must be prescribed at any location on the boundary.  

\subsection{Semi-discrete formulation}

The weak formulation for the fourth-order problem of \cref{ibvp4} reads:
\begin{align}
&\text{For a.e. }t\in\TDom\text{, find } \scalar \in H^2_{\scalar_D,g}(\Dom) \text{ and } \pddd{t}\scalar \in H^{-2}(\Dom) \text{ s.t. } \forall\,  \scalartest \in H^2_{0,0}(\Dom): \nonumber\\
    &\quad \begin{cases} 
    \displaystyle \big\langle \pddd{t}\scalar ,  \rho\scalartest \big\rangle + \int\limits_\Dom  \kappa \Delta \scalar\, \Delta \scalartest \dDom = - \int\limits_{\bdyDom_q} \! q \, \scalartest \dbdy - \int\limits_{\bdyDom_m} \! m \, \nabla\scalartest\cdot\normal \dbdy  \,, \\
    \displaystyle \phi\big|_{t=0} = \phi_0 \, , \\[0.2cm]
    \displaystyle \pdd{t}\phi\big|_{t=0} = \dot{\scalar}_{0} \, ,
    \end{cases}  \label{weakform4}
\end{align}
with $ H^2_{\scalar_D,g}(\Dom)$ the $H^2(\Dom)$ Sobolev space for which each member satisfies the trace equalities from \cref{ibvp4_Ess2} and \cref{ibvp4_Ess1}, $H^{-2}(\Dom)$ is the corresponding dual space, and $ \big\langle \cdot ,  \cdot \big\rangle$ denotes the pairing between them. If $\pddd{t}\scalar \in L^2(\Dom)$, then $\big\langle \pddd{t}\scalar ,  \rho\scalartest \big\rangle = \int_\Dom  \rho\pddd{t}\scalar\,\scalartest \dDom $.

In an immersed setting, both essential constraints need to be enforced weakly. A semi-discrete formulation may again be written in the general form of:
\begin{align}
    \mat{\tilde{M}} \, \ddd{t}\numvec{\hat{\scalar}} + \mat{\tilde{K}} \, \numvec{\hat{\scalar}} = \numvec{\tilde{F}} \,, 
\end{align}
where, as in the second-order case, the mass matrix may or may not be lumped and may or may not contain a ghost-penalty term. The stiffness matrix, in its most generic form, consists of the following terms:
\begin{align}\label{K4}
    \begin{split}
         \big[\, \mat{K} \,\big]_{ij}& = \int_{\Dom} \kappa \Delta N_i  \Delta N_j  \d{\Dom}  + \int_{\IF_g} \! \kappa \gammaK \, \jump{\partial^{k+1}_n N_i } \jump{\partial^{k+1}_n N_j }  \dbdy   \\[-0.1cm]
         &\hspace{-0.75cm}
         + \int_\bdyDirch  \!\!\! \kappa \beta_\scalar\, N_i N_j \dbdy  + \int_{\bdyDom_N} \!\!\! \kappa \beta_g \, \nabla N_i\cdot\normal \,   \nabla N_j\cdot\normal  \dbdy \\[-0.1cm]
         &+ \int_\bdyDirch \!\!\! \nabla\big( \kappa \Delta N_i\big)\cdot\normal \,  N_j\dbdy  
         + \int_\bdyDirch \!\!\! \nabla \big( \kappa \Delta N_j\big) \cdot\normal \, N_i \dbdy  \\[-0.1cm]
         &\hspace{1cm}
          - \int_{\bdyDom_N} \!\!\! \nabla N_j \cdot \normal \, \kappa \Delta N_i\dbdy  - \int_{\bdyDom_N} \!\!\! \nabla N_i \cdot\normal  \, \kappa \Delta N_j \dbdy \,.
    \end{split}
\end{align}
The second term is the ghost-penalty term, the third and fourth term are penalty terms for the enforcement of the two essential boundary conditions, and the last four terms are the consistency and symmetry terms to complete the Nitsche formulations.

For the penalty formulation, the following parameter scalings are applicable:
\begin{subequations}
\begin{alignat}{3}
    \beta_\scalar \big|_\DomEl = \bbeta_\scalar h_\DomEl^{-3} \,, \\
    \beta_g \big|_\DomEl  = \bbeta_g h_\DomEl^{-1} \,. 
\end{alignat}
\end{subequations}
For the Nitsche formulation without ghost-penalty stabilization on the stiffness matrix, the positive definiteness (coercivity) requirement necessitates a minimum penalty value. We make use of:
\begin{subequations}
\begin{alignat}{3}
    &\beta_\scalar \big|_\DomEl = 3 \sup\limits_{\scalartest^h\in\scalarhspace} \frac{ \norm{\nabla(\Delta \scalartest^h) \cdot \normal }^2_{\bdyDirch\cap \DomEl}  }{ \norm{ \Delta \scalartest^h }^2_{\Dom\cap \DomEl}  } =: \bbeta_\scalar \big|_\DomEl\, (\hfrac\, h_\DomEl)^{-3} \,, \label{betac1}\\
    &\beta_g  \big|_\DomEl = 3 \sup\limits_{\scalartest^h\in\scalarhspace} \frac{ \norm{\Delta \scalartest^h }^2_{\bdyDom_N\cap \DomEl}  }{ \norm{ \Delta \scalartest^h }^2_{\Dom\cap \DomEl}  } =: \bbeta_g \big|_\DomEl\, (\hfrac\, h_\DomEl)^{-1}\,. \label{betac2}
\end{alignat}
\end{subequations}
When the stiffness matrix includes ghost-penalty stabilization, then we use the following parameters scalings, which satisfy the dimensional consistency requirement:
\begin{subequations}
\begin{alignat}{3}
    \gammaK \big|_\DomEl &= \bar{\gamma} \,  h_{\DomEl}^{2p-3} \\
    \beta_\scalar \big|_\DomEl &= \bbeta_\scalar h_\DomEl^{-3} \,, \\
    \beta_g \big|_\DomEl &= \bbeta_g h_\DomEl^{-1}\,.
\end{alignat}
\end{subequations}

\subsection{Analysis of the critical time-step size}
\label{ssec:analysis4th}
By following the methodology outlined in \cref{ssec:analysis2nd}, we derive the scaling relations for the four cut functions depicted in \cref{fig:cut_funcs}. The results for the various formulations are collected in \cref{tab:R1_B,tab:R1b_B} for the two corner-cut functions, and in \cref{tab:R2_B,tab:R3_B} for the two sliver-cut functions. The results for all four cut cases are combined and summarized in \cref{tab:Rconc_B}. The main conclusions drawn for the second-order equation carry over to this fourth-order problem: 
\begin{itemize}[itemsep=0.1pt]
    \item The use of a consistent mass matrix without mass scaling results in an undesirable cut-size dependent critical time-step size for all formulations. 
    \item Lumping the mass matrix is insufficient to yield a cut-size independent scheme for spline basis functions of order $p=2$ or even $p=3$. To ensure a cut-size independent critical time-step size, either $p=4$ basis function need to be used or ghost mass needs to be added.
    \item For both penalty formulations, the maximum eigenvalue scales linearly with the non-dimensionalized penalty parameter $\bbeta$, irrespective of the order of the basis functions and/or the use of mass-scaling.
    \item Both Nitsche formulations with local penalty parameters (i.e., without a ghost-stiffness term) are inapplicable.
    \item Both Nitsche formulations require a ghost-mass term as well as a ghost-stiffness term to ensure a cut-size independent critical time-step size.
\end{itemize}
The only difference compared to the results of the second-order problem, is that, in the case of mass lumping without mass scaling, this fourth-order problem requires at least quartic basis functions. For the analogous formulation for the second-order problem, quadratic basis functions were sufficient. This appears to hint at the more general condition $p\geq s$ (or, following \cref{rmk:cont}, $k+1\geq s$), with $s$ the order of the spatial partial differential operator.

\begin{table}[H]
\caption{Resulting scaling with cut size for the first corner-cut function for the fourth-order problem.}\label{tab:R1_B}
\scriptsize
\begin{tabular}{|l|c|ccc|}\hline
                                                                                                &                                                 			    &  \begin{tabular}[c]{@{}c@{}}Consistent mass  \\ (\cref{consistentmass}) \end{tabular}               & \begin{tabular}[c]{@{}c@{}} Lumped mass  \\ (\cref{lumpedmass}) \end{tabular}                      & \begin{tabular}[c]{@{}c@{}}Consistent/lumped  mass\\and ghost mass (\cref{ghostmass})\\ with $\gammaM = \bgammaM h^{2p+1}$ \end{tabular} \\\hline
                                                                                                &    $\qquad\,\, \diagdown\mathcal{O}(\tilde{M})$ 			    & \multirow{2}{*}{$\hfrac^{2pd+d}$}  & \multirow{2}{*}{$\hfrac^{pd+d}$}     & \multirow{2}{*}{$\bgammaM\hfrac^{0}$}                                                           \\
                                                                                                &               $\mathcal{O}(\tilde{K})\diagdown$ 			    &                                    &                                      &                                                                                           	  \\\hline
                                                                                                &                                                				& \redcell                           & \greencell                           & \greencell                                                                              		  \\[-0.1cm]
\begin{tabular}[c]{@{}l@{}}   Neumann \end{tabular}                                             &               $\hfrac^{2pd+d-4}$                    			& \redcell $\hfrac^{-4}$             & \greencell  $^*\,\,\hfrac^{pd-4} $   & \greencell $\bgammaM^{-1} \hfrac^{2pd+d-4}$       											  \\[0.3cm]
\begin{tabular}[c]{@{}l@{}}Penalty on $\scalar$ \\ $\quad \beta_\scalar = \bbeta_\scalar h_K^{-3}$ \end{tabular}        &               $\hfrac^{2pd+d-4}$ 			                    & \redcell $\hfrac^{-4}$             & \greencell  $^*\,\,\hfrac^{pd-4}$    & \greencell $\bgammaM^{-1} \hfrac^{2pd+d-4}$   												  \\[0.3cm]
\begin{tabular}[c]{@{}l@{}}Nitsche on $\scalar$ \\ $\quad \beta_\scalar = \bbeta_\scalar(\hfrac h_K)^{-3}$ \end{tabular}&               $\bbeta_\scalar\hfrac^{2pd+d-4}$  			            & \redcell $\bbeta_\scalar\hfrac^{-4}$       & \yellowcell  $^*\,\,\bbeta_\scalar\hfrac^{pd-4}$  & \greencell $\bbeta_\scalar\bgammaM^{-1} \hfrac^{2pd+d-4}$                                           \\[0.3cm]
\begin{tabular}[c]{@{}l@{}}Penalty on $\nabla\scalar\cdot\normal$ \\ $\quad \beta_g = \bbeta_g h_K^{-1}$ \end{tabular}        &               $\hfrac^{2pd+d-4}$ 			                    & \redcell $\hfrac^{-4}$             & \greencell  $^*\,\,\hfrac^{pd-4}$    & \greencell $\bgammaM^{-1} \hfrac^{2pd+d-4}$   												  \\[0.3cm]
\begin{tabular}[c]{@{}l@{}}Nitsche on $\nabla\scalar\cdot\normal$ \\ $\quad \beta_g = \bbeta_g(\hfrac h_K)^{-1}$ \end{tabular}&               $\bbeta_g\hfrac^{2pd+d-4}$  			            & \redcell $\bbeta_g\hfrac^{-4}$       & \yellowcell  $^*\,\,\bbeta_g\hfrac^{pd-4}$  & \greencell $\bbeta_g\bgammaM^{-1} \hfrac^{2pd+d-4}$                                           \\[0.3cm]
\begin{tabular}[c]{@{}l@{}}Nitsche on $\scalar$, $\beta_\scalar = \bbeta_\scalar h_K^{-3}$ \\ on $\nabla\scalar\cdot\normal$, $\beta_g = \bbeta_g(\hfrac h_K)^{-1}$ \\  and ghost $\gammaK = \bgammaK h^{2p-1}$ \end{tabular} & $ \bgammaK \hfrac^0 $          & \redcell $\bgammaK\hfrac^{-2pd-d}$ & \redcell $\bgammaK\hfrac^{-pd-d}$    & \greencell $\bgammaK\bgammaM^{-1} \hfrac^0$ 												      \\[0.3cm]\hline
\end{tabular}\\[0.2cm]
\raggedright $\hspace{2cm}^{*\phantom{*}}\,$: Unstable for $d=1$ and $p=2$.\\
\raggedright $\hspace{2cm}^{**}\,$: Unstable for $d=1$ and $p=2$, and may dominate for $d=2$ and $p=2$.
\end{table}

\begin{table}[H]
\caption{Resulting scaling with cut size for the second corner-cut function for the fourth-order problem.}\label{tab:R1b_B}
\scriptsize
\begin{tabular}{|l|c|ccc|}\hline
                                                                                                &                                                 			    &  \begin{tabular}[c]{@{}c@{}}Consistent mass  \\ (\cref{consistentmass}) \end{tabular}               & \begin{tabular}[c]{@{}c@{}} Lumped mass  \\ (\cref{lumpedmass}) \end{tabular}                      & \begin{tabular}[c]{@{}c@{}}Consistent/lumped  mass\\and ghost mass (\cref{ghostmass})\\ with $\gammaM = \bgammaM h^{2p+1}$ \end{tabular} \\\hline
                                                                                                                                 &  $\qquad\,\, \diagdown\mathcal{O}(\tilde{M})$    & \multirow{2}{*}{$\hfrac^{0}$}              & \multirow{2}{*}{$\hfrac^{0}$}             & \multirow{2}{*}{$\hfrac^{0}$}                                                                       \\
                                                                                                                                 & $\mathcal{O}(\tilde{K})\diagdown$                &                                            &                                           &                                                                                                     \\\hline
                                                                                                                                 &                                                  &   \greencell                               & \greencell                                & \greencell                                                                                          \\[-0.1cm]
\begin{tabular}[c]{@{}l@{}}   Neumann \end{tabular}                                                                              & $\hfrac^{0}$                                     &   \greencell $\hfrac^{0}$                  & \greencell  $\hfrac^{0}$                  & \greencell $ \hfrac^{0}$                                                                             \\[0.3cm]
\begin{tabular}[c]{@{}l@{}}Penalty on $\scalar$ \\ $\quad \beta_\scalar = \bbeta_\scalar h_K^{-3}$ \end{tabular}                 & $\bbeta_\scalar \hfrac^{0}$                      &   \yellowcell $^*\,\,\bbeta_\scalar \hfrac^{0}$  & \yellowcell  $^*\,\,\bbeta_\scalar \hfrac^{0}$  & \yellowcell $^*\,\,\bbeta_\scalar \hfrac^{0}$                                                                       \\[0.3cm]
\begin{tabular}[c]{@{}l@{}}Nitsche on $\scalar$ \\ $\quad \beta_\scalar = \bbeta_\scalar(\hfrac h_K)^{-3}$ \end{tabular}         & $\bbeta_\scalar \hfrac^{d-4}$                    &   \redcell $\bbeta_\scalar \hfrac^{d-4}$    & \redcell  $\bbeta_\scalar \hfrac^{d-4}$    & \redcell $\bbeta_\scalar \hfrac^{d-4}$                                                              \\[0.3cm]
\begin{tabular}[c]{@{}l@{}}Penalty on $\nabla\scalar\cdot\normal$ \\ $\quad \beta_g = \bbeta_g h_K^{-1}$ \end{tabular}           & $\bbeta_g \hfrac^{0}$                            &   \yellowcell $^{*}\,\,\bbeta_g \hfrac^{0}$  & \yellowcell  $^{*}\,\,\bbeta_g \hfrac^{0}$  & \yellowcell $^{*}\,\,\bbeta_g \hfrac^{0}$                                                                       \\[0.3cm]
\begin{tabular}[c]{@{}l@{}}Nitsche on $\nabla\scalar\cdot\normal$ \\ $\quad \beta_g = \bbeta_g(\hfrac h_K)^{-1}$ \end{tabular}   & $\bbeta_g \hfrac^{d-2}$                            &   \greencell $^{**}\,\,\bbeta_g \hfrac^{d-2}$         & \greencell  $^{**}\,\,\bbeta_g \hfrac^{d-2}$       & \greencell $^{**}\,\,\bbeta_g  \hfrac^{d-2}$                                                                   \\[0.3cm]
\begin{tabular}[c]{@{}l@{}}Nitsche on $\scalar$, $\beta_\scalar = \bbeta h_K^{-3}$ \\ on $\nabla\scalar\cdot\normal$, $\beta_g = \bbeta(\hfrac h_K)^{-1}$ \\  and ghost $\gammaK = \bgammaK h^{2p-1}$ \end{tabular}   & $\bbeta \hfrac^0 $  &   \greencell $^{***}\,\,\bbeta\hfrac^{0}$  & \greencell  $^{***}\,\,\bbeta\hfrac^{0}$ & \greencell $^{***}\,\,\bbeta \hfrac^0$                                                         \\[0.3cm] \hline
\end{tabular}\\[0.2cm]
\raggedright $\hspace{2cm}^{*}\phantom{^*\,\&\,^{***}}\,$: Adverse scaling with $\bbeta$.\\
\raggedright $\hspace{2cm}^{**}\phantom{\,\&\,^{***}}\,$: Unstable for $d=1$.\\
\raggedright $\hspace{2cm}^{**}\,\&\,^{***}\,$: Adverse scaling with $\bbeta$, but only small $\bbeta$ is required for sufficient accuracy.
\end{table}

\newpage
$$$$\\[-1cm]

\begin{table}[H]
\caption{Resulting scaling with cut size for the first sliver-cut function for the fourth-order problem.}\label{tab:R2_B}
\scriptsize
\begin{tabular}{|l|c|ccc|}\hline
                                                                              &                                                                         &  \begin{tabular}[c]{@{}c@{}}Consistent mass  \\ (\cref{consistentmass}) \end{tabular}               & \begin{tabular}[c]{@{}c@{}} Lumped mass  \\ (\cref{lumpedmass}) \end{tabular}                      & \begin{tabular}[c]{@{}c@{}}Consistent/lumped  mass\\and ghost mass (\cref{ghostmass})\\ with $\gammaM = \bgammaM h^{2p+1}$ \end{tabular} \\\hline
                                                                              &    $\qquad\,\, \diagdown\mathcal{O}(\tilde{M})$                         & \multirow{2}{*}{$\hfrac^{2p+1}$}    & \multirow{2}{*}{$\hfrac^{p+1}$}   & \multirow{2}{*}{$\bgammaM\hfrac^{0}$}                                                \\
                                                                              &               $\mathcal{O}(\tilde{K})\diagdown$                         &                                     &                                   &                                                                                        \\\hline
                                                                              &                                                                         &   \redcell                          & \orangecell                       & \greencell                                                                              \\[-0.1cm]
\begin{tabular}[c]{@{}l@{}}   Neumann \end{tabular}                           &                                $\hfrac^{2p-3}$                          &   \redcell $\hfrac^{-4}$            & \orangecell  $\hfrac^{p-4}$       & \greencell $\bgammaM^{-1}\hfrac^{2p-3}$                                                   \\[0.3cm]
\begin{tabular}[c]{@{}l@{}}Penalty on $\scalar$ \\ $\quad \beta_\scalar = \bbeta_\scalar h_K^{-3}$ \end{tabular}  &                    $\hfrac^{2p-3}$                          &   \redcell $\hfrac^{-4}$            & \orangecell  $\hfrac^{p-4}$       & \greencell $\bgammaM^{-1}\hfrac^{2p-3}$                                                   \\[0.3cm]
\begin{tabular}[c]{@{}l@{}}Nitsche on $\scalar$ \\ $\quad \beta_\scalar = \bbeta_\scalar(\hfrac h_K)^{-3}$ \end{tabular} &             $\bbeta_\scalar\hfrac^{2p-3} $                   &   \redcell $\bbeta_\scalar\hfrac^{-4}$      & \orangecell  $\bbeta_\scalar\hfrac^{p-4}$ & \greencell $\bbeta_\scalar\bgammaM^{-1}\hfrac^{2p-3}$                                             \\[0.3cm]
\begin{tabular}[c]{@{}l@{}}Penalty on $\nabla\scalar\cdot\normal$ \\ $\quad \beta_g = \bbeta_g h_K^{-1}$ \end{tabular}  &                    $\hfrac^{2p-3}$                          &   \redcell $\hfrac^{-4}$            & \orangecell * $\hfrac^{p-4}$, $\bbeta_g \hfrac^{p-3}$        & \greencell $\bgammaM^{-1}\hfrac^{2p-3}$                                                   \\[0.3cm]
\begin{tabular}[c]{@{}l@{}}Nitsche on $\nabla\scalar\cdot\normal$ \\ $\quad \beta_g = \bbeta_g(\hfrac h_K)^{-1}$ \end{tabular} &             $\bbeta_g\hfrac^{2p-3} $                   &   \redcell $\bbeta_g\hfrac^{-4}$      & \orangecell  $\bbeta_g\hfrac^{p-4}$ & \greencell $\bbeta_g\bgammaM^{-1}\hfrac^{2p-3}$                                             \\[0.3cm]
\begin{tabular}[c]{@{}l@{}}Nitsche on $\scalar$, $\beta_\scalar = \bbeta_\scalar h_K^{-3}$ \\ on $\nabla\scalar\cdot\normal$, $\beta_g = \bbeta_g(\hfrac h_K)^{-1}$ \\  and ghost $\gammaK = \bgammaK h^{2p-1}$ \end{tabular} & $ \bgammaK\hfrac^0 $   &   \redcell $\bgammaK\hfrac^{-2p-1}$ & \redcell $\bgammaK\hfrac^{-p-1}$  & \greencell $\bgammaK\bgammaM^{-1}\hfrac^0$                                               \\[0.3cm]\hline
\end{tabular}\\[0.2cm]
\raggedright $\hspace{2cm}^*\,$: Even though the $(p-4)$-scaling is more severe, the $(p-3)$-scaling may dominate for \\
\raggedright $\hspace{2.35cm}$ small $\hfrac$ and/or large $\bbeta_g$.\\
\end{table}

$$$$\\[-1.5cm]

\begin{table}[H]
\caption{Resulting scaling with cut size for the second sliver-cut function for the fourth-order problem.}\label{tab:R3_B}
\scriptsize
\begin{tabular}{|l|c|ccc|}\hline
                                                                                                                                 &                                                  &  \begin{tabular}[c]{@{}c@{}}Consistent mass  \\ (\cref{consistentmass}) \end{tabular}               & \begin{tabular}[c]{@{}c@{}} Lumped mass  \\ (\cref{lumpedmass}) \end{tabular}                      & \begin{tabular}[c]{@{}c@{}}Consistent/lumped  mass\\and ghost mass (\cref{ghostmass})\\ with $\gammaM = \bgammaM h^{2p+1}$ \end{tabular} \\\hline
                                                                                                                                 &  $\qquad\,\, \diagdown\mathcal{O}(\tilde{M})$    & \multirow{2}{*}{$\hfrac^{0}$}             & \multirow{2}{*}{$\hfrac^{0}$}           & \multirow{2}{*}{$\hfrac^{0}$}                                                                       \\
                                                                                                                                 & $\mathcal{O}(\tilde{K})\diagdown$                &                                           &                                         &                                                                                                    \\\hline
                                                                                                                                 &                                                  &   \greencell                              & \greencell                              & \greencell                                                                                         \\[-0.1cm]
\begin{tabular}[c]{@{}l@{}}   Neumann \end{tabular}                                                                              & $\hfrac^{0}$                                     &   \greencell $\hfrac^{0}$                 & \greencell  $\hfrac^{0}$                & \greencell $\hfrac^{0}$                                                                             \\[0.3cm]
\begin{tabular}[c]{@{}l@{}}Penalty on $\scalar$ \\ $\quad \beta_\scalar = \bbeta_\scalar h_K^{-3}$ \end{tabular}                                         & $\bbeta_\scalar \hfrac^{0}$                              &   \yellowcell $^*\,\,\bbeta_\scalar \hfrac^{0}$   & \yellowcell  $^*\,\,\bbeta_\scalar \hfrac^{0}$  & \yellowcell $^*\,\,\bbeta_\scalar \hfrac^{0}$                                                                       \\[0.3cm]
\begin{tabular}[c]{@{}l@{}}Nitsche on $\scalar$ \\ $\quad \beta_\scalar = \bbeta_\scalar(\hfrac h_K)^{-3}$ \end{tabular}                                 & $\bbeta_\scalar \hfrac^{-3}$                             &   \redcell $\bbeta_\scalar \hfrac^{-3}$           & \redcell  $\bbeta_\scalar \hfrac^{-3}$          & \redcell $\bbeta_\scalar  \hfrac^{-3}$                                                                         \\[0.3cm]
\begin{tabular}[c]{@{}l@{}}Penalty on $\nabla\scalar\cdot\normal$ \\ $\quad \beta_g = \bbeta_g h_K^{-1}$ \end{tabular}                                         & $\bbeta_g \hfrac^{0}$                              &   \yellowcell $^*\,\,\bbeta_g \hfrac^{0}$   & \yellowcell  $^*\,\,\bbeta_g \hfrac^{0}$  & \yellowcell $^*\,\,\bbeta_g \hfrac^{0}$                                                                       \\[0.3cm]
\begin{tabular}[c]{@{}l@{}}Nitsche on $\nabla\scalar\cdot\normal$ \\ $\quad \beta_g = \bbeta_g(\hfrac h_K)^{-1}$ \end{tabular}                                 & $\bbeta_g \hfrac^{-1}$                             &   \redcell $\bbeta_g \hfrac^{-1}$           & \redcell  $\bbeta_g \hfrac^{-1}$          & \redcell $\bbeta_g  \hfrac^{-1}$                                                                         \\[0.3cm]
\begin{tabular}[c]{@{}l@{}}Nitsche on $\scalar$, $\beta_\scalar = \bbeta h_K^{-3}$ \\ on $\nabla\scalar\cdot\normal$, $\beta_g = \bbeta(\hfrac h_K)^{-1}$ \\  and ghost $\gammaK = \bgammaK h^{2p-1}$ \end{tabular}   & $\bbeta \hfrac^0 $                      &   \greencell $^{**}\,\,\bbeta\hfrac^{0}$  & \greencell  $^{**}\,\,\bbeta\hfrac^{0}$ & \greencell $^{**}\,\,\bbeta \hfrac^0$                                                         \\[0.3cm] \hline
\end{tabular}\\[0.2cm]
\raggedright $\hspace{2cm}^{*\phantom{*}}\,$: Adverse scaling with $\bbeta$.\\
\raggedright $\hspace{2cm}^{**}\,$: Adverse scaling with $\bbeta$, but only small $\bbeta$ is required for sufficient accuracy.
\end{table}

\begin{table}[H]
\caption{Overview of the stability characteristics for the fourth-order problem. Worst-case scaling of the cut scenarios considered in \cref{tab:R1_B,tab:R2_B,tab:R3_B}.}\label{tab:Rconc_B}
\scriptsize
\begin{tabular}{|l|ccc|}\hline
                                                                              &  \begin{tabular}[c]{@{}c@{}}Consistent mass  \\ (\cref{consistentmass}) \end{tabular}               & \begin{tabular}[c]{@{}c@{}} Lumped mass  \\ (\cref{lumpedmass}) \end{tabular}                      & \begin{tabular}[c]{@{}c@{}}Consistent/lumped  mass\\and ghost mass (\cref{ghostmass})\\ with $\gammaM = \bgammaM h^{2p+1}$ \end{tabular} \\\hline
                                                                              &                            \redcell   & \orangecell                    & \greencell                                                                               \\[-0.1cm]
\begin{tabular}[c]{@{}l@{}}   Neumann \end{tabular}                           &                            \redcell   & \orangecell Unstable for $p\in\{2,3\}$ & \greencell                                                        \\[0.3cm]
\begin{tabular}[c]{@{}l@{}}Penalty on $\scalar$ \\ $\quad \beta_\scalar = \bbeta_\scalar h_K^{-3}$ \end{tabular}  &                            \redcell   & \orangecell \begin{tabular}[c]{@{}c@{}}Unstable for $p\in\{2,3\}$\\ Scales with $\bbeta$  \end{tabular}  & \yellowcell    Scales with $\bbeta_\scalar$                                                     \\[0.3cm]
\begin{tabular}[c]{@{}l@{}}Nitsche on $\scalar$ \\ $\quad \beta_\scalar = \bbeta_\scalar(\hfrac h_K)^{-3}$ \end{tabular} &                      \redcell   & \redcell                       & \redcell                                                                                  \\[0.3cm]
\begin{tabular}[c]{@{}l@{}}Penalty  on $\nabla\scalar\cdot\normal$ \\ $\quad \beta_g = \bbeta_g h_K^{-1}$ \end{tabular}  &                            \redcell   & \orangecell \begin{tabular}[c]{@{}c@{}}Unstable for $p\in\{2,3\}$\\ Scales with $\bbeta$  \end{tabular}  & \yellowcell    Scales with $\bbeta_g$                                                     \\[0.3cm]
\begin{tabular}[c]{@{}l@{}}Nitsche  on $\nabla\scalar\cdot\normal$ \\ $\quad \beta_g = \bbeta_g (\hfrac h_K)^{-1}$ \end{tabular} &                      \redcell   & \redcell                       & \redcell                                                                                  \\[0.3cm]
\begin{tabular}[c]{@{}l@{}}Nitsche on $\scalar$, $\beta_\scalar = \bbeta h_K^{-3}$ \\ on $\nabla\scalar\cdot\normal$, $\beta = \bbeta_g(\hfrac h_K)^{-1}$ \\  and ghost $\gammaK = \bgammaK h^{2p-1}$ \end{tabular} & \redcell & \redcell                       & \greencell \begin{tabular}[c]{@{}c@{}}Requires $\bgammaM\sim \bgammaK$\\ and small $\bbeta$\end{tabular}                                                 \\[0.3cm] \hline
\end{tabular}
\end{table}

\section{Numerical experiments}
\label{sec:NumExp}

In this section, we present the results of our numerical experiments, which are designed to accomplish two objectives. First, in \cref{ssec:numexpscaling}, we verify the derived scaling relations and conclusions from \cref{sec:SecondOrder,sec:FourthOrder}. Secondly, we assess whether the presence of the ghost mass negatively impacts the error behavior and convergence characteristics of the explicit scheme. To this end, we examine a linear vibrating drum in \cref{ssec:drum} and a linear Kirchhoff-Love shell in \cref{ssec:shell}.

\subsection{Verification of time-step size scaling}
\label{ssec:numexpscaling}

\begin{figure}[!t]
    \centering
    \includegraphics[width=0.4\linewidth,trim=80 20 68 18, clip,]{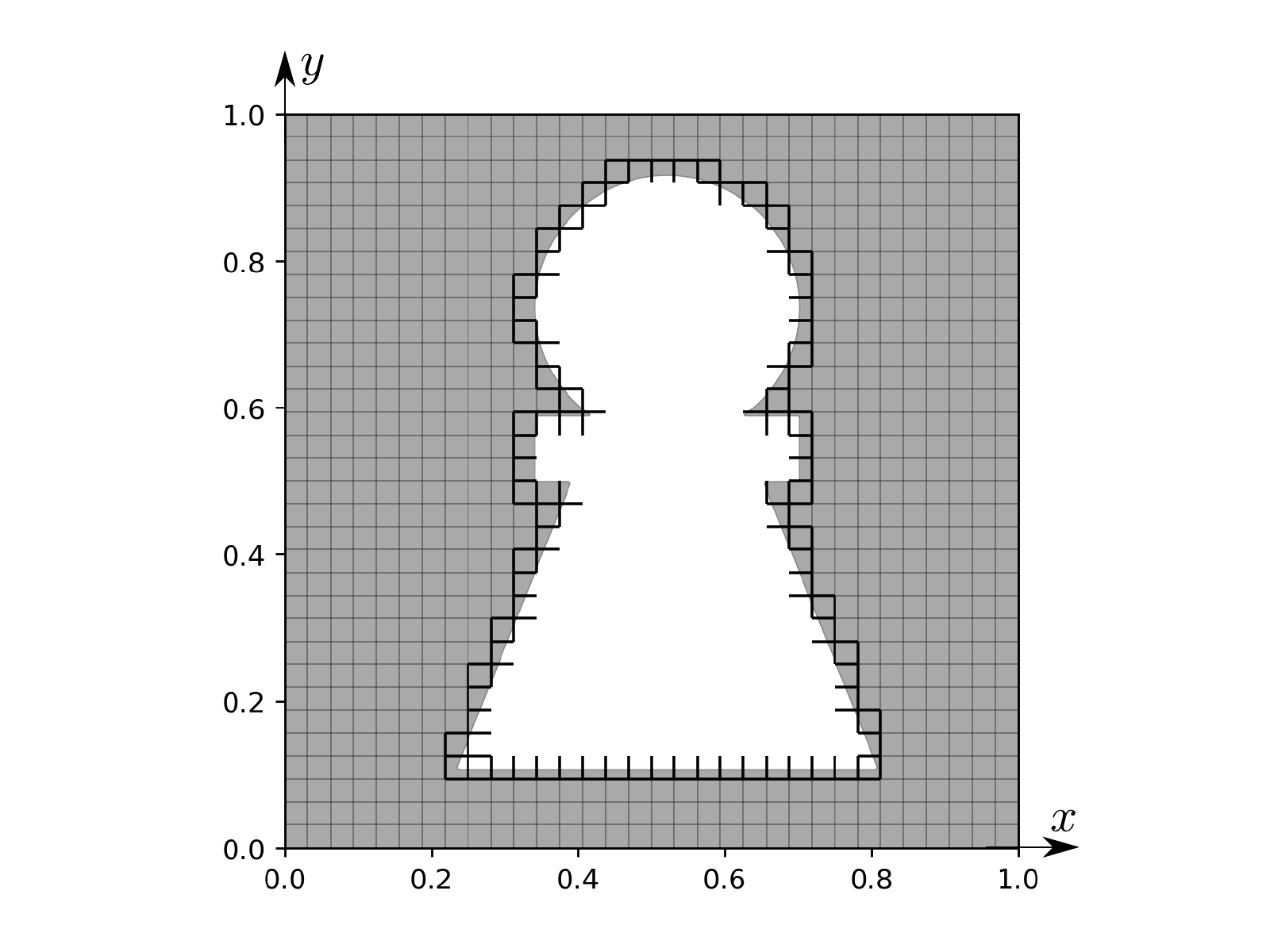}
    \caption{Domain, cut-out, mesh and ghost faces.}
    \label{fig:2Dcase}
\end{figure}

To numerically verify the scaling relations collected in \cref{tab:R1,tab:R1b,tab:R2,tab:R3,tab:R1_B,tab:R1b_B,tab:R2_B,tab:R3_B}, we consider the two-dimensional domain illustrated in \cref{fig:2Dcase}. The figure shows the mesh that is used for all subsequent computations, with the ghost faces highlighted. The domain cut-out involves straight edges, curved edges and positive and negative corner cuts, such that a wide variety of cut configurations may occur. To generate different cases, we randomly displace the cut-out within the domain by a distance between $-h_\DomEl$ and $h_\DomEl$ in the $x$ and $y$ directions. For each new domain, we compute the critical time-step, as defined in \cref{Dtcrit}, for the different formulations and polynomial orders. All computations are performed with row-sum mass-lumped mass matrices (apart from the ghost-mass term, as addressed in Remark~\ref{rmk:ghostmasslump}).

\subsubsection{Neumann boundaries}

\begin{figure}[!b]
\vspace{-0.5cm}
    \centering
    \subfloat[$p=1$.]{
        \begin{tikzpicture}
            \draw (0, 0) node[inner sep=0] {\includegraphics[width=0.48\linewidth]{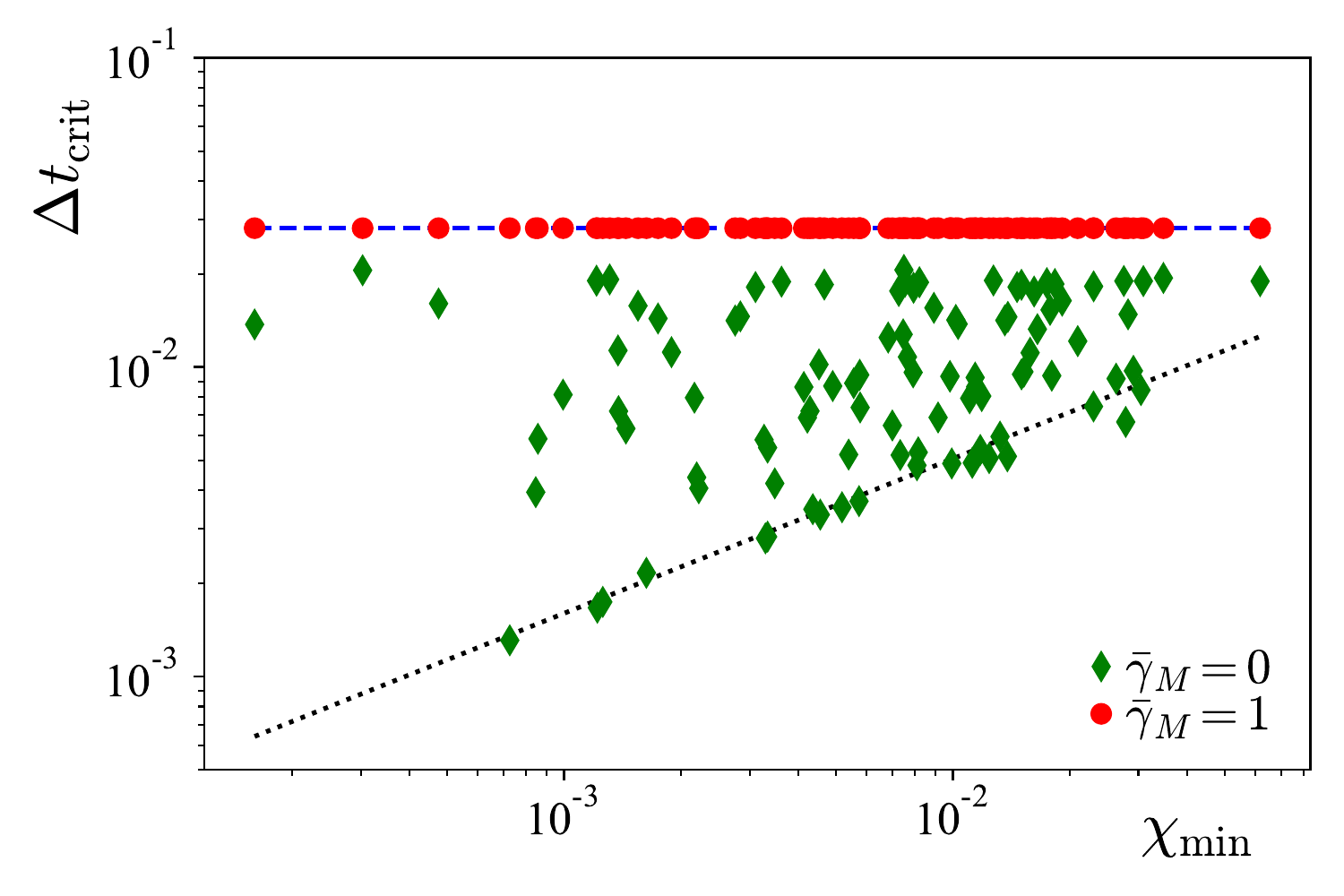}\label{fig:Dt2_Neumann_P1}};
            \draw (-1.8, -1.4) node [fill=white,inner sep=1pt] {$\frac{1}{2}$};
        \end{tikzpicture}}
    \subfloat[$p=2$.]{\includegraphics[width=0.48\linewidth]{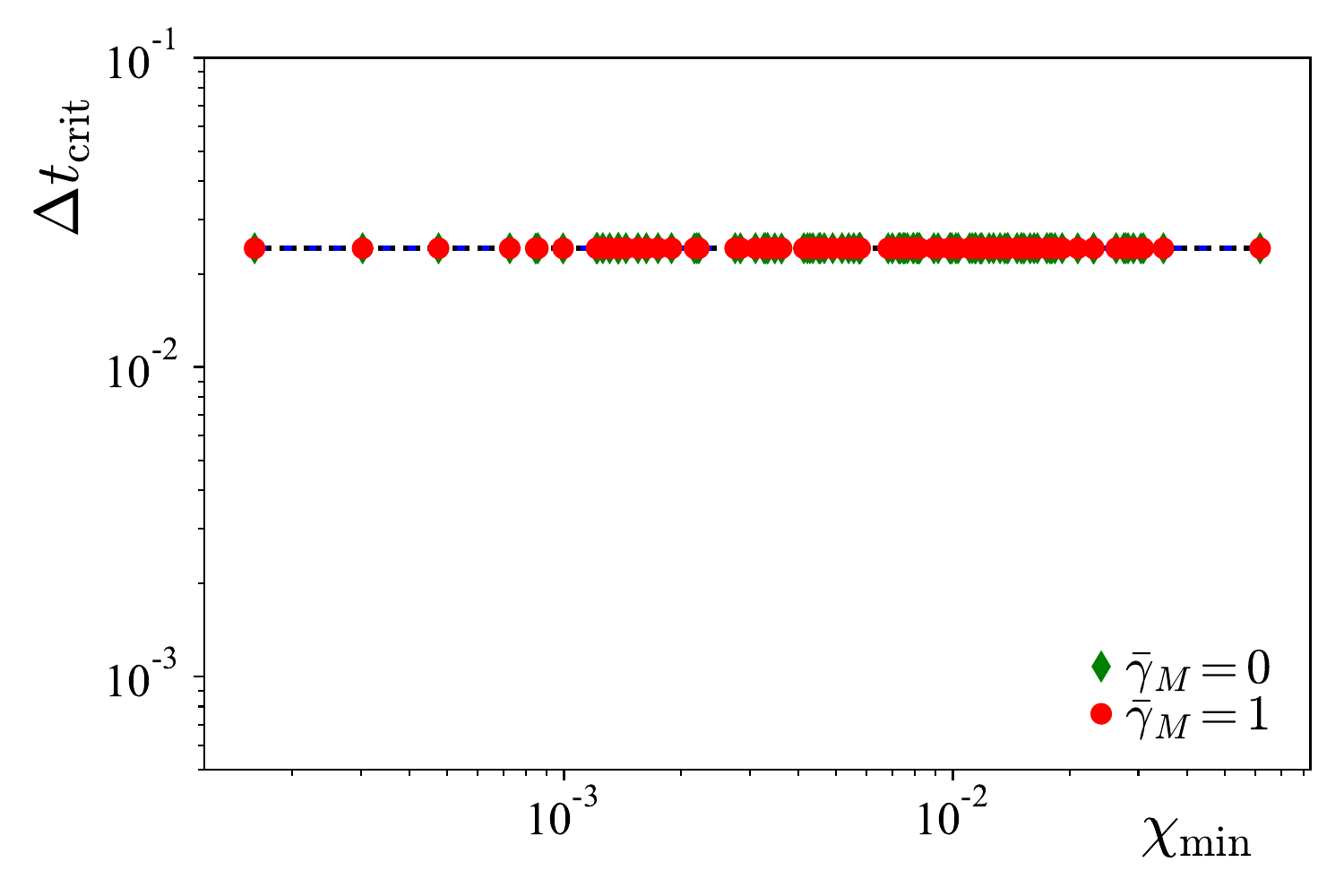}\label{fig:Dt2_Neumann_P2}}
    \caption{Neumann formulation of the second-order equation: first row in \cref{tab:R1,tab:R1b,tab:R2,tab:R3,tab:Rconc}. Critical time-step size dependency on the minimal element size fraction for 100 perturbations of the domain of \cref{fig:2Dcase}.}
    \label{fig:Dt2_Neumann}
\end{figure}

\begin{figure}[!t]
\vspace{1cm}
    \centering
    \subfloat[$p=2$.]{
        \begin{tikzpicture}
            \draw (0, 0) node[inner sep=0] {\includegraphics[width=0.48\linewidth]{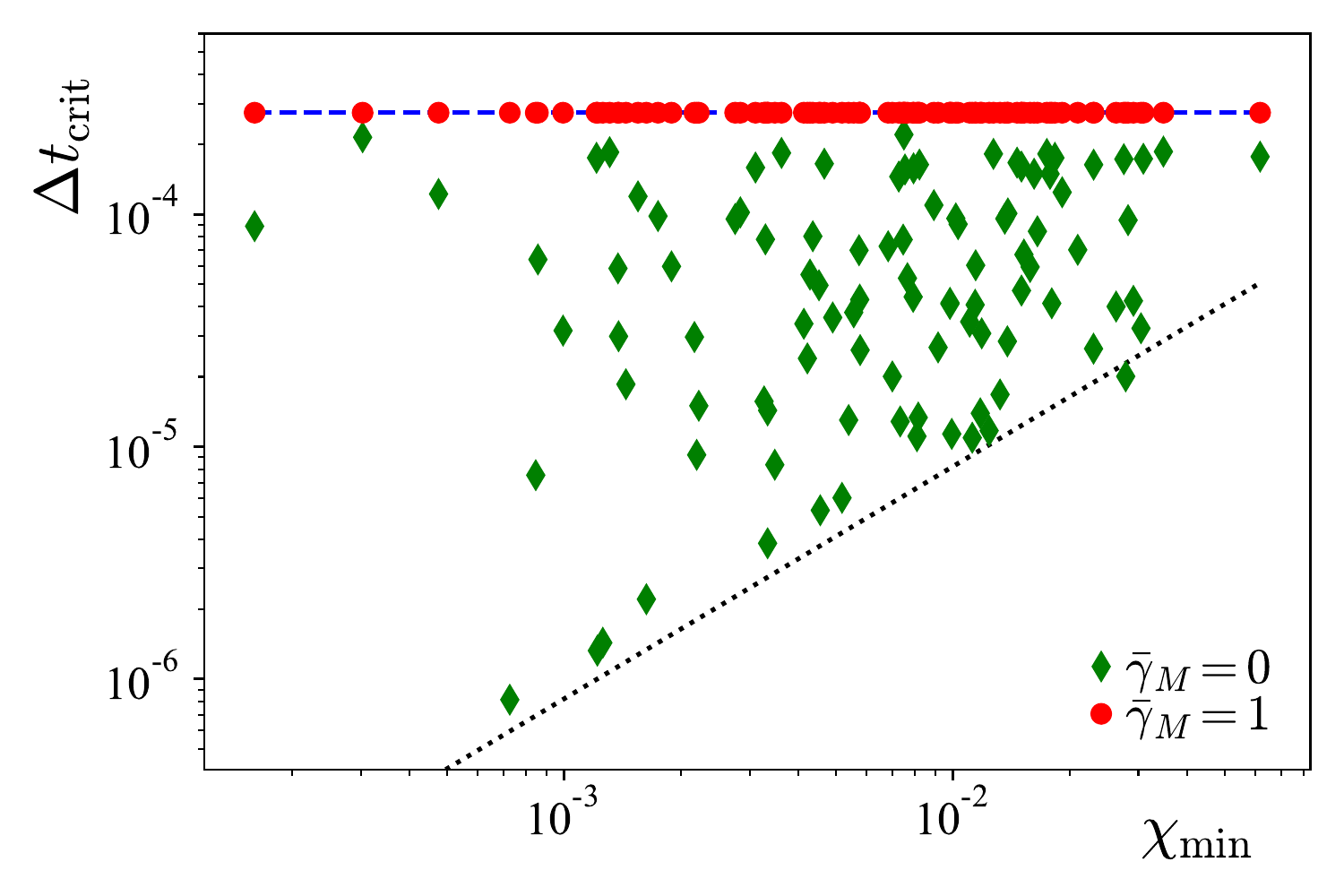}\label{fig:Dt4_Neumann_P2}};
            \draw (-0.65, -1.45) node [fill=white,inner sep=1pt] {$1$};
        \end{tikzpicture}}
    \subfloat[$p=3$.]{
        \begin{tikzpicture}
            \draw (0, 0) node[inner sep=0] {\includegraphics[width=0.48\linewidth]{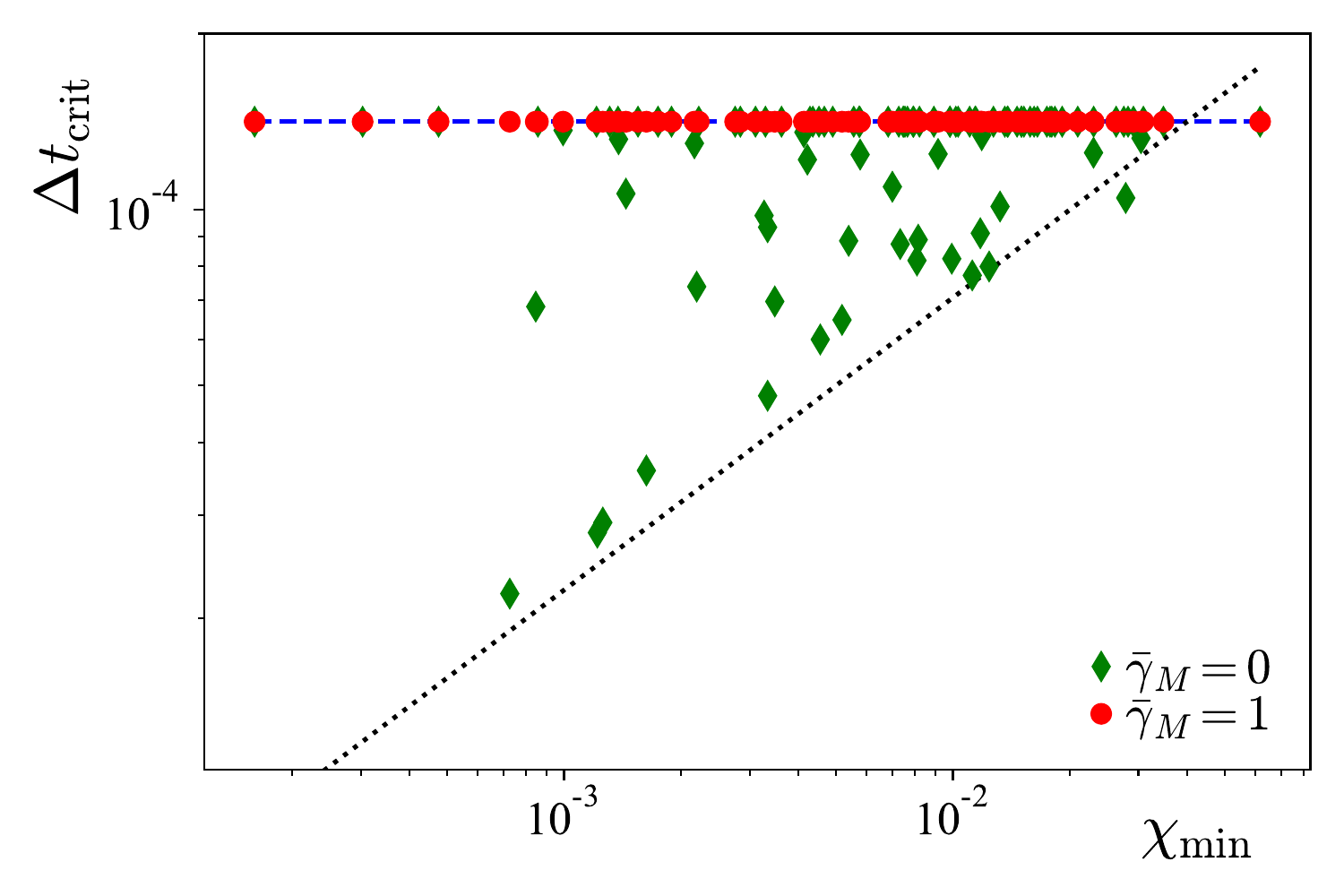}\label{fig:Dt4_Neumann_P3}};
            \draw (-1.45, -1.4) node [fill=white,inner sep=1pt] {$\frac{1}{2}$};
        \end{tikzpicture}}
    \caption{Neumann formulation of the fourth-order equation: first row in \cref{tab:R1_B,tab:R1b_B,tab:R2_B,tab:R3_B,tab:Rconc_B}. Critical time-step size dependency on the minimal element size fraction for 100 perturbations of the domain of \cref{fig:2Dcase}.}
    \label{fig:Dt4_Neumann}
\end{figure}

\Cref{fig:Dt2_Neumann,fig:Dt4_Neumann} present the results for the case where the entire cut-boundary is a Neumann boundary, for the second- and fourth-order equation, respectively. \Cref{fig:Dt2_Neumann_P1,fig:Dt4_Neumann_P2} involve basis functions with the lowest permitted polynomial order ($p=1$, respectively $p=2$), and in \cref{fig:Dt2_Neumann_P2,fig:Dt4_Neumann_P3} this order is incremented by one. The blue dashed lines in these figures correspond to the case without a cut-out, representing the optimal achievable value. We note that these `optimal' values do still suffer from the usual boundary outliers due to the repeating knots at the exterior boundary \cite{Cottrell2006,Deng2021,LeidingerPhD}, but believe that this represents the relevant comparison. The black dotted lines indicate the anticipated scaling relations derived in the preceding analysis sections. 

For a row-sum lumped mass matrix without ghost mass, the scaling of the maximum eigenvalue is predicted to be $\hfrac^{p-s}$, according to \cref{tab:R2,tab:R2_B}, where $s\in\{2,4\}$ is the order of the spatial differential operator. The critical time-step size relates to the maximum eigenvalue according to $\Delta t_\text{crit} \propto (\lambda_{\text{max}}^h)^{-\frac{1}{2}}$, and should thus scale as $\hfrac^{\frac{1}{2}(s-p)}$. As the black dotted lines indicate, this is indeed the trend that the bottom green diamonds follow. Not all data points follow this trend, since the smallest $\hfrac$ in the domain (on the horizontal axis) may correspond to a corner-cut function instead, which, according to \cref{tab:R1,tab:R1_B}, does not cause adverse scaling.

According to the results presented in \cref{tab:Rconc,tab:Rconc_B}, the cut-size dependency of the critical time-step size can be mitigated by raising the polynomial to $p\geq s$, or by incorporation the ghost-mass term. The latter approach is confirmed to be effective in all plots, as indicated by the red markers, and the former approach is demonstrated in \cref{fig:Dt2_Neumann_P2}. For unfavourably cut elements, both methods have the potential to increase the critical time-step size by orders of magnitude. For the second-order equation, quadratic basis functions are sufficient to achieve this effect, but quartic basis functions are required to eliminate the scaling for the fourth-order (shell-type) equation. In the event that the use of such high-order basis functions is not feasible, the ghost-mass term serves as an alternative solution.


\subsubsection{Penalty formulations}

\begin{figure}[H]
$$$$\\[-2cm]
    \centering
    \subfloat[$p=1$.]{
        \begin{tikzpicture}
            \draw (0, 0) node[inner sep=0] {\includegraphics[width=0.48\linewidth]{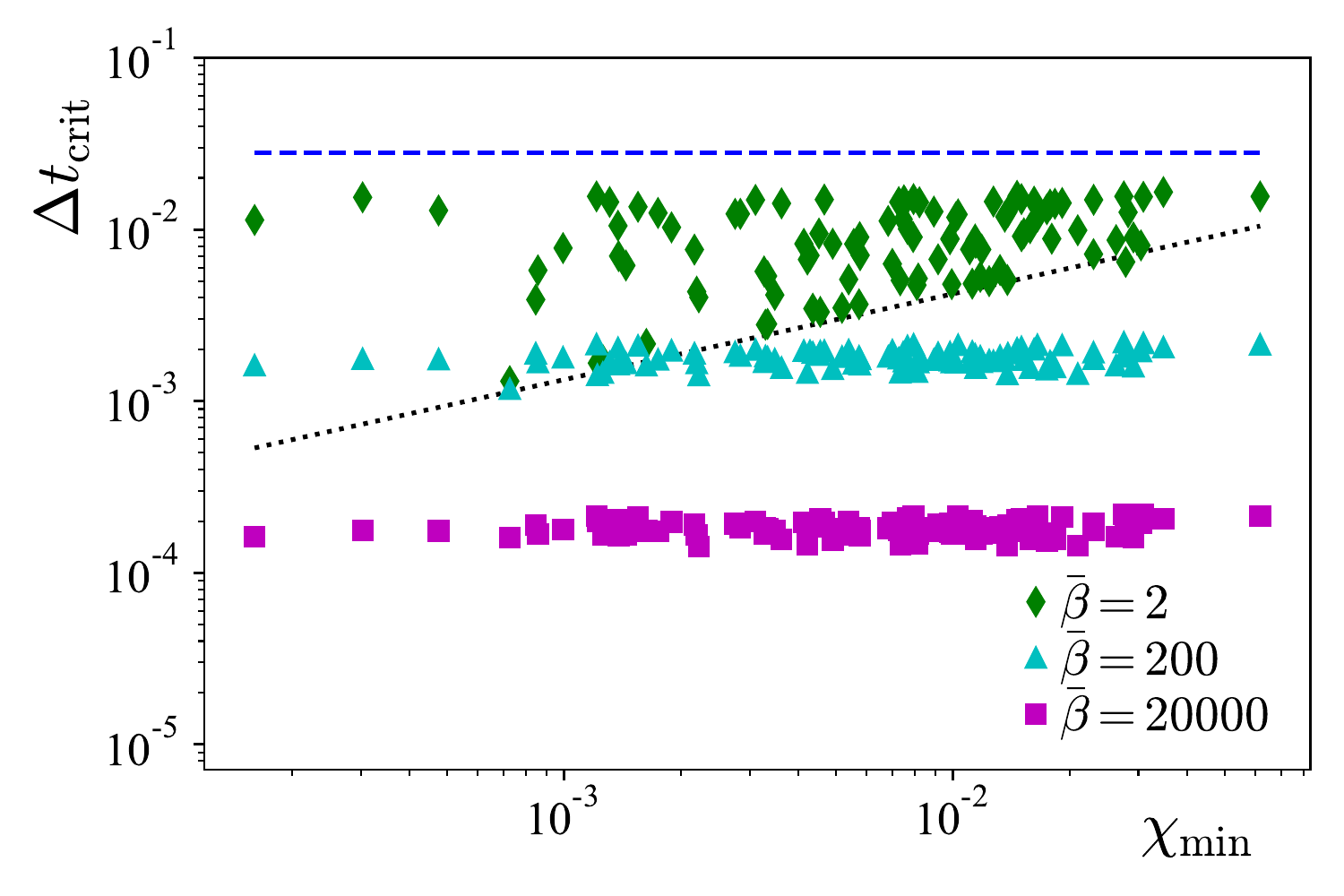}\label{fig:Dt2_penalty_P1}};
            \draw (-1.9, 0.1) node [fill=white,inner sep=1pt] {$\frac{1}{2}$};
        \end{tikzpicture}}
    \subfloat[$p=2$.]{\includegraphics[width=0.48\linewidth]{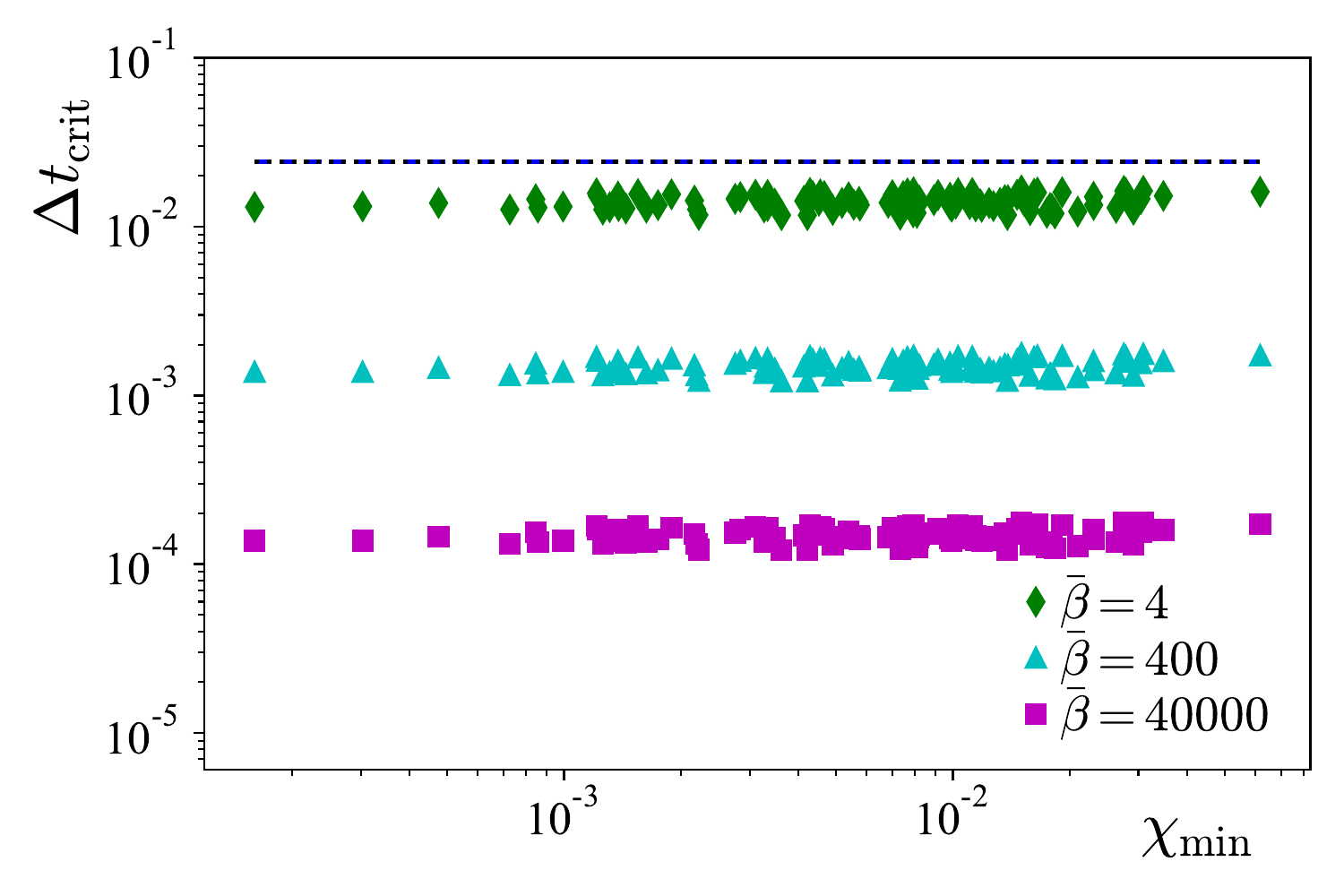}\label{fig:Dt2_penalty_P2}}
    \caption{Penalty formulations without ghost mass for the second-order problem: second row and second column in \cref{tab:R1,tab:R1b,tab:R2,tab:R3,tab:Rconc}. Critical time-step size dependency on the minimal element size fraction for 100 perturbations of the domain of \cref{fig:2Dcase}.}
    \label{fig:Dt2_penalty}
\end{figure}

\begin{figure}[H]
$$$$\\[-2cm]
    \centering
    \subfloat[On $\phi$, $p=2$. ]{
        \begin{tikzpicture}
            \draw (0, 0) node[inner sep=0] {\includegraphics[width=0.48\linewidth]{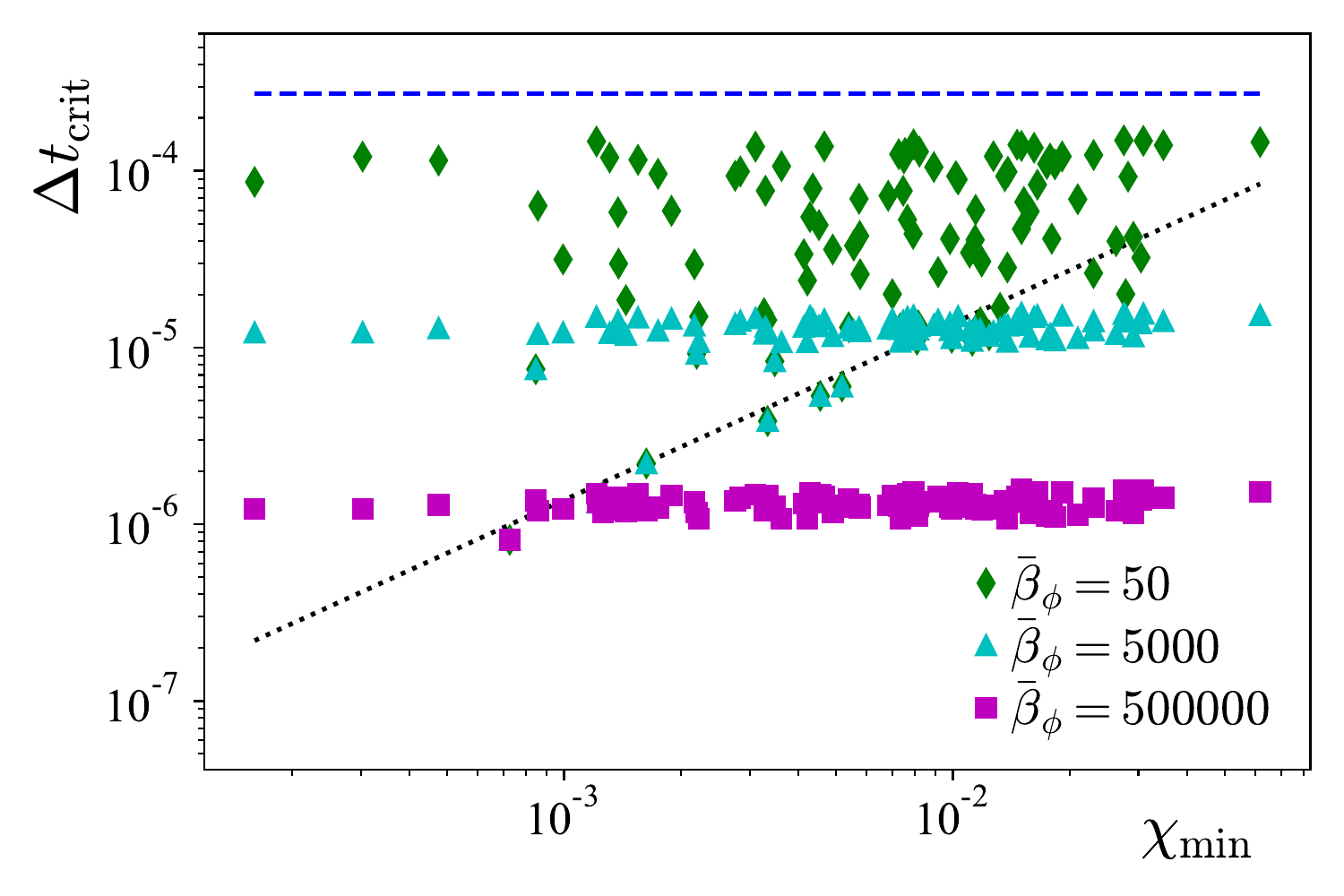}\label{fig:Dt4_penaltyPhi_P2}};
            \draw (-1.8, -0.85) node [fill=white,inner sep=1pt] {$1$};
        \end{tikzpicture}}
    \subfloat[On $\phi$, $p=3$.]{
        \begin{tikzpicture}
            \draw (0, 0) node[inner sep=0] {\includegraphics[width=0.48\linewidth]{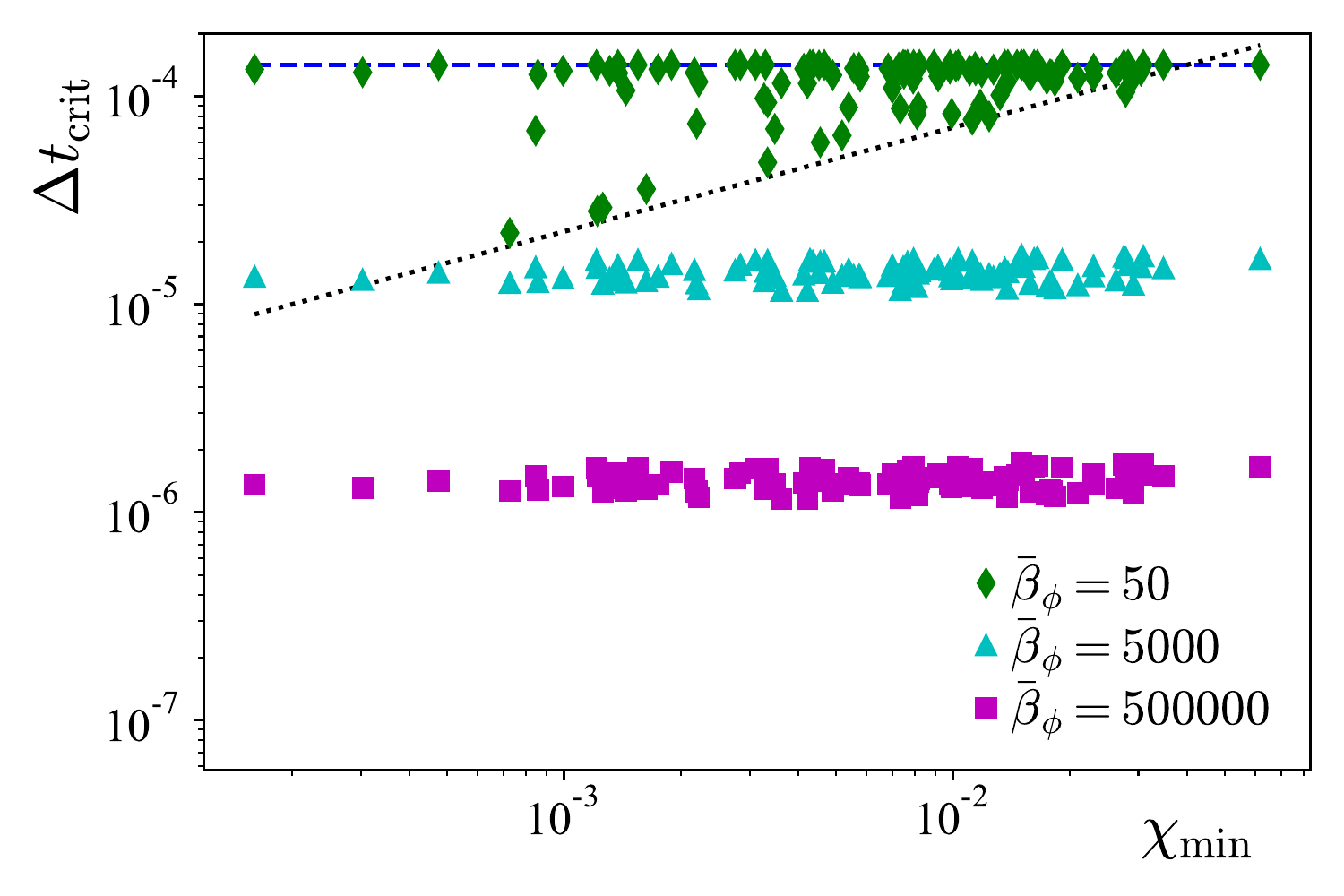}\label{fig:Dt4_penaltyPhi_P3}};
            \draw (-2.15, 0.9) node [fill=white,inner sep=1pt] {$\frac{1}{2}$};
        \end{tikzpicture}}\\
    \subfloat[On $\nabla\phi\cdot\normal$, $p=2$.]{
        \begin{tikzpicture}
            \draw (0, 0) node[inner sep=0] {\includegraphics[width=0.48\linewidth]{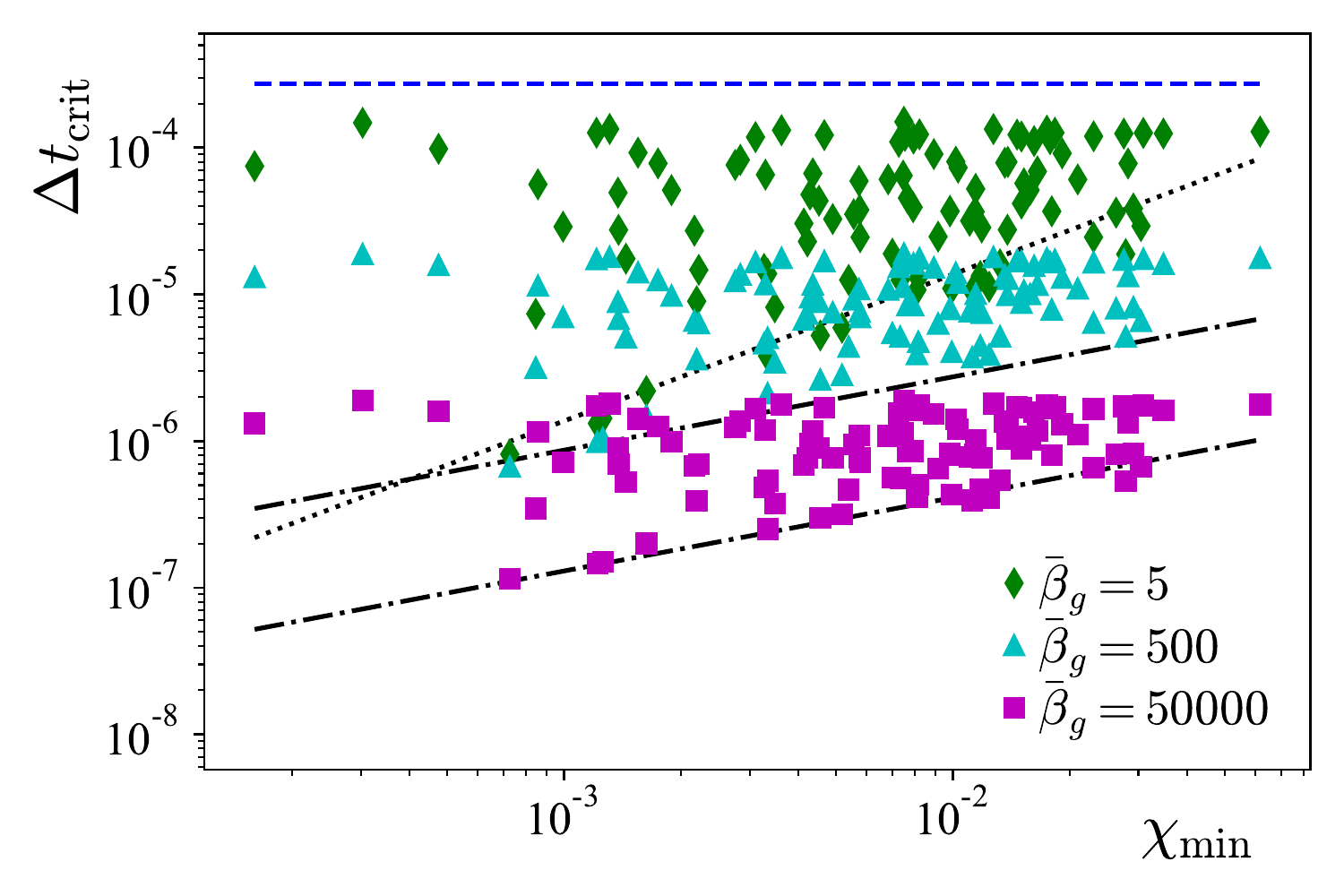}\label{fig:Dt4_penaltydPhi_P2}};
            \draw (-2.25, -0.55) node [fill=white,inner sep=1pt] {$1$};
            \draw (-1.5, -0.9) node [fill=white,inner sep=1pt] {$\frac{1}{2}$};
        \end{tikzpicture}}
    \subfloat[On $\nabla\phi\cdot\normal$, $p=3$.]{
        \begin{tikzpicture}
            \draw (0, 0) node[inner sep=0] {\includegraphics[width=0.48\linewidth]{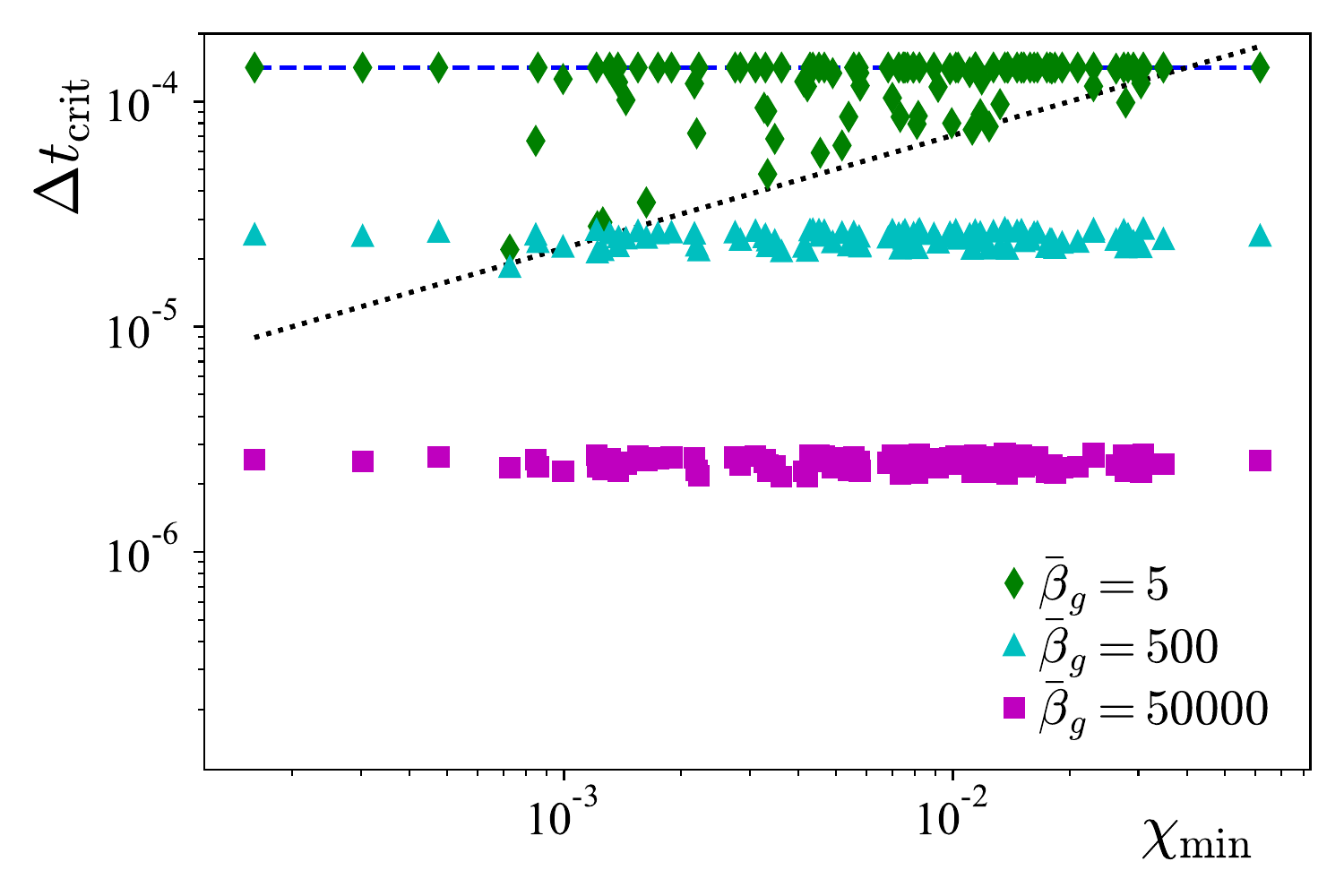}\label{fig:Dt4_penaltydPhi_P3}};
            \draw (-2.15, 0.7) node [fill=white,inner sep=1pt] {$\frac{1}{2}$};
        \end{tikzpicture}}
    \caption{Penalty formulations without ghost mass for the fourth-order problem: second and fourth row and second column in \cref{tab:R1_B,tab:R1b_B,tab:R2_B,tab:R3_B,tab:Rconc_B}. Critical time-step size dependency on the minimal element size fraction for 100 perturbations of the domain of \cref{fig:2Dcase}.}
    \label{fig:Dt4_penalty}
\end{figure}

If the cut-boundary is a Dirichlet boundary instead, and a penalty method is used for the enforcement of the constraints, then, according to \cref{tab:R3,tab:R3_B}, the maximum eigenvalue scales with $\bbeta$. This is confirmed in \cref{fig:Dt2_penalty,fig:Dt4_penalty} for basis functions with different polynomial orders for the second- and fourth-order problems, respectively. For $p\leq s$ ($s\in\{2,4\}$ the order of the spatial differential operator) the first sliver-cut function (\cref{tab:R2,tab:R2_B}) also introduces a scaling of the critical time-step with $\hfrac^{\frac{1}{2}(s-p)}$ in the same manner as for the Neumann boundary formulation. Both types of scaling are observed in the various subfigures: for $p<s$ and small $\bbeta$ the markers follow the black dotted line, but as $\bbeta$ is increased the second sliver-cut functions produce the highest eigenvalues and the dependency on $\hfrac$ is no-longer dominant in the considered collection of cut cases. 

We notice that the data-points in \cref{fig:Dt4_penaltydPhi_P2} are more scattered than those in the other figures, which can be attributed to the two types of scaling of the first sliver-cut function identified in \cref{tab:R2_B}. Additional dash-dotted lines are added to highlight the second, $(\bbeta_g^{-\frac{1}{2}} \hfrac^{\frac{1}{2}(3-p)})$-order, scaling.

\begin{figure}[H]
$$$$\\[-2cm]
    \centering
    \subfloat[$p=1$.]{\includegraphics[width=0.48\linewidth]{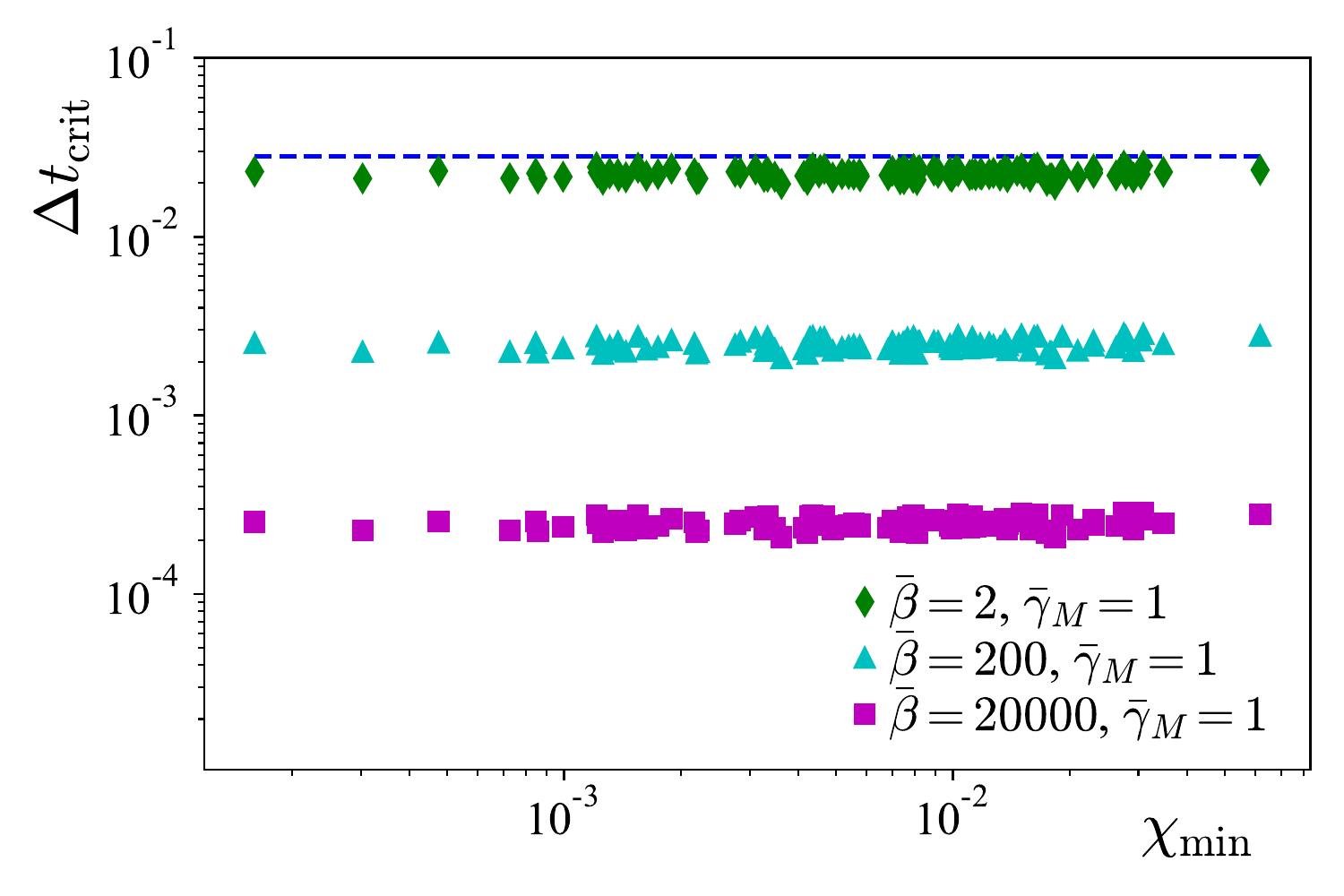}\label{fig:Dt2_penaltyMg_P1}}
    \subfloat[$p=2$.]{\includegraphics[width=0.48\linewidth]{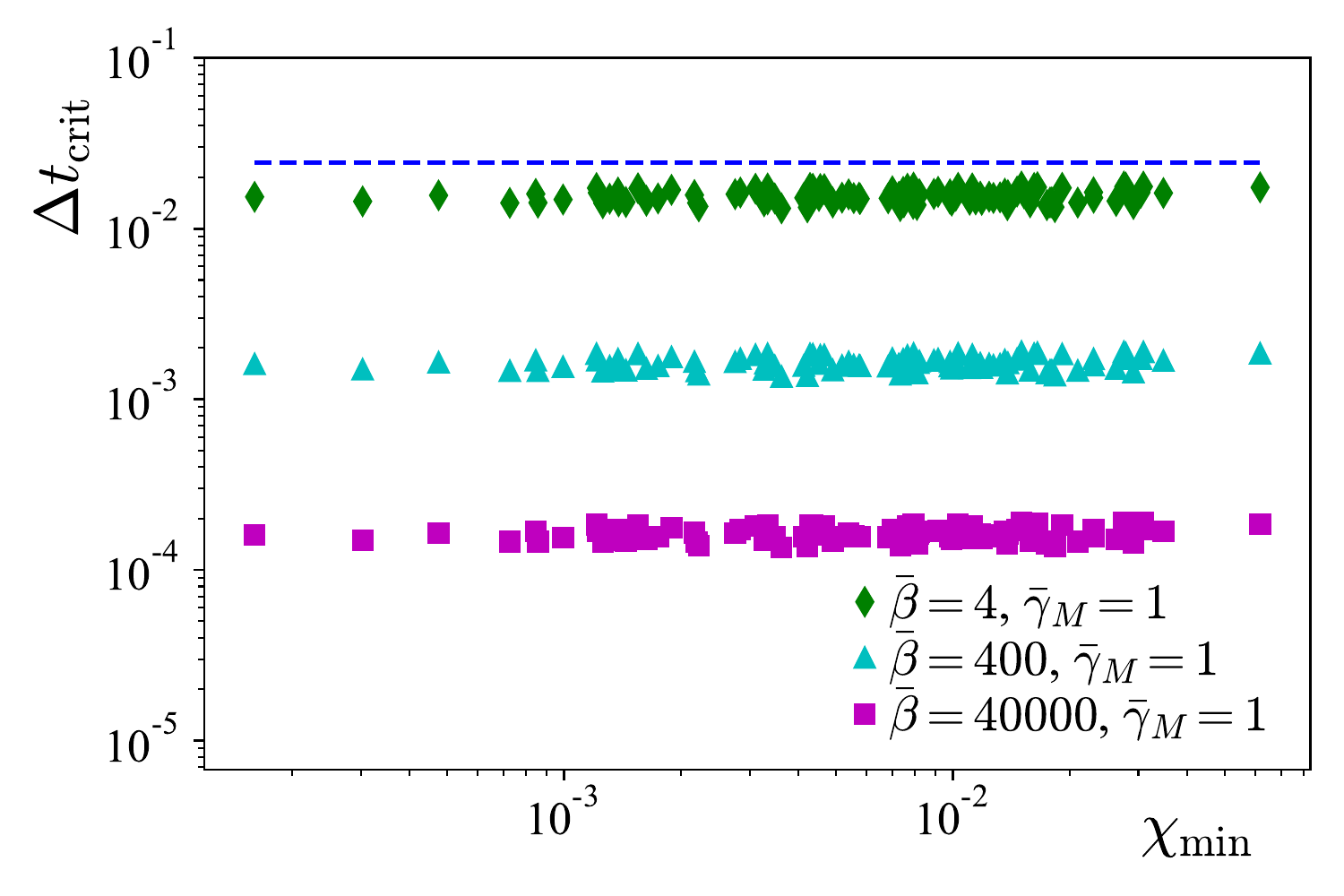}\label{fig:Dt2_penaltyMg_P2}}
    \caption{Penalty formulations with ghost mass for the second-order problem: second row and third column in \cref{tab:R1,tab:R1b,tab:R2,tab:R3,tab:Rconc}. Critical time-step size dependency on the minimal element size fraction for 100 perturbations of the domain of \cref{fig:2Dcase}.}
    \label{fig:Dt2_penaltyMg}
\end{figure}

\begin{figure}[H]
$$$$\\[-2cm]
    \centering
    \subfloat[$p=2$.]{\includegraphics[width=0.48\linewidth]{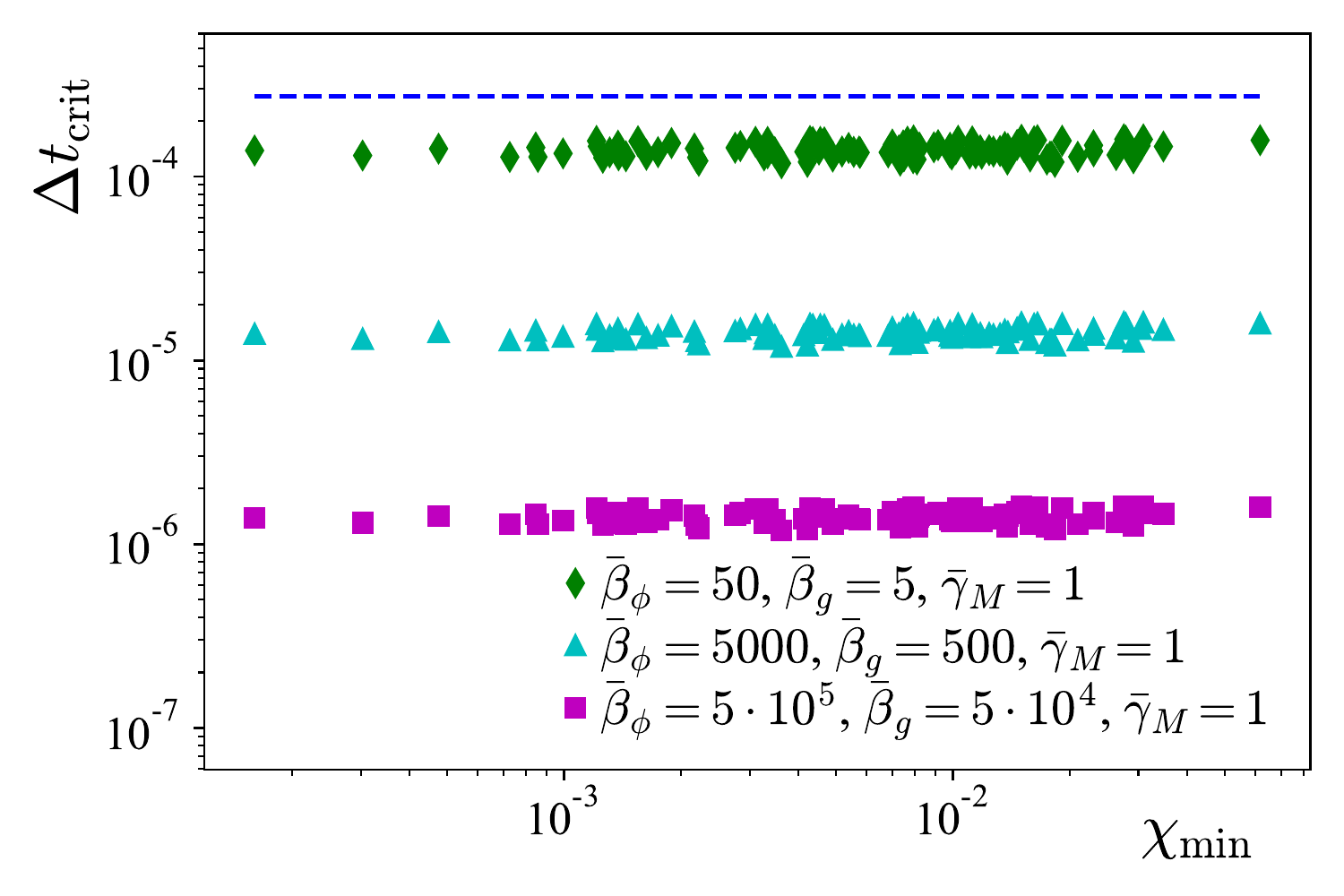}\label{fig:Dt4_penaltyMg_P2}}
    \subfloat[$p=3$.]{\includegraphics[width=0.48\linewidth]{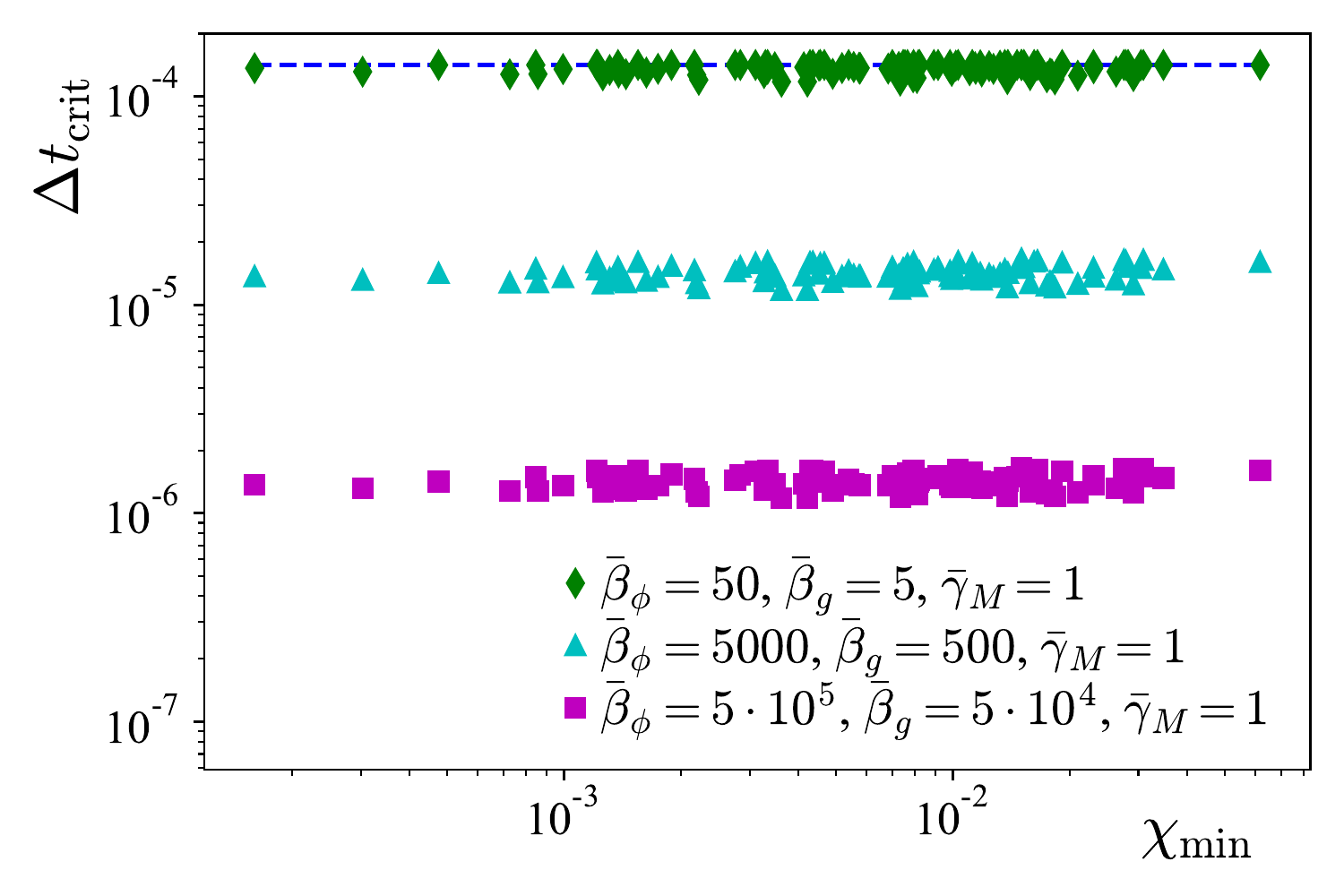}\label{fig:Dt4_penaltyMg_P3}}
    \caption{Penalty formulations for both $\phi$ and $\nabla\phi\cdot\normal$ with ghost mass for the fourth-order problem: second and fourth row and third column in \cref{tab:R1_B,tab:R1b_B,tab:R2_B,tab:R3_B,tab:Rconc_B}. Critical time-step size dependency on the minimal element size fraction for 100 perturbations of the domain of \cref{fig:2Dcase}.}
    \label{fig:Dt4_penaltyMg}
\end{figure}

Addition of the ghost-mass term only suppresses the $\hfrac^{\frac{1}{2}(s-p)}$-scaling, as it did for the pure Neumann formulation. This means that significant critical time-step size improvements can be achieved for small $\bbeta$ values, or when $p<s$ (from \cref{fig:Dt2_penalty_P1} to \cref{fig:Dt2_penaltyMg_P1}, from \cref{fig:Dt4_penaltyPhi_P2,fig:Dt4_penaltydPhi_P2} to \cref{fig:Dt4_penaltyMg_P2}, and from \cref{fig:Dt4_penaltyPhi_P3,fig:Dt4_penaltydPhi_P3} to \cref{fig:Dt4_penaltyMg_P3}), but not for larger penalty values, or when $p\geq s$.

\subsubsection{Nitsche formulations}

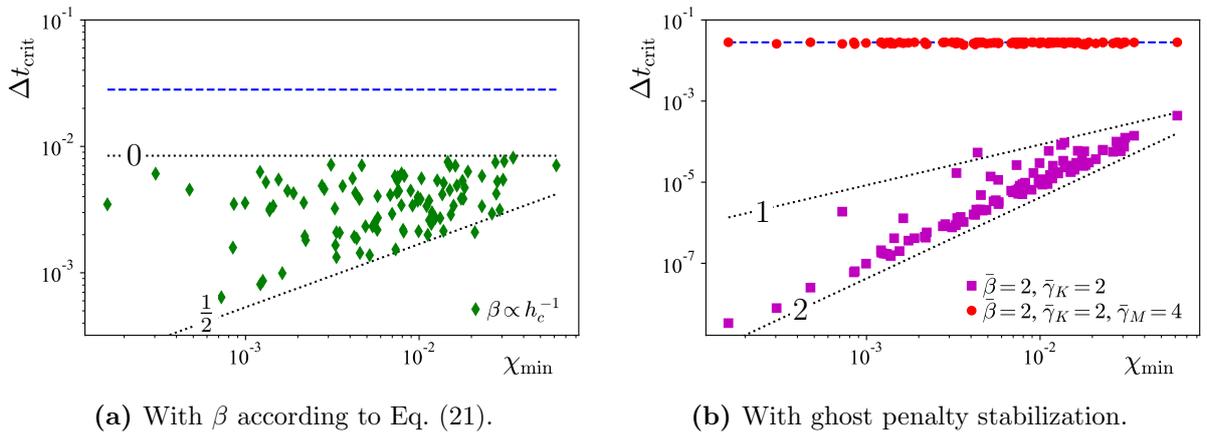
\begin{figure}[H]
\vspace{0.25cm}
    \centering
    \subfloat[With $\beta$ according to \cref{betac}.] {
        \begin{tikzpicture}
            \draw (0, 0) node[inner sep=0] {\includegraphics[width=0.48\linewidth]{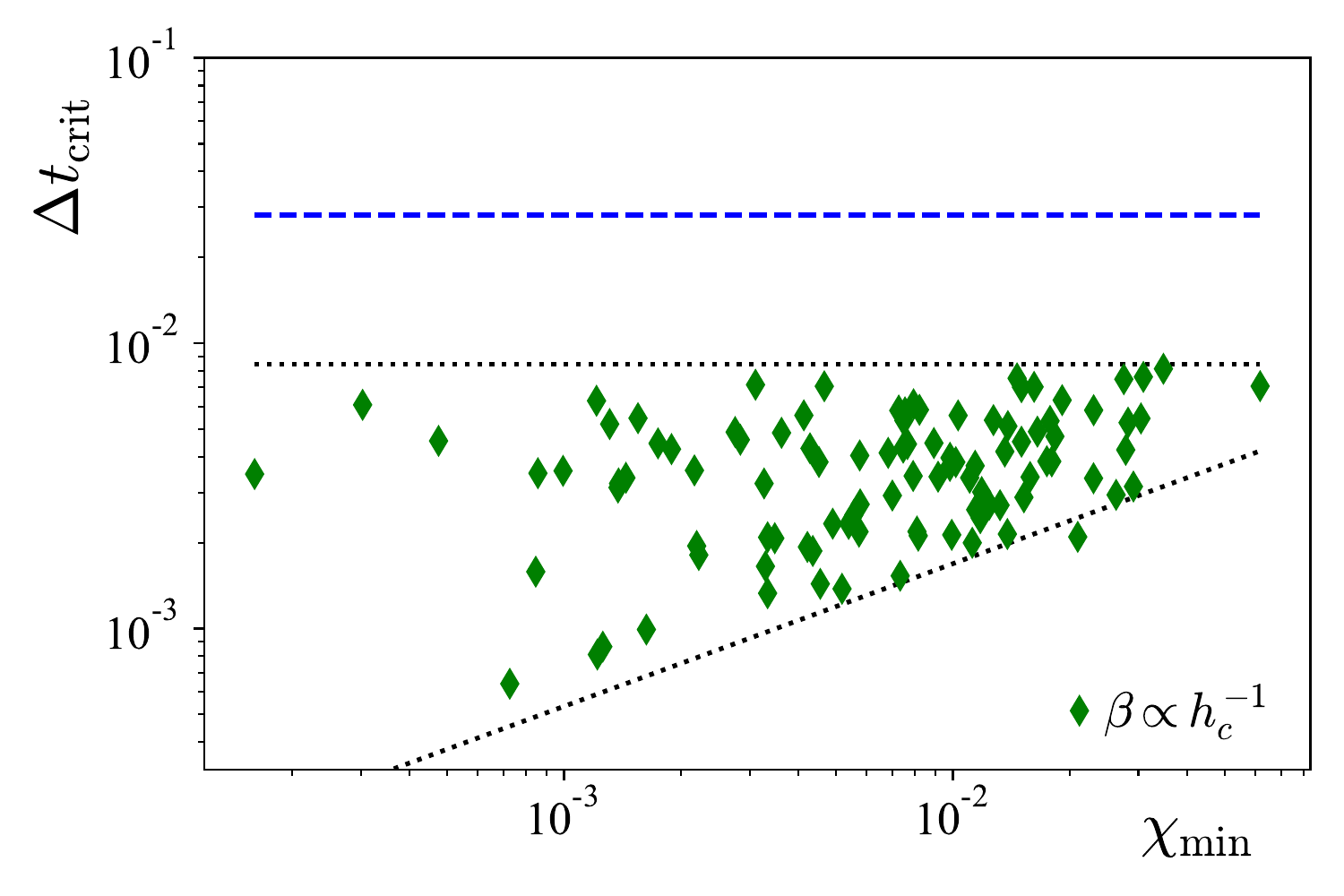}\label{fig:Dt2_Nitsche_a}};
            \draw (-1.15, -1.6) node [fill=white,inner sep=1pt] {$\frac{1}{2}$};
            \draw (-2.1, 0.5) node [fill=white,inner sep=1pt] {$0$};
        \end{tikzpicture} }
    \subfloat[With ghost penalty stabilization.]{
        \begin{tikzpicture}
            \draw (0, 0) node[inner sep=0] {\includegraphics[width=0.48\linewidth]{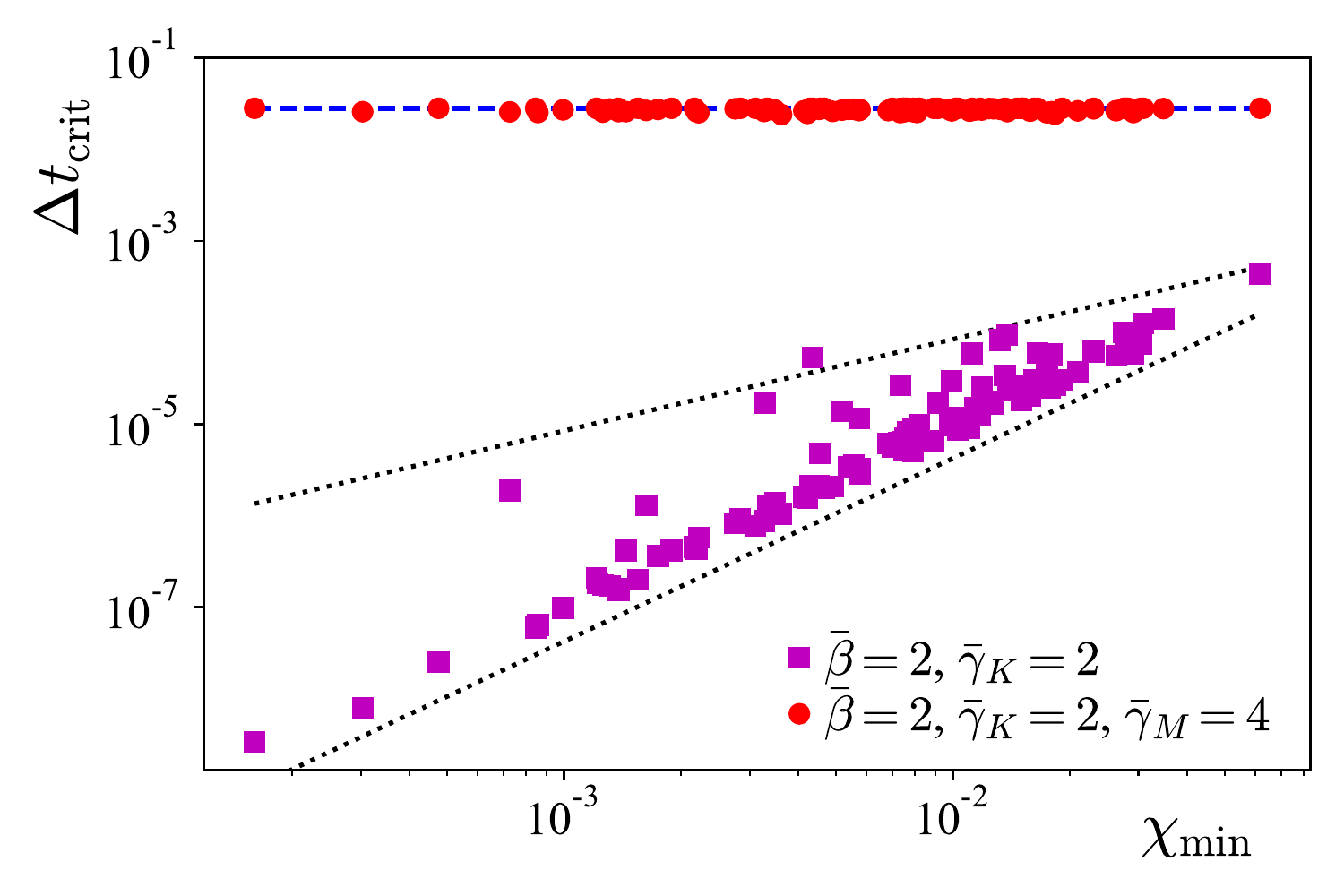}\label{fig:Dt2_Nitsche_b}};
            \draw (-2.0, -0.25) node [fill=white,inner sep=1pt] {$1$};
            \draw (-1.5, -1.55) node [fill=white,inner sep=1pt] {$2$};
        \end{tikzpicture}}
    \caption{Nitsche formulations for the second-order problem, for $p=1$: third and fourth row and second and third column in \cref{tab:R1,tab:R1b,tab:R2,tab:R3,tab:Rconc}. Critical time-step size dependency on the minimal element size fraction for 100 perturbations of the domain of \cref{fig:2Dcase}.}
    \label{fig:Dt2_Nitsche}
\end{figure}

\begin{figure}[H]
$$$$\\[-2cm]
    \centering
    \subfloat[On $\scalar$, with $\bbeta_\scalar$ according to \cref{betac1}.] {
        \begin{tikzpicture}
            \draw (0, 0) node[inner sep=0] {\includegraphics[width=0.48\linewidth]{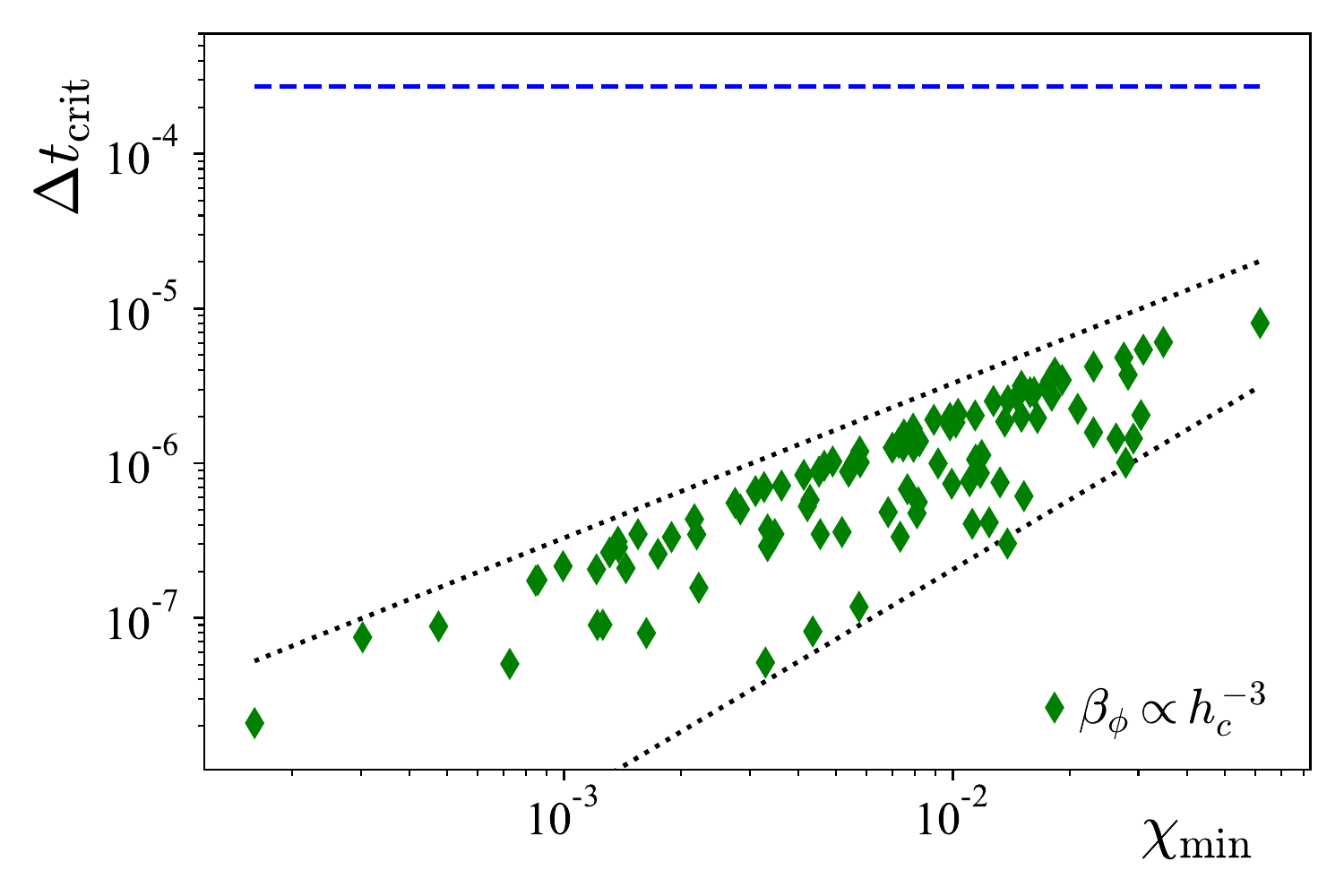}\label{fig:Dt4_Nitsche_a}};
            \draw (-2.2, -1.1) node [fill=white,inner sep=1pt] {$1$};
            \draw (0.2, -1.55) node [fill=white,inner sep=1pt] {$\frac{3}{2}$};
        \end{tikzpicture}}
    \subfloat[On $\scalar$, with ghost penalty stabilization.]{
        \begin{tikzpicture}
            \draw (0, 0) node[inner sep=0] {\includegraphics[width=0.48\linewidth]{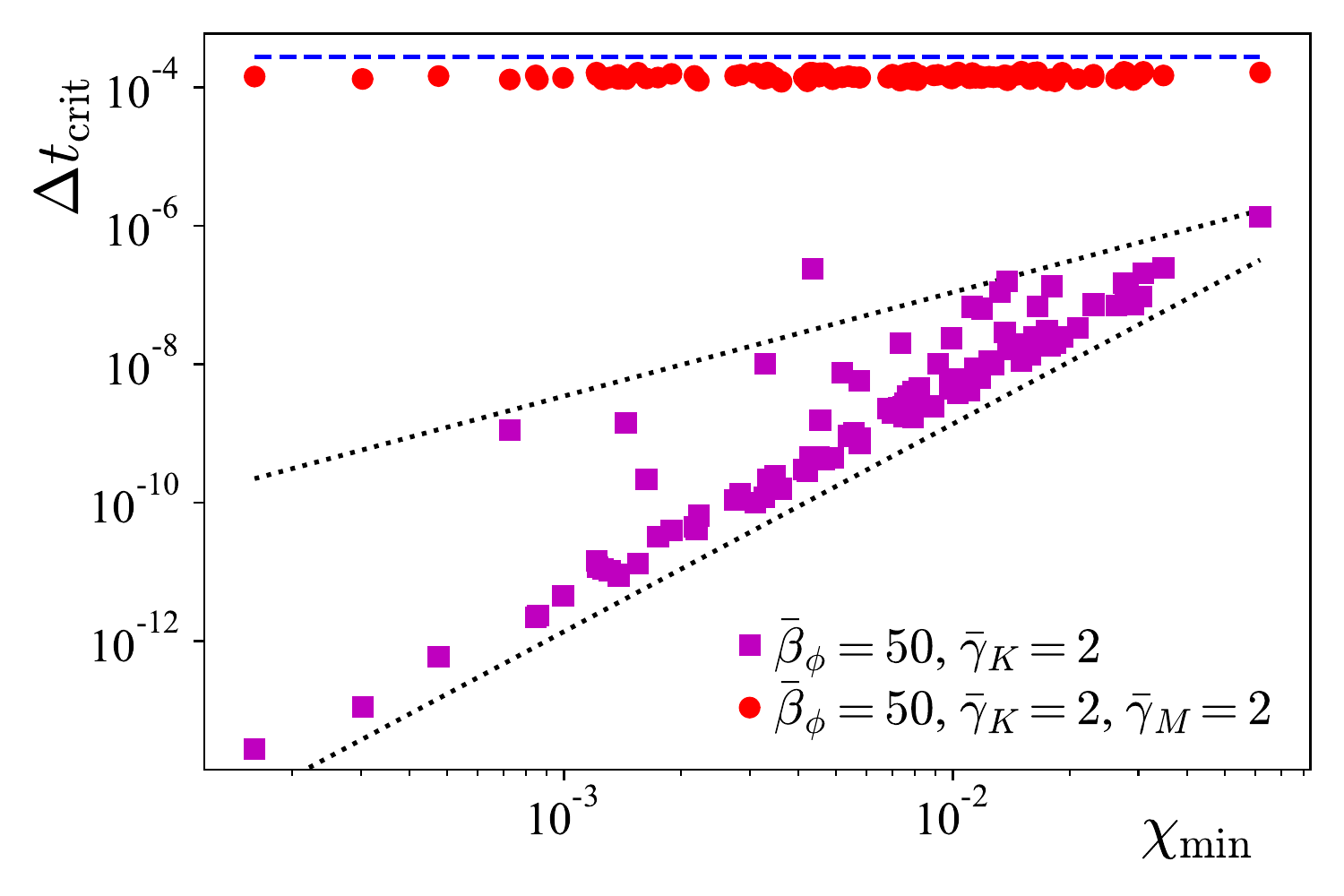}\label{fig:Dt4_Nitsche_b}};
            \draw (-1.9, -0.05) node [fill=white,inner sep=1pt] {$\frac{3}{2}$};
            \draw (-1.5, -1.55) node [fill=white,inner sep=1pt] {$3$};
        \end{tikzpicture}}\\
    \subfloat[On $\nabla\scalar\cdot\normal$, with $\bbeta_g$ according to \cref{betac2}.] {
        \begin{tikzpicture}
            \draw (0, 0) node[inner sep=0] {\includegraphics[width=0.48\linewidth]{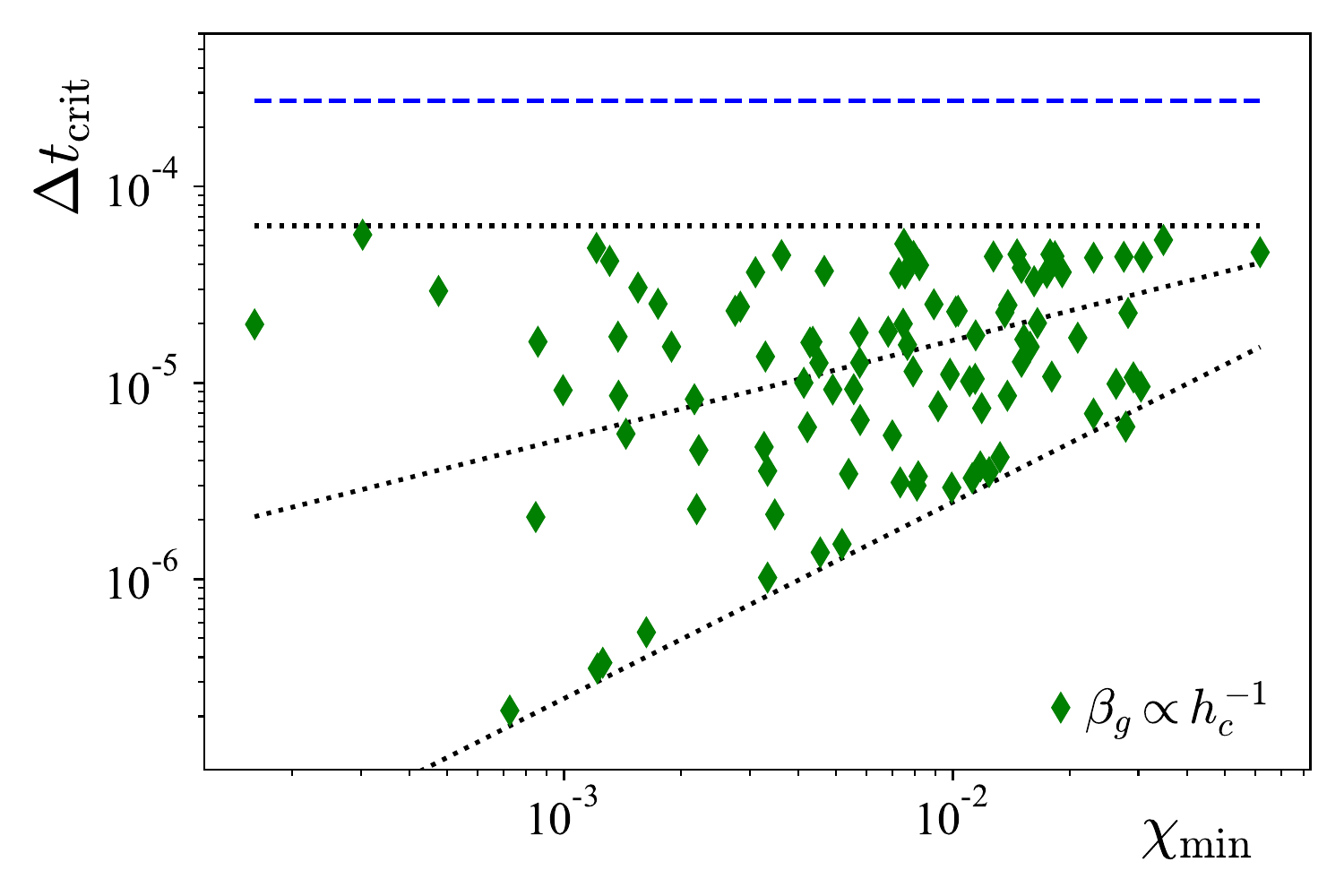}\label{fig:Dt4_Nitsche_c}};
            \draw (-2.2, 1.3) node [fill=white,inner sep=1pt] {$0$};
            \draw (-1.8, -0.25) node [fill=white,inner sep=1pt] {$\frac{1}{2}$};
            \draw (-1.2, -1.7) node [fill=white,inner sep=1pt] {$1$};
        \end{tikzpicture}}
    \subfloat[On $\nabla\scalar\cdot\normal$, with ghost penalty stabilization.]{
        \begin{tikzpicture}
            \draw (0, 0) node[inner sep=0] {\includegraphics[width=0.48\linewidth]{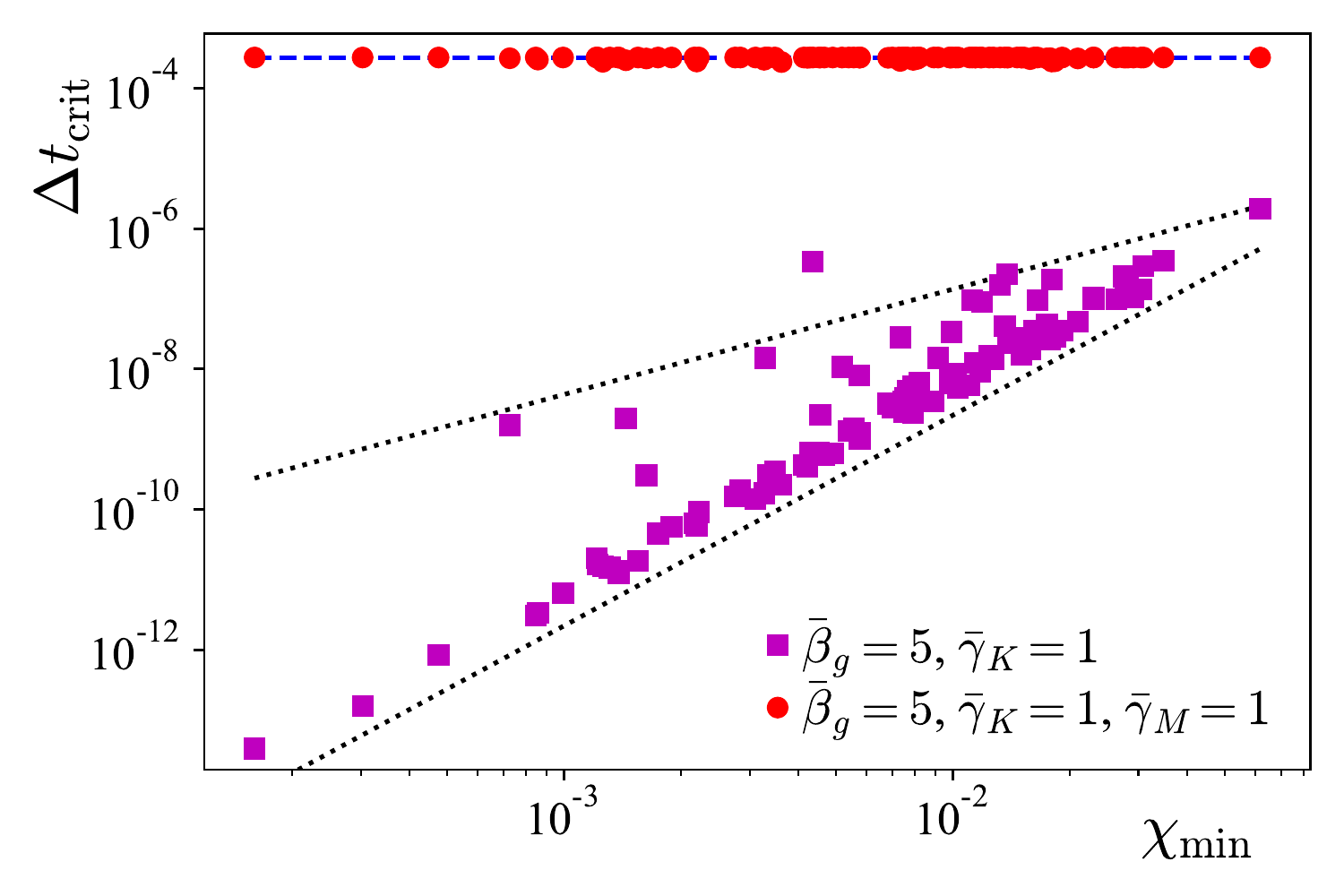}\label{fig:Dt4_Nitsche_d}};
            \draw (-1.9, -0.05) node [fill=white,inner sep=1pt] {$\frac{3}{2}$};
            \draw (-1.5, -1.55) node [fill=white,inner sep=1pt] {$3$};
        \end{tikzpicture}}
    \caption{Nitsche formulations for the fourth-order problem, for $p=2$: third, fourth and fifth row and second and third column in \cref{tab:R1_B,tab:R1b_B,tab:R2_B,tab:R3_B,tab:Rconc_B}. Critical time-step size dependency on the minimal element size fraction for 100 perturbations of the domain of \cref{fig:2Dcase}.}
    \label{fig:Dt4_Nitsche}
\end{figure}
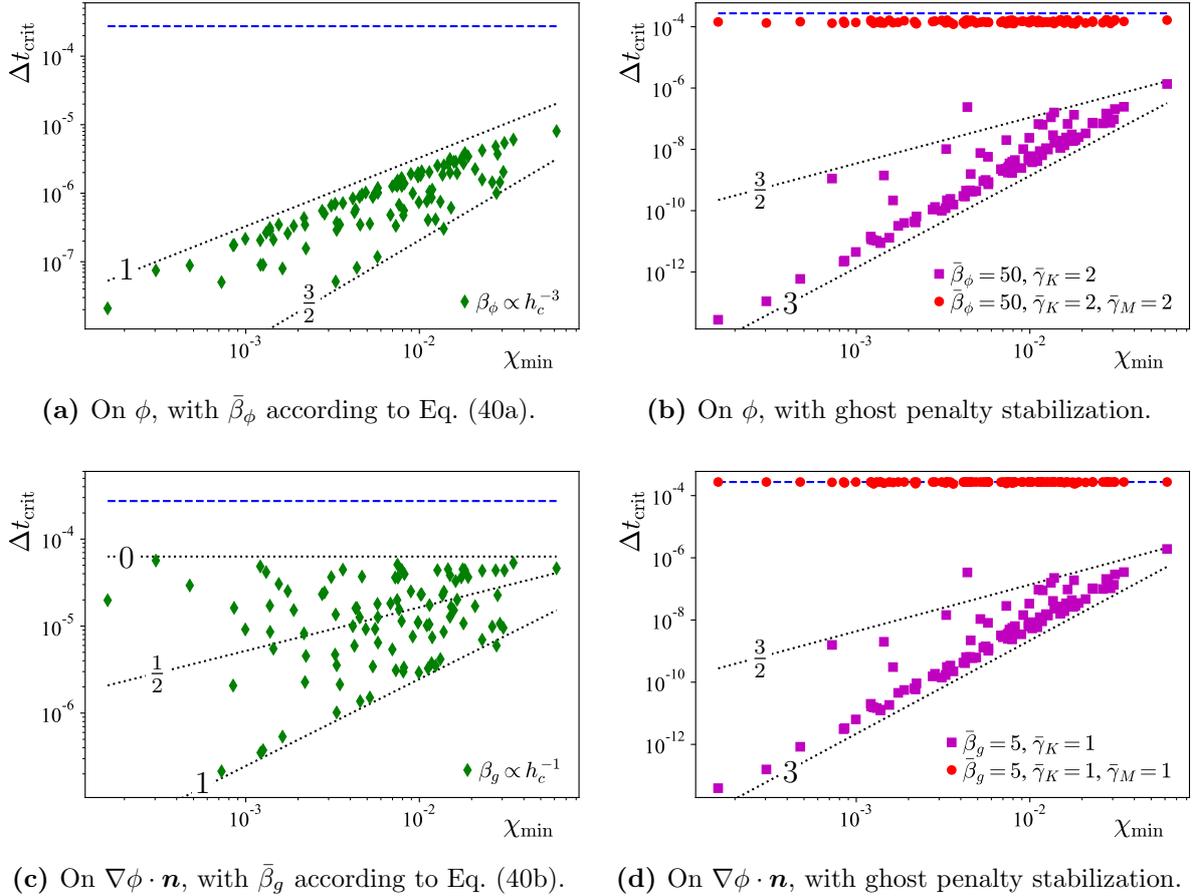

The results of the various Nitsche formulations are presented in \cref{fig:Dt2_Nitsche,fig:Dt4_Nitsche} for the second- and fourth-order problem, respectively. \Cref{fig:Dt2_Nitsche_a,fig:Dt4_Nitsche_a,fig:Dt4_Nitsche_c} correspond to Nitsche formulations with penalty parameters that are chosen per the inverse estimates of \cref{betac,betac1,betac2}, which thus involve the local element size measure $h_c = \hfrac\, h_\DomEl$. 
The four functions on cut elements induce different scaling orders. According to \cref{tab:R1,tab:R1_B}, the first corner-cut function causes a ($\hfrac^{\frac{1}{2}(s-pd)}$)-order scaling. For the considered combinations of $s$, $p$ and $d$, this reduces to a scaling with $\hfrac^0$ for all cases. This zeroth-order scaling is observed in \cref{fig:Dt2_Nitsche_a,fig:Dt4_Nitsche_c}, but not in \cref{fig:Dt4_Nitsche_a}. For the Nitsche formulation used in \cref{fig:Dt4_Nitsche_a}, the scaling induced by the \textit{second} corner-cut function dominates, as found in \cref{tab:R1b_B}. Furthermore, the first sliver-cut function induces a scaling of order $\hfrac^{\frac{1}{2}(s-p)}$ according to \cref{tab:R2,tab:R2_B}, and the second sliver-cut function induces a scaling of order $\hfrac^{\frac{1}{2}}$ (\cref{fig:Dt2_Nitsche_a,fig:Dt4_Nitsche_c}) or $\hfrac^{\frac{3}{2}}$ (\cref{fig:Dt4_Nitsche_a}) according to \cref{tab:R3,tab:R3_B}. These latter scalings are independent of $p$, and therefore increasing the polynomial order will not improve the results. All predicted scaling trends are plotted in each of the three subfigures (sometimes overlapping) confirming that they are indeed bounds of the obtained critical time-step size.

In \cref{fig:Dt2_Nitsche_b,fig:Dt4_Nitsche_b,fig:Dt4_Nitsche_d}, the results are plotted for the Nitsche formulations with penalty parameter that scale with the inverse of $h_\DomEl$. These formulations incorporate a ghost-penalty term in the stiffness matrix to guarantee a its positive definiteness. According to \cref{tab:R1,tab:R1_B}, a function on a corner-cut element induces a detrimental cut-size dependency of the critical time-step size of order $\hfrac^{\frac{1}{2}(pd+d)}$. Similarly, the first sliver-cut function causes a scaling of the order $\hfrac^{\frac{1}{2}(p+1)}$ as stated in \cref{tab:R2,tab:R2_B}. Both these trends are plotted in each of the respective subfigures. The positive exponent $p$ in these relations implies that this scaling will only worsen with increasing polynomial order. However, the addition of the ghost-mass term eliminates the cut-size dependency altogether. The highest eigenvalue does still scale with $\bbeta \hfrac^0$. In \cref{fig:Dt4_Nitsche_b} we observe that the red markers exhibit a smaller critical time-step size than the maximally achievable value indicated by the dashed-blue line. This indicates that the penalty value required to stabilized the Nitsche formulation is not sufficiently small to avoid affecting the highest eigenvalues.

\subsection{Convergence of a linear pre-stressed membrane}
\label{ssec:drum}

The results of the analysis on the scaling of the critical time-step size with the cut-element size indicate that addition of ghost mass can, in certain cases, significantly increase the critical time-step size. In particular, it enables a Nitsche formulation with a cut-size independent critical time-step size. Of course, the added ghost mass should not come at the cost of a severe accuracy reduction. In this section, we study the impact of the added ghost mass on the solution error for a linear pre-stressed membrane, i.e., for the second-order wave equation of \cref{sec:SecondOrder}. We consider the same geometry description as before, depicted in \cref{fig:2Dcase}, and take as an exact solution a simple standing sine wave:
\begin{align}
    \scalar_{\text{exact}}(t,x,y) = \cos( \sqrt{2} \pi t ) \sin( \pi x ) \sin(\pi y) \,,
\end{align}
from which we infer the required initial and boundary conditions. In the following, we compute one full period of this oscillation, for which we use a Newmark-type central difference method for time integration. This limits the optimal convergence rates to second-order.

\begin{figure}[!b]
\vspace{0.25cm}
    \centering
    \subfloat[$H^1_0$-error.]{\includegraphics[width=0.48\linewidth]{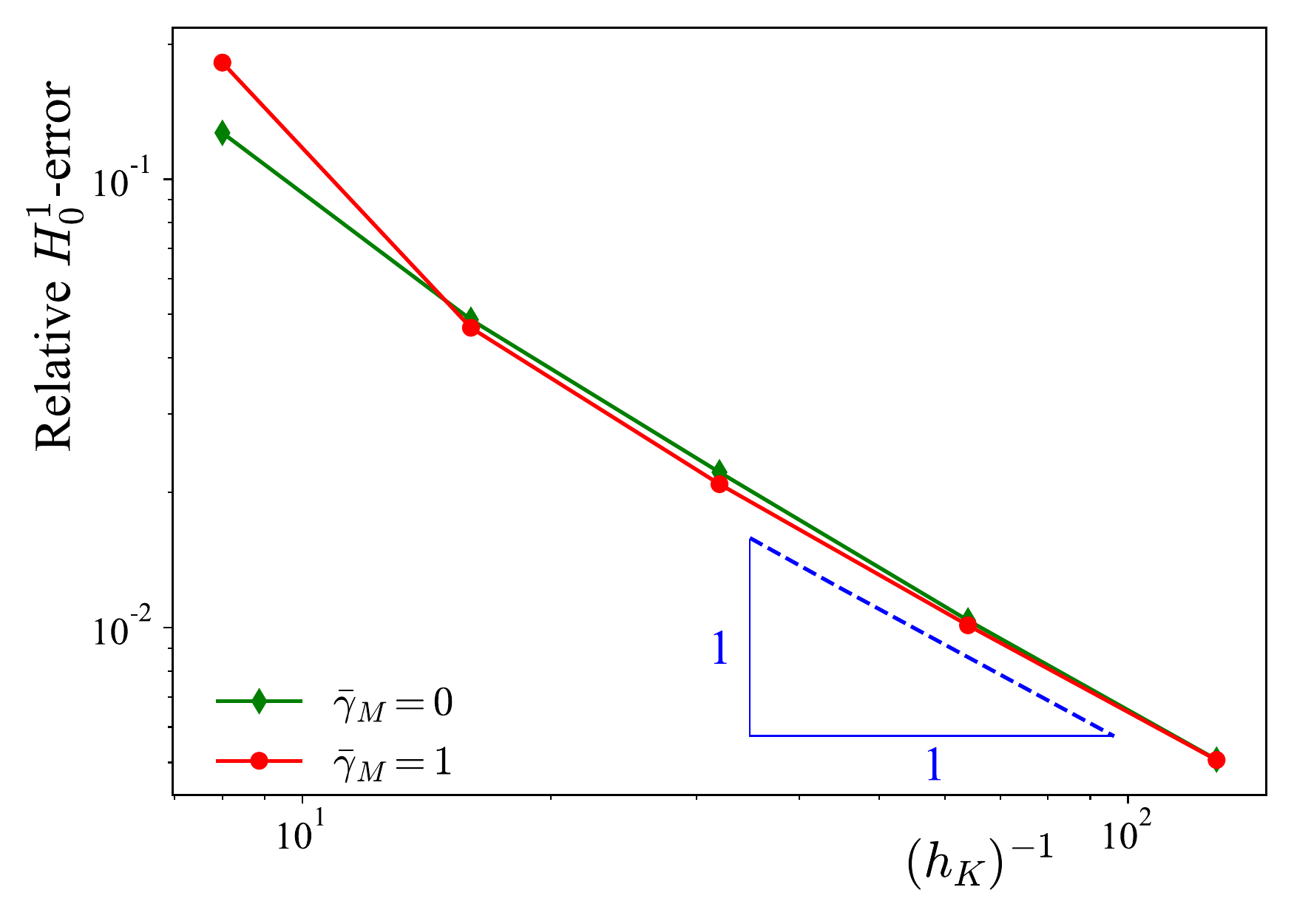}\label{fig:dconvP1Neum_a}}
    \subfloat[$L2$-error.]{\includegraphics[width=0.48\linewidth]{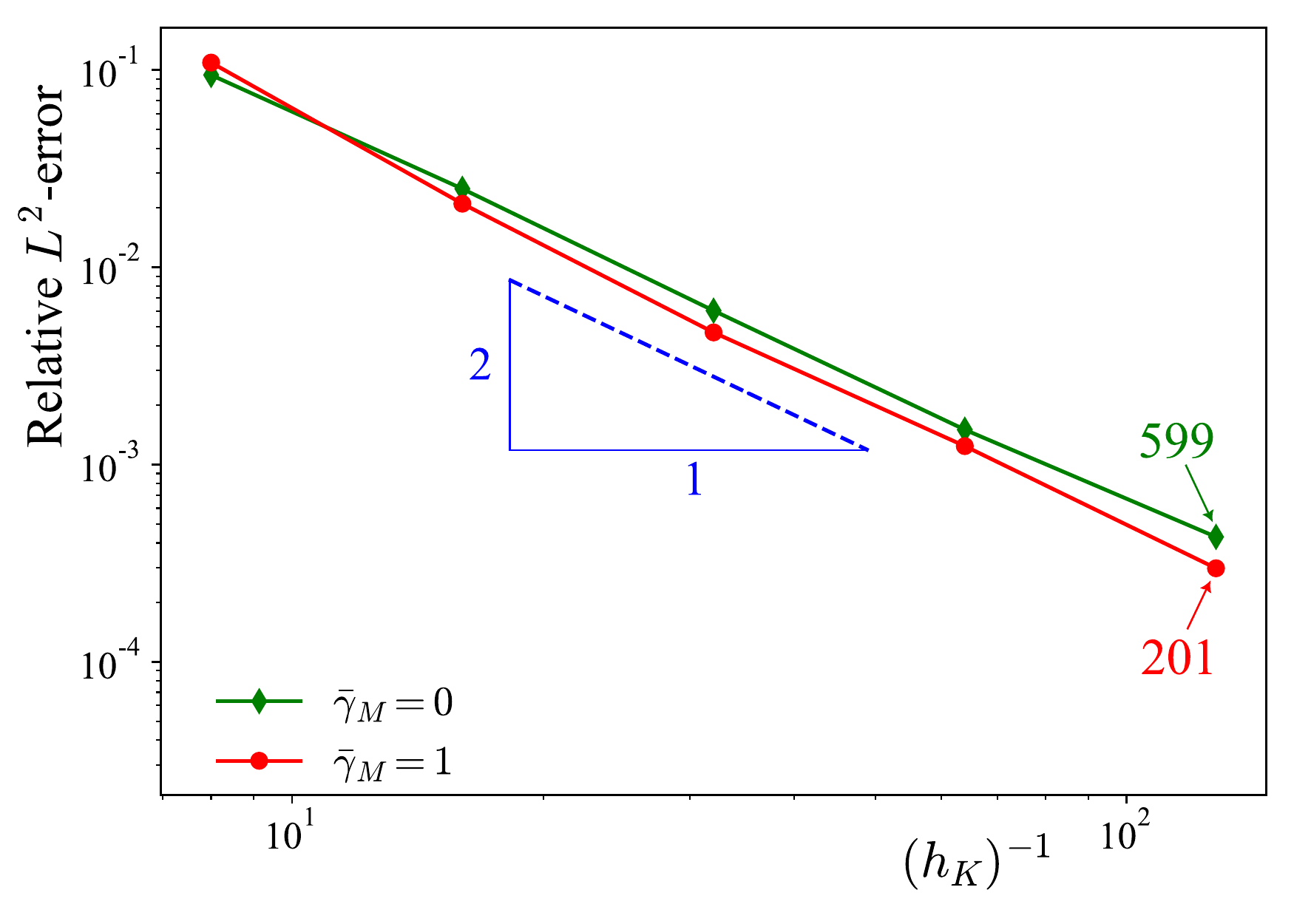}\label{fig:convP1Neum_b}}
    \caption{Error convergence for a vibrating pre-stressed membrane on the domain of \cref{fig:2Dcase} after one full period of oscillation, for the Neumann case with and without ghost mass, with $p=1$.}
    \label{fig:convP1Neum}
\end{figure}

First, we investigate the scenario where the entire cut-out represents a Neumann boundary. According to \cref{tab:Rconc}, and as verified in \cref{fig:Dt2_Neumann_P1,fig:Dt2_Neumann_P1}, the addition of ghost mass only affects the critical time-step size for the $p=1$ order of basis functions. \Cref{fig:convP1Neum} shows the convergence curves of the relative $H^{1}_{0}$ and $L^2$-errors for the formulations with and without a ghost-mass term. The dashed-blue lines indicate the optimal orders of convergence. The numbers attached to the final markers in the $L^2$-error graphs denote the total number of time steps required to carry out the corresponding simulation. For this particular example, we observe a factor 3 reduction of the required number of time steps on the most refined grid. Despite this reduction, we observe that the addition of the ghost-mass term does not compromise the accuracy of the approximation.

\enlargethispage{-0.5cm}

\begin{figure}[!t]
\vspace{-0.5cm}
    \centering
    \subfloat[$H^1_0$-error, $p=1$.]{\includegraphics[width=0.48\linewidth]{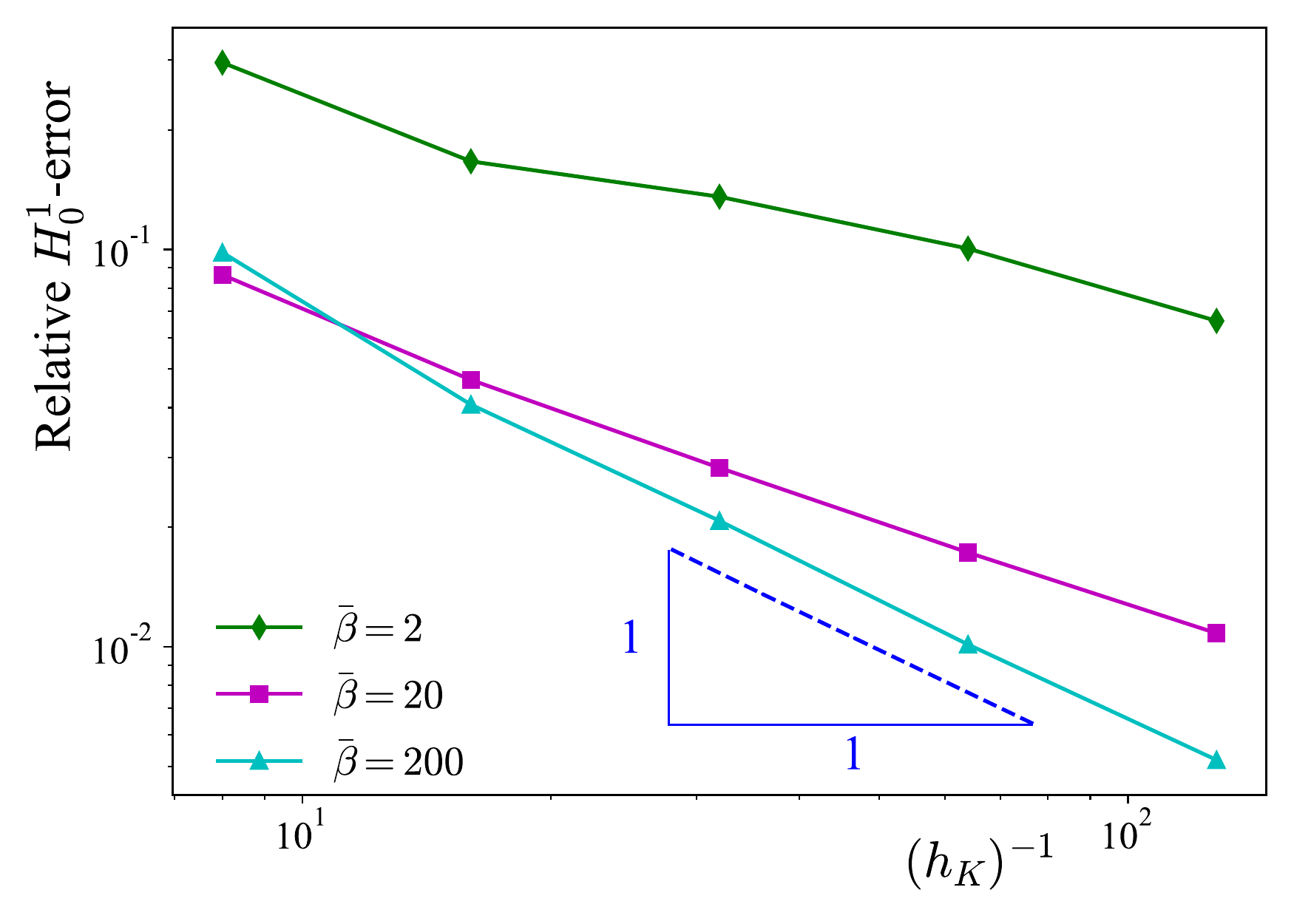}\label{fig:dconvP1Penalty_a}}
    \subfloat[$L2$-error, $p=1$.]{\includegraphics[width=0.48\linewidth]{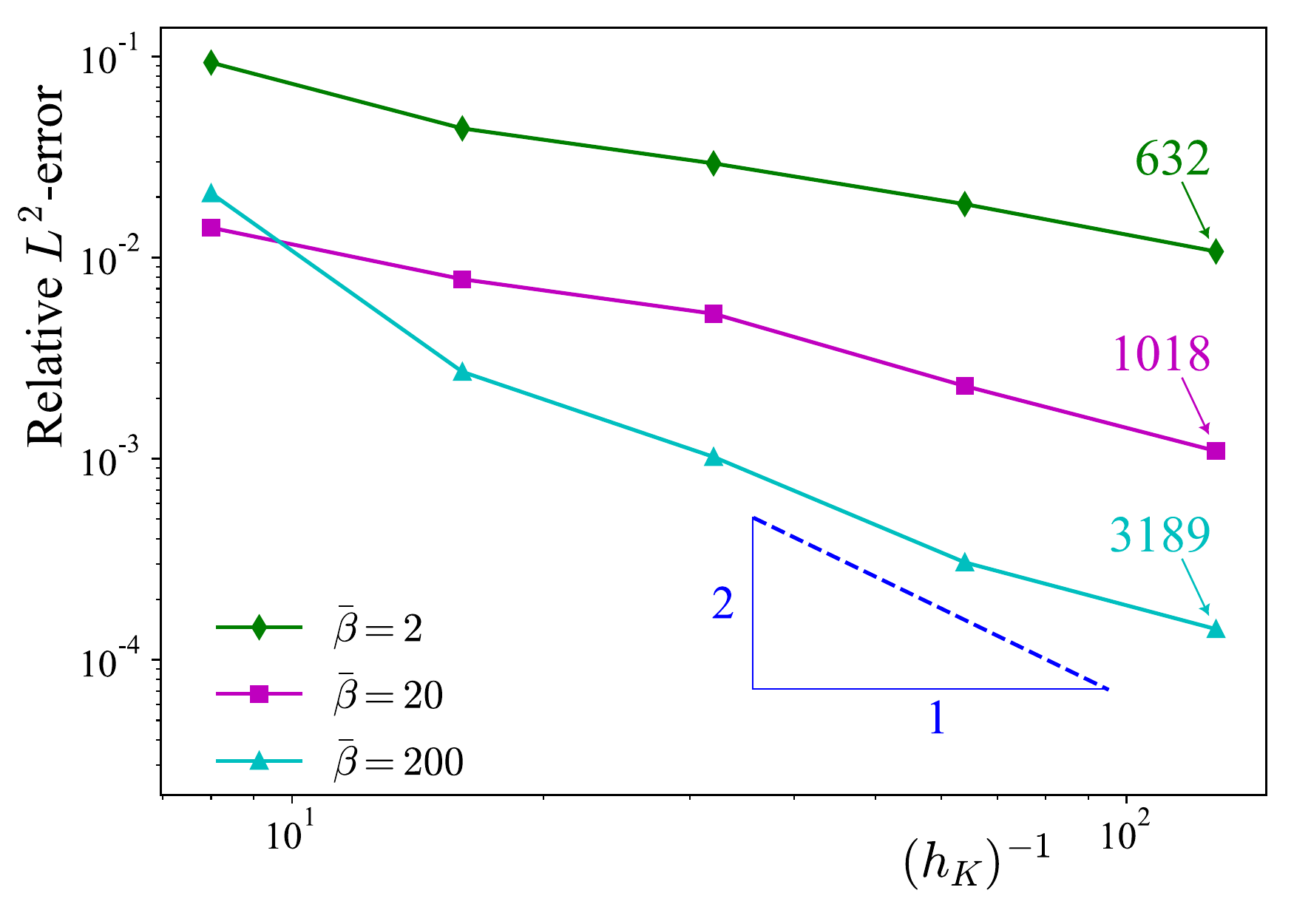}\label{fig:convP1Penalty_b}}\\
    \subfloat[$H^1_0$-error, $p=2$.]{\includegraphics[width=0.48\linewidth]{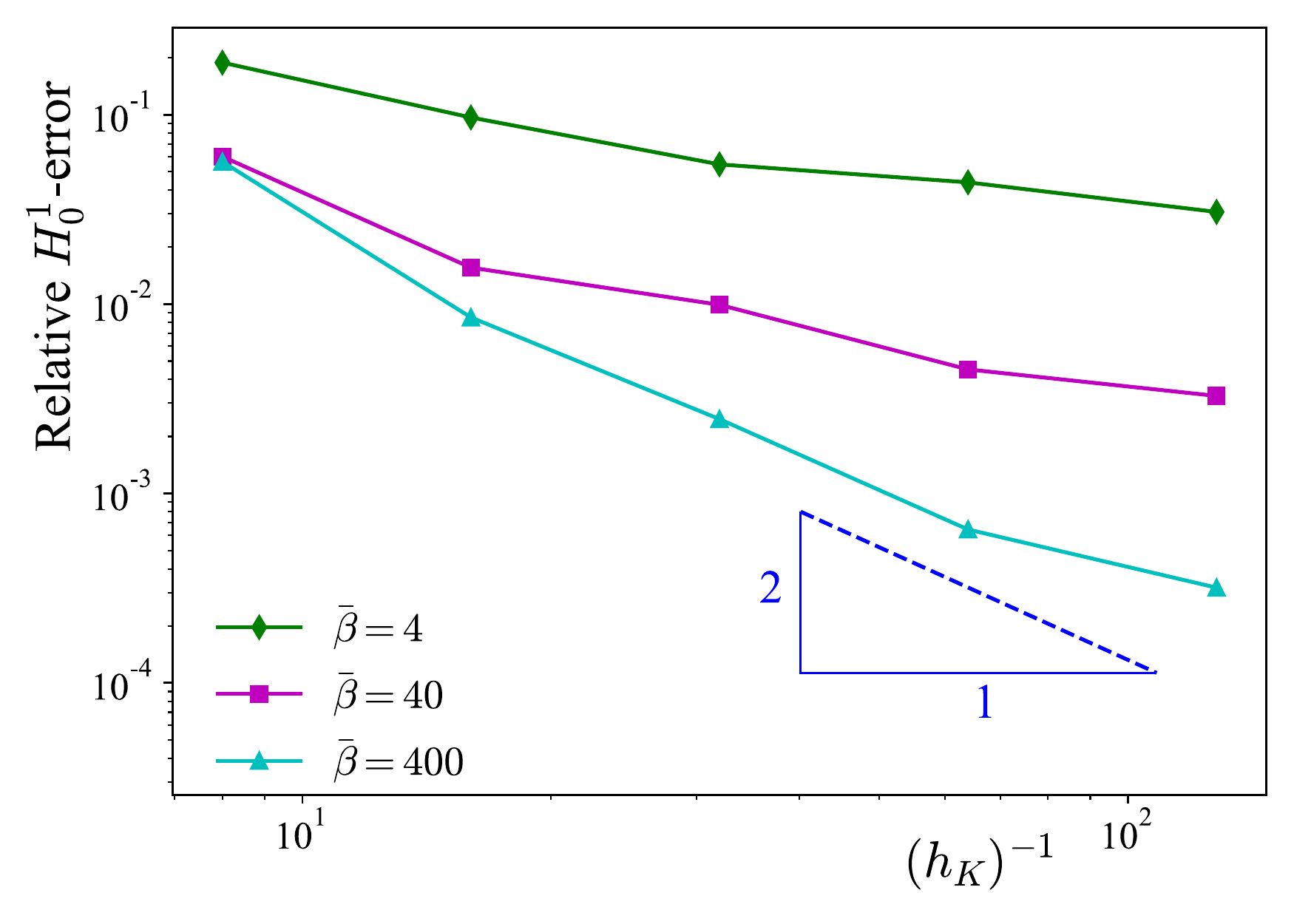}\label{fig:dconvP2Penalty_a}}
    \subfloat[$L2$-error, $p=2$.]{\includegraphics[width=0.48\linewidth]{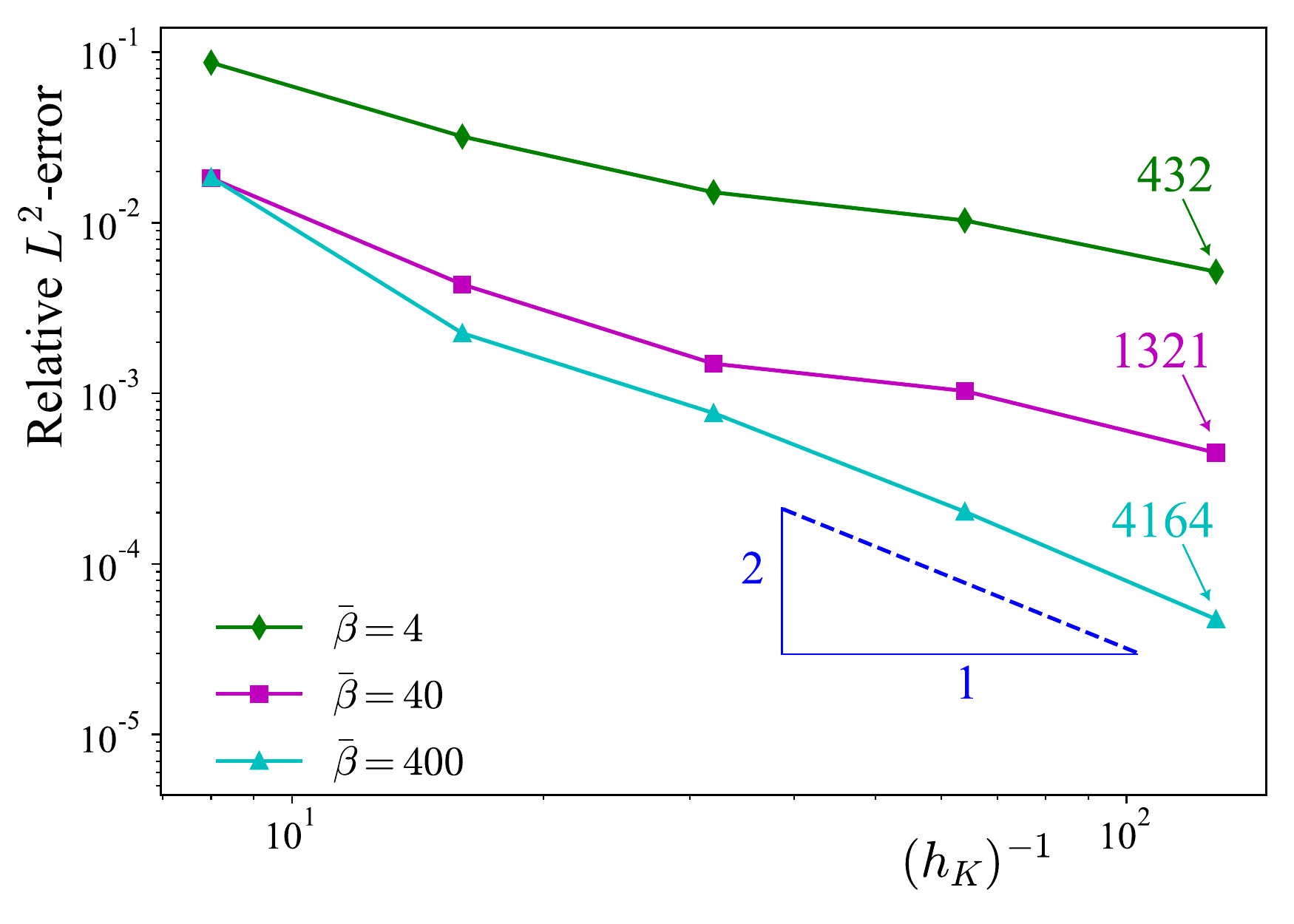}\label{fig:convP2Penalty_b}}
    \caption{Error convergence for a vibrating pre-stressed membrane on the domain of \cref{fig:2Dcase} after one full period of oscillation, for penalty methods with different penalty factors.\\}
    \label{fig:convPenalty}
\end{figure}

Next, we consider penalty methods with different penalty factors, without ghost mass. 
The error convergence curves for $p=1$ and $p=2$ are shown in \cref{fig:convPenalty}. For both polynomial orders, the smallest chosen penalty parameter only marginally impacts the critical time-step size, as may be ascertained from \cref{fig:Dt2_penalty}. As a result, the $p=1$ simulation on the most refined grid requires $\sim 5\%$ more time steps than the corresponding Neumann case plotted in \cref{fig:convP1Neum}. At the same time, the error has increased significantly. This is in part due to the variational inconsistency of the penalty formulation, which causes a loss of the optimal (spatial) convergence rate.
The errors may be reduced by several orders of magnitude by increasing the penalty parameter. However, an increase of the penalty parameter comes at the cost of a reduced critical time-step size. In particular, we observe a factor 5 (for $p=1$) and 10 (for $p=2$) increase of the required total number of time steps for a factor 100 increase in $\bbeta$. At the highest penalty level, the optimal orders of convergence appear to be achieved in some of the subfigures, but these convergence rates drop past a certain mesh refinement level.

\begin{figure}[!t]
    \centering
    \subfloat[$H^1_0$-error, $p=1$.]{\includegraphics[width=0.48\linewidth]{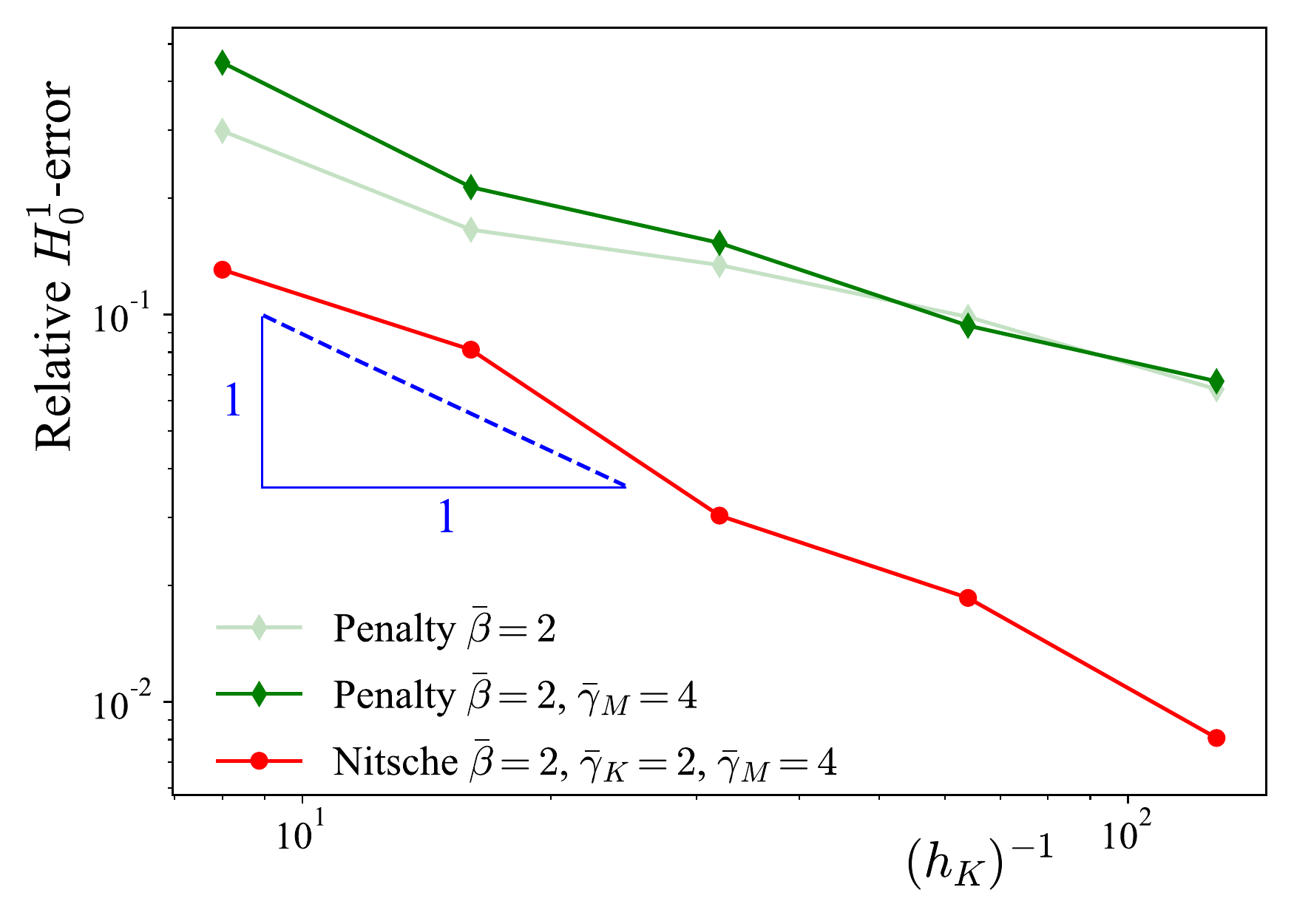}\label{fig:dconvP1Nitsche_a}}
    \subfloat[$L2$-error, $p=1$.]{\includegraphics[width=0.48\linewidth]{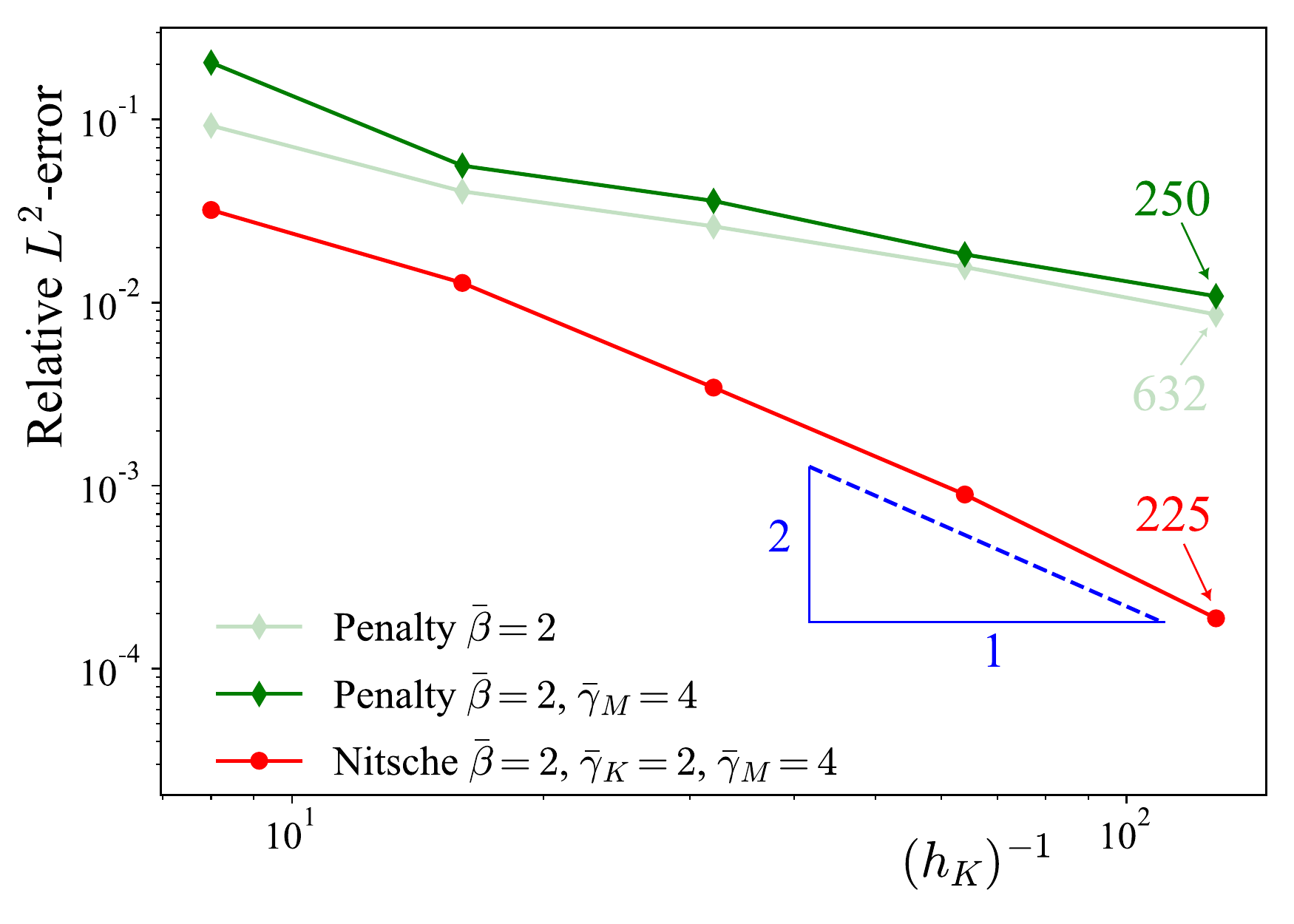}\label{fig:convP1Nitsche_b}}\\
    \subfloat[$H^1_0$-error, $p=2$.]{\includegraphics[width=0.48\linewidth]{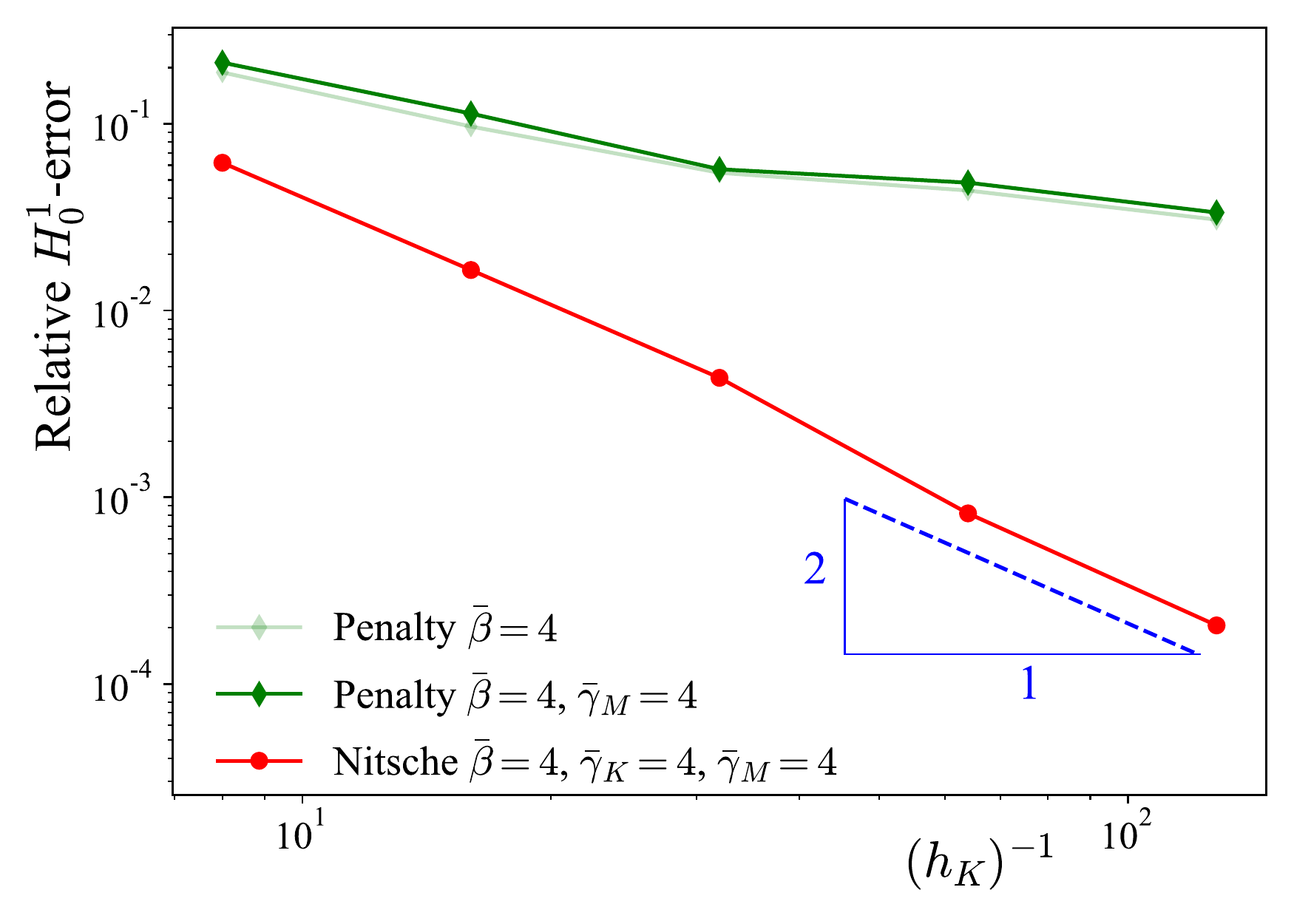}\label{fig:dconvP2Nitsche_a}}
    \subfloat[$L2$-error, $p=2$.]{\includegraphics[width=0.48\linewidth]{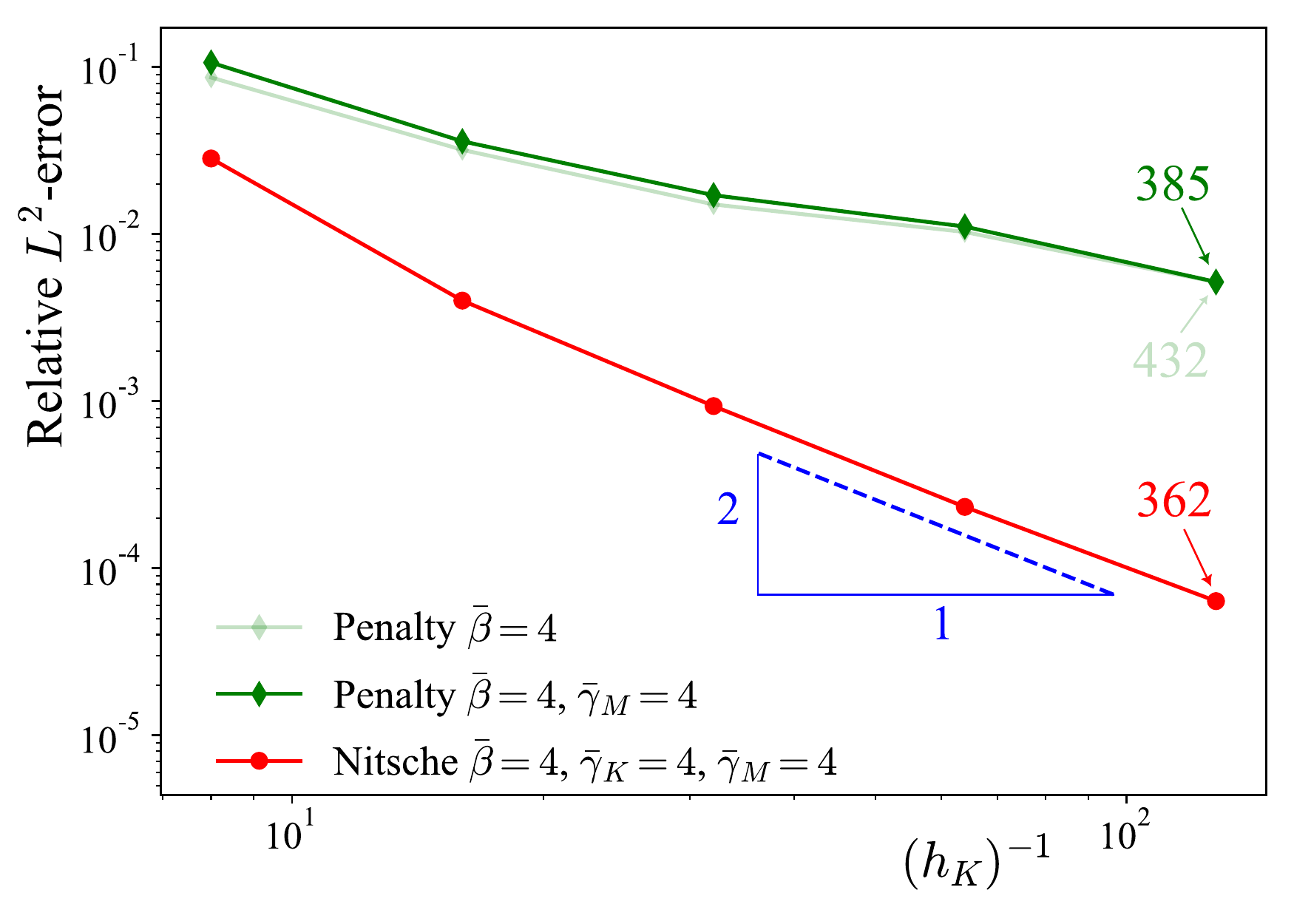}\label{fig:convP2Nitsche_b}}
    \caption{Error convergence for a vibrating pre-stressed membrane on the domain of \cref{fig:2Dcase} after one full period of oscillation, for Nitsche or penalty enforcement of Dirichlet constraints.}
    \label{fig:convNitsche}
\end{figure}

Finally, we compare the performance of a penalty method with that of a Nitsche formulation. 
To ensure an equitable comparison, both formulations are augmented with ghost mass and make use of the same penalty parameter values. As a result, the simulation with Nitsche's method and the simulation with the penalty method require a comparable number of times steps. The relative $H^{1}_{0}$ and $L^2$-errors are plotted in \cref{fig:convNitsche}, for $p=1$ and $p=2$. For reference, the results of the penalty methods without ghost mass from \cref{fig:convPenalty} are overlaid. As anticipated, the addition of the ghost mass to the penalty method for these moderate penalty parameter values reduces the required number of time steps significantly for $p=1$ and only little for $p=2$. In both cases, the ghost mass only marginally affects accuracy. A jump in solution quality is achieved by switching to the Nitsche formulation. Due to the variational consistency of the formulation, both the $H^1_0$ and $L^2$-errors converge optimally, yielding error reductions by orders of magnitude on the finer meshes. For $p=1$, the error values even closely resemble those of the pure Neumann case in \cref{fig:convP1Neum}.

\subsection{Transient response of a linear Kirchhoff-Love shell}
\label{ssec:shell}

\begin{figure}[t]
\vspace{-0.5cm}
    \centering
    \subfloat[Problem description.]{\includegraphics[width=0.38\linewidth]{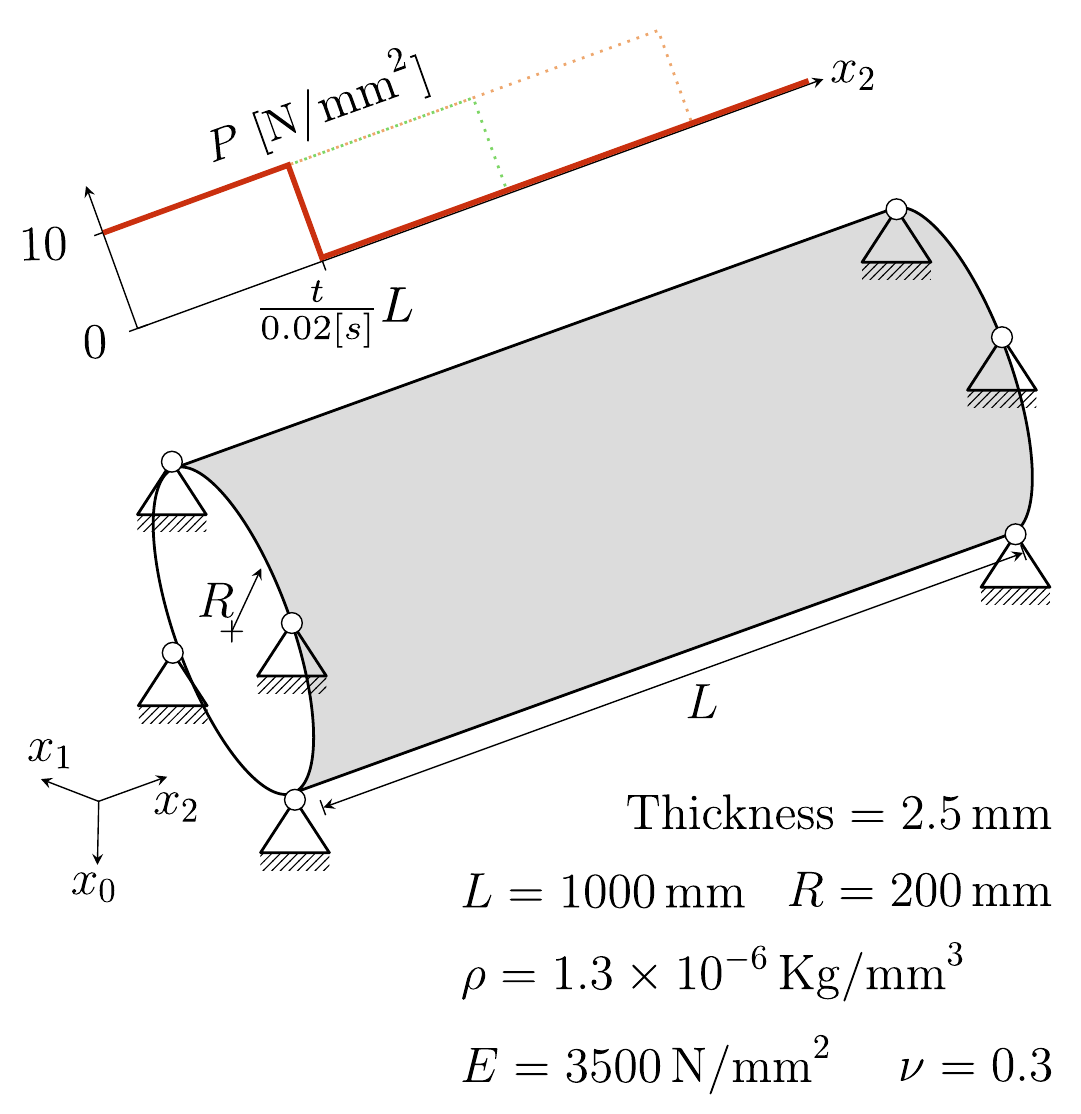}\label{fig:KirchhoffLoveModel_a}}\hfill
    \subfloat[Transient solution.]{\includegraphics[width=0.58\linewidth,trim=90 90 90 45, clip]{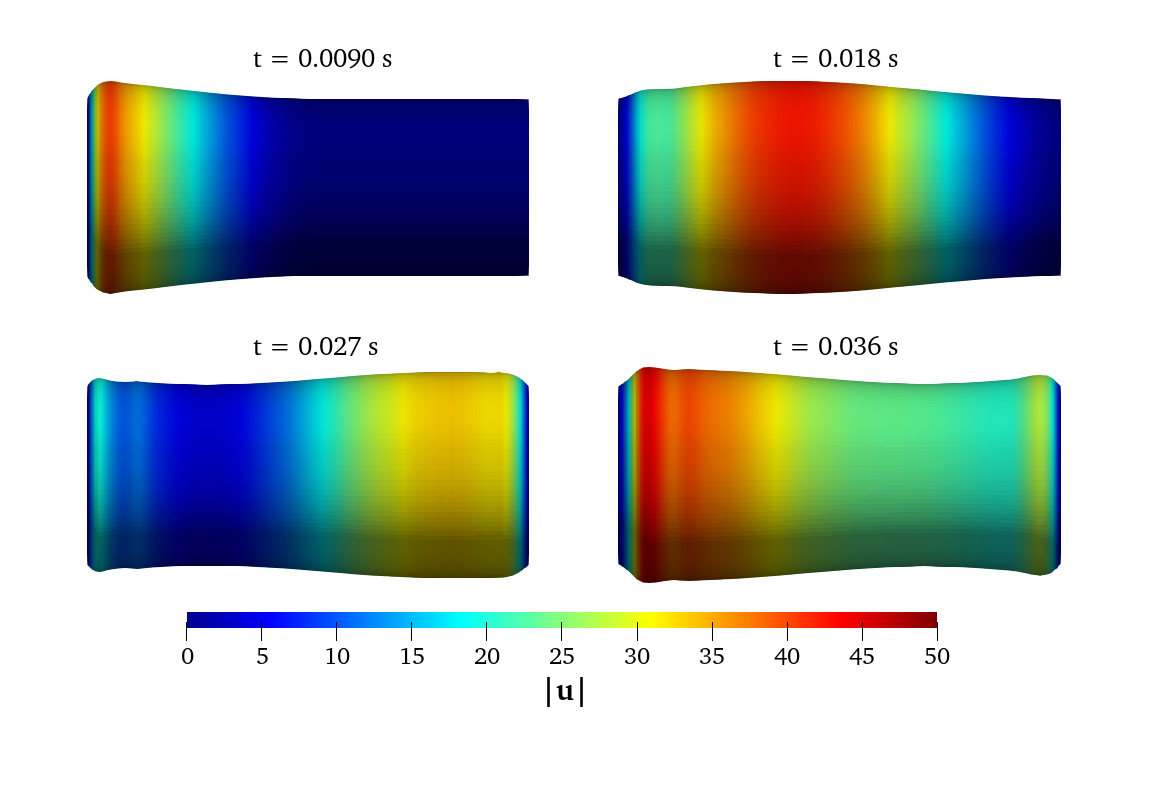}\label{fig:KirchhoffLoveModel_b}}
    \caption{Kirchhoff-Love shell model problem.}
\end{figure}

To demonstrate the practical applications of ghost mass and Nitsche's method, we examine their usage in a transient simulation of a Kirchhoff-Love shell  \cite{Bischoff2004,Kiendl2009}. As a case study, we consider a pressure shock-wave propagating through a pipe-segment with pinned support on both ends. \Cref{fig:KirchhoffLoveModel_a} shows the geometry, the material parameters and the load-function. The shock-wave travels through the entire pipe-segment in $0.02$ seconds. \Cref{fig:KirchhoffLoveModel_b} shows the resulting displacement field of a reference computation, computed on a $30\times 30$ grid of $p=2$ polynomial B-splines. 

\begin{figure}[t]
    \centering
    \includegraphics[width=0.7\linewidth]{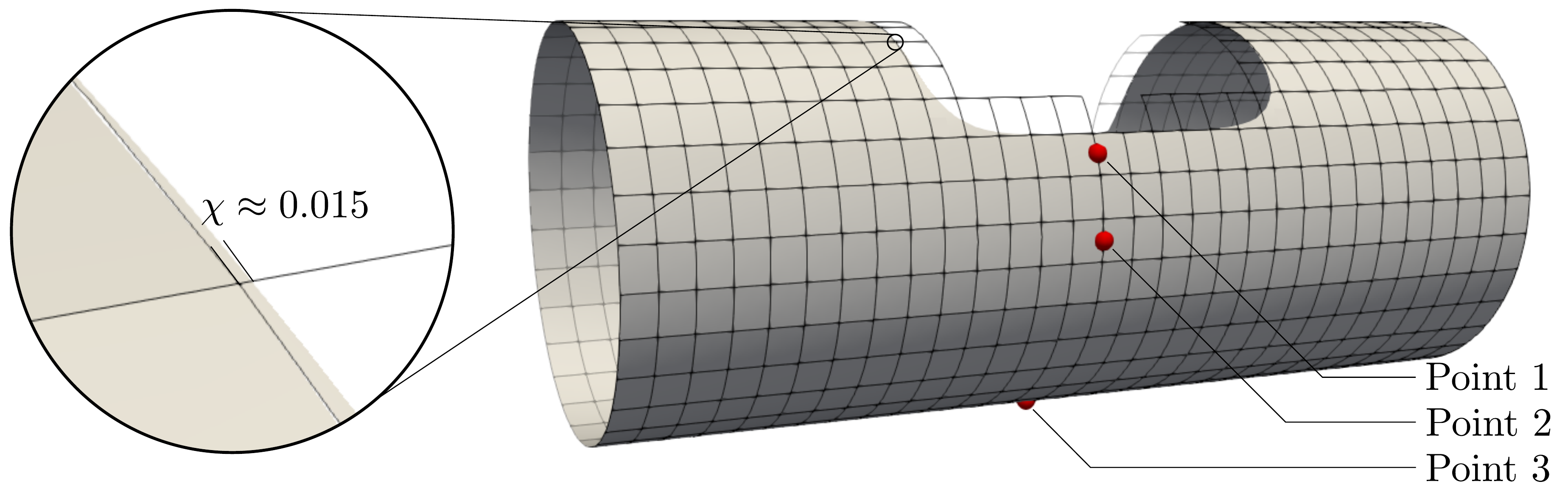}
    \caption{Immersed domain and mesh of $30\times 30$ elements, with filleted rectangular cut-out. Indicated points are measurement locations referenced in \cref{fig:KL_pointvals}.}
    \label{fig:KirchhoffLoveModelImmersed}
\end{figure}

We then create a rectangular cut-out with corner fillets at the center of the cylinder, as depicted in \cref{fig:KirchhoffLoveModelImmersed}, and treat the corresponding interior edges as immersed boundaries (or, trimmed patches). To repeat the simulation of the traveling shock-wave, we use the reference solution depicted in \cref{fig:KirchhoffLoveModel_b} as the manufactured solution for prescribing boundary conditions at these internal edges. Specifically, we extract the in-plane and out-of-plane displacement vectors and the normal rotation of the shell at the integration points around the cut-out. These three fields represent essential boundary conditions for the Kirchhoff-Love shell model. In the immersogeometric computations, we enforce all conditions with either a penalty method or with Nitsche's method. The exact Nitsche formulation that we adopt can be found in~\cite{Benzaken2021}. In this formulation, the three conditions all require a respective penalty parameter; $\beta_1$ for the out-of-plane displacement, $\beta_2$ for the normal rotation, and $\beta_3$ for the in-plane displacement. In \cite{Benzaken2021}, $\beta_1$ is referred to as $C^S_{\mathrm{pen},1}$, $\beta_2$ as $C^S_{\mathrm{pen},3}$, and $\beta_3$ as~$C^S_{\mathrm{pen},4}$.

\Cref{tab:KL} presents the critical time-step sizes for the various formulation for the example of \cref{fig:KirchhoffLoveModelImmersed}, with a background mesh consisting of $30\times 30$ quadratic B-spline elements, and with the penalty parameters chosen as $\beta_1 = 50$, $\beta_2 = 15$, and $\beta_3 = 2.5$.
The ghost penalty parameters for the contributions to the stiffness matrix (both for the in-plane displacement and the out-of-plane displacement field) and to the mass matrix are all set to 0.1. At these penalty values, the smallest eigenvalue of the system of equations arising from Nitsche's method is positive. As the critical time-step sizes in the table show, the use of ghost mass for the penalty method increases the critical time-step size by more than 15 times. With the Nitsche formulation, we retrieve a critical time-step size equal to that of the uncut case, which is the optimal value.

\begin{figure}[!b]
    \centering
    \includegraphics[width=0.55\linewidth,trim=10 12 10 10, clip]{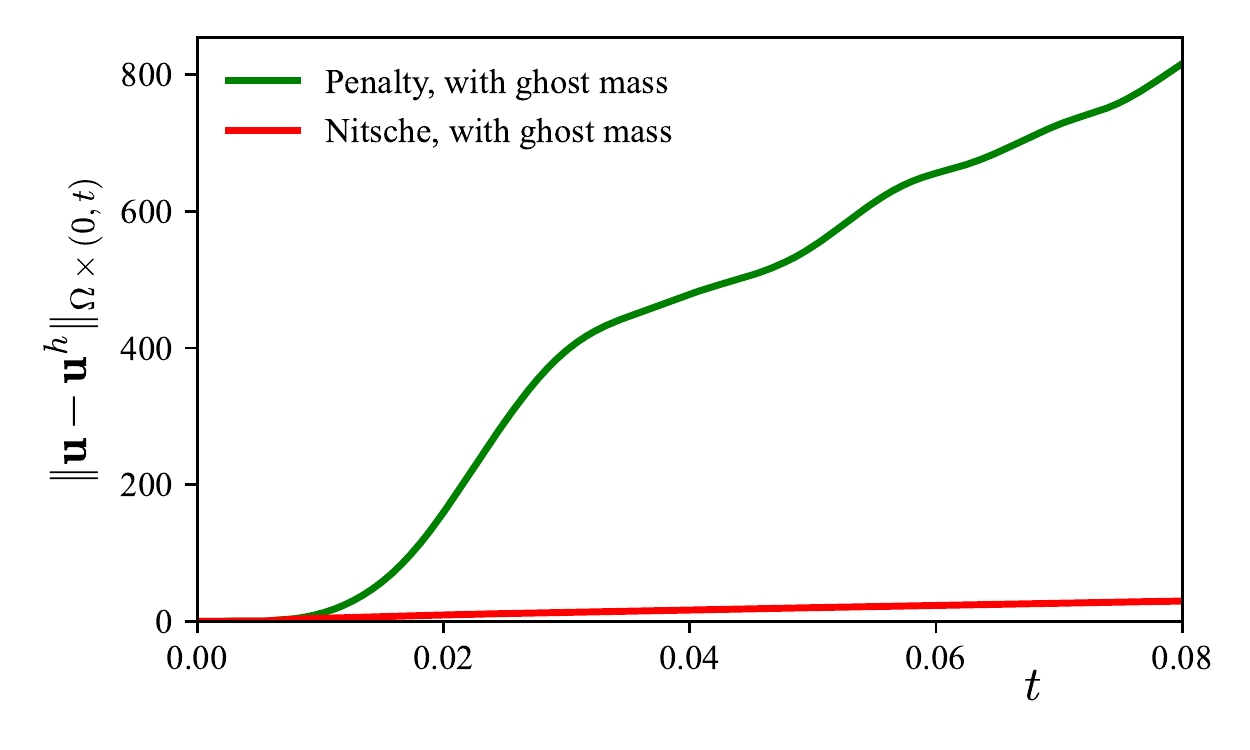}
    \caption{Evolution of the error in the Bochner norm for the penalty and Nitsche formulations of the Kirchhoff-Love shell.}
    \label{fig:KL_bochner}
\end{figure}
\begin{figure}
\vspace{-1cm}
    \centering
    \subfloat[Penalty formulation, with ghost mass.]{\includegraphics[width=0.5\linewidth,trim=230 230 140 200, clip]{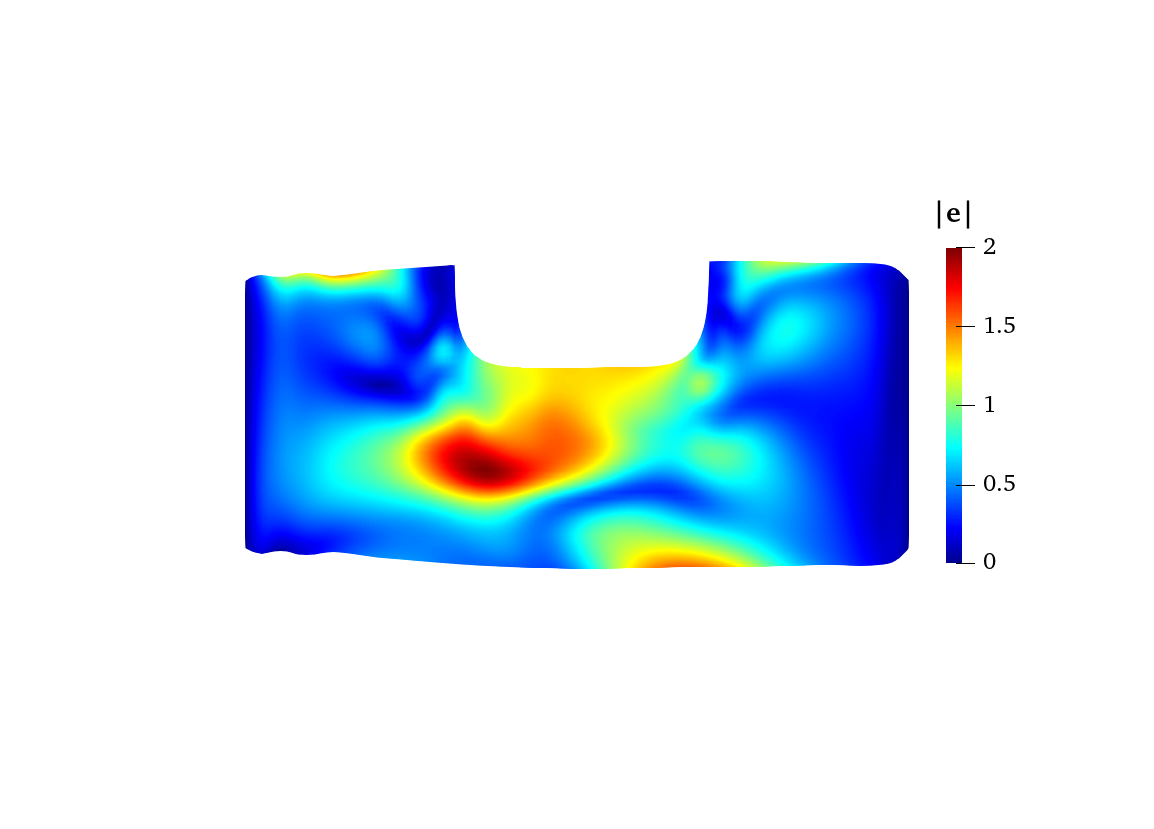}\label{fig:KL_error_a}}\hfill
    \subfloat[Nitsche formulation, with ghost mass.]{\includegraphics[width=0.5\linewidth,trim=230 230 140 200, clip]{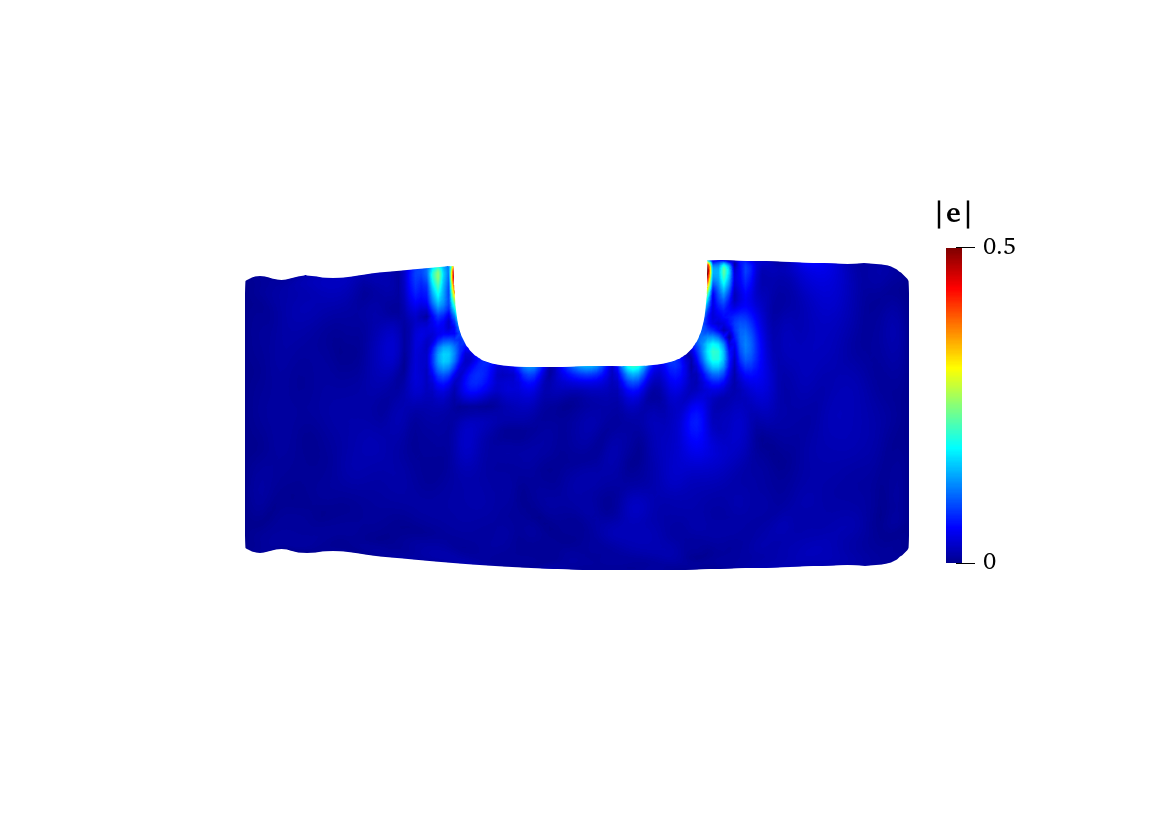}\label{fig:KL_error_b}}
    \caption{Magnitude of the displacement error in mm, at $t=0.045$.}
    \label{fig:KL_error}
\end{figure}

The traveling shock-wave computations are executed at the respective critical time-step size of each formulation. The evolution of the error over time is captured by the following time-dependent Bochner norm:
\begin{align}
    \norm{\vec{u}-\vec{u}^h}_{\Omega\times(0,t)} = \frac{1}{t} \int\limits_{0}^t \norm{\vec{u}(\tau,\cdot)-\vec{u}^h(\tau,\cdot)}_{L^2(\Omega)} \,\text{d}\tau \,.
\end{align}
The progression of this error is plotted in \cref{fig:KL_bochner} for both the penalty formulation and the Nitsche formulation. The adoption of the Nitsche formulation results in a reduction of the error by well over an order of magnitude. The manifestation of the error is depicted in \cref{fig:KL_error} for time $t=0.045$. The deviation from the expected physical response is significantly more pronounced for the penalty formulation than for the Nitsche formulation (also note the change in colorbar in this regard). In particular, the error distribution in \cref{fig:KL_error_a} indicates that penalty enforcement of the boundary conditions not only locally impacts the solution, but in fact disturbs the complete solution field. In contrast, the far-field impact of the Nitsche-based weak boundary condition enforcement is almost negligible, owing to its variational consistency.

\begin{table}[!t]
\caption{Critical time-step sizes for the different Kirchhoff-Love shell simulations.}
\begin{tabular}{|l|c|c|c|c|}\hline& & & &  \\[-0.4cm]
Method                         & \hspace{0.5cm}Uncut\hspace{0.5cm} & \hspace{0.2cm}\small\begin{tabular}[c]{@{}c@{}}Penalty \\ formulation \end{tabular}\hspace{0.2cm}    & \hspace{0cm}\small\begin{tabular}[c]{@{}c@{}}Penalty formulation, \\ with ghost mass \end{tabular}\hspace{0cm} & \small \hspace{0cm}\begin{tabular}[c]{@{}c@{}}Nitsche formulation, \\ with ghost mass \end{tabular}\hspace{0cm} \\[0.2cm] \hline\\[-0.45cm] \hline  & & & &  \\[-0.35cm]
$\Delta t_\text{crit}$ {[}ms{]}& 0.475 & 0.0311 & 0.393      & 0.475    \\[-0.35cm]        & & & & \\
\hline
\end{tabular}
\label{tab:KL}
\end{table}

\begin{figure}
    \centering
    \subfloat[Penalty formulation, with ghost mass.]{\includegraphics[width=0.49\linewidth,trim=10 10 10 10, clip]{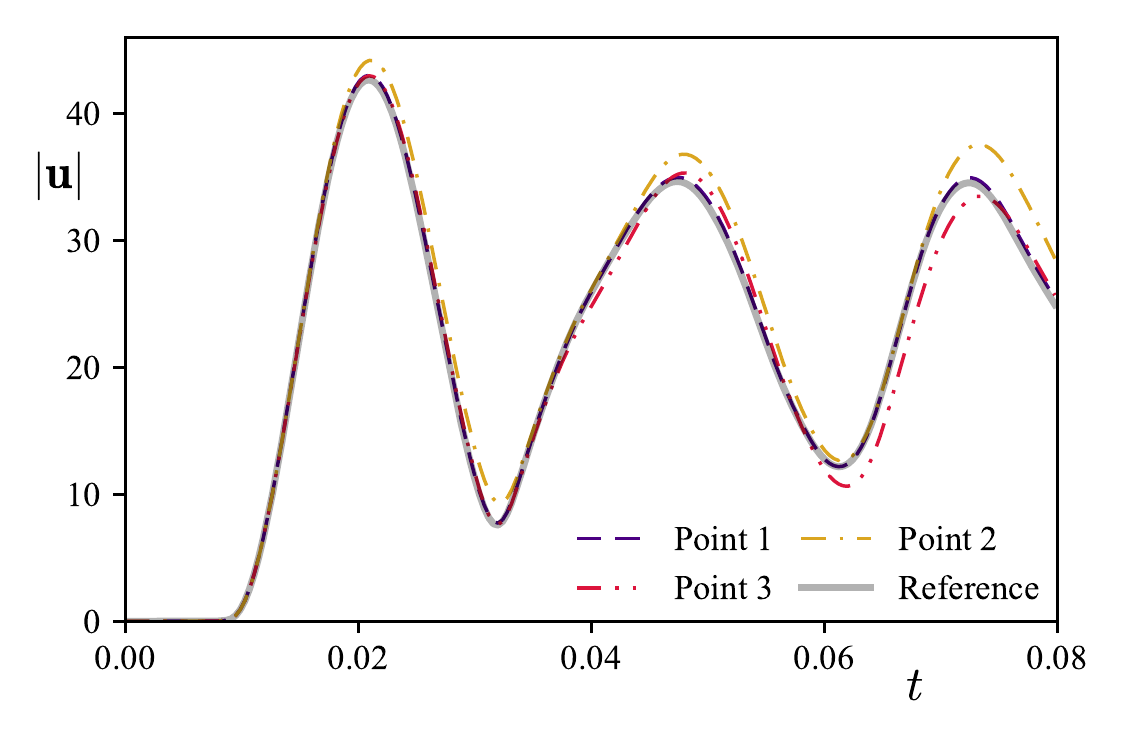}\label{fig:KL_pointvals_a}}\hfill
    \subfloat[Nitsche formulation, with ghost mass.]{\includegraphics[width=0.49\linewidth,trim=10 10 10 10, clip]{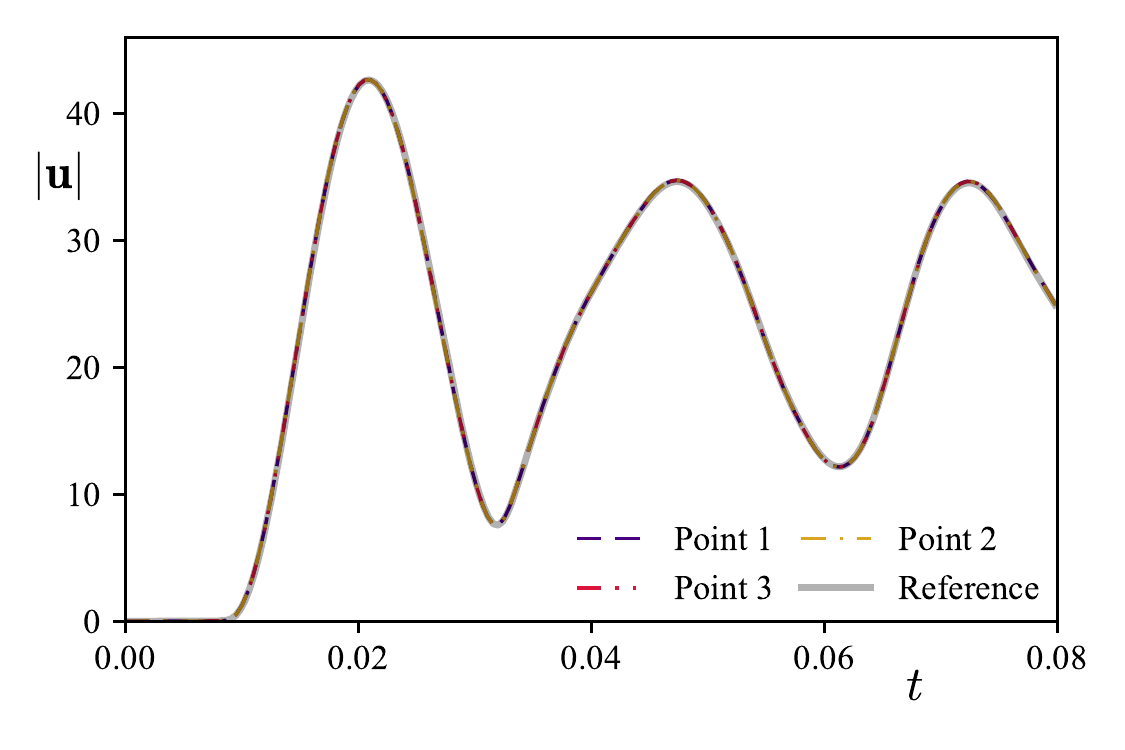}\label{fig:KL_pointvals_b}}
    \caption{Displacement magnitude over time at the three point locations indicated in \cref{fig:KirchhoffLoveModelImmersed}.} \label{fig:KL_pointvals}
\end{figure}

It can also be inferred from the error distribution that the axisymmetric nature of the true solution, as depicted in \cref{fig:KirchhoffLoveModel_b}, is lost in both finite element approximations. This loss is quantitatively assessed in \Cref{fig:KL_pointvals}, which plots the displacement magnitudes at the three locations along the span of the cylinder indicated in \cref{fig:KirchhoffLoveModelImmersed}. The results for the penalty method are shown in \cref{fig:KL_pointvals_a}, while those for the Nitsche formulation are shown in \cref{fig:KL_pointvals_a}. For the penalty method, the displacement magnitudes at the three points deviate increasingly from the reference solution, and also also diverge with respect to each other. This erroneous behavior is significantly suppressed when the Nitsche formulation is adopted, as evidenced by the nearly overlapping displacement curves for all points.

\section{Conclusion and outlook}
\label{sec:Conclusion}

Immersogeometric explicit analysis offers a powerful framework for streamlining the design-to-analysis pipeline for crash-test type simulations. To ensure a robust and reliable pipeline, it is imperative that the critical time-step size of the explicit time-stepping scheme is not affected by the size of the cut (or ``trimmed'') elements. In this article, we have studied the dependency of the critical time-step size on the cut-size for different types of boundary conditions and different methods of enforcement. The formulations that we investigated include a pure Neumann problem, penalty enforcement of Dirichlet constraints, and Nitsche's method for enforcement of Dirichlet constraints. To ensure a positive definite stiffness matrix when Nitsche's method is used, we considered a formulation where the Nitsche penalty parameter is determined based on an element-local eigenvalue problem, leading to a cut-size dependent penalty parameter, and a formulation with an additional ghost-penalty term, which then permits a cut-size independent penalty parameter. For each formulation, we considered a consistent mass matrix and a row-sum lumped mass matrix, both with and without additional ghost-penalty based mass scaling (``ghost mass'').

We have found that a formulation with a consistent mass matrix, without any form of mass scaling, always suffers from adverse scaling of the critical time-step size with element cut-size. As was first observed in \cite{LeidingerPhD}, this problematic scaling is mitigated by row-sum mass lumping, but only when the polynomial order of the maximum regularity splines is sufficiently high. Our analysis confirms this observation: for a second-order and a fourth-order equation we show that mass lumping alone is insufficient when the order of the basis functions is lower than the order of the spatial differential operator. To enable `lower-order' discretization (i.e., up to cubics for shell-type equations), an additional mass scaling is required. The cut-size dependency vanishes for our proposed addition of ghost mass.

Our analysis also shows that penalty enforcement of Dirichlet conditions suffers from the drawback that the critical time-step size scales with the penalty parameter after a threshold value, confirming the observations in \cite{Leidinger2019,LeidingerPhD}. While this is also the case for Nitsche's method for enforcement of Dirichlet conditions, the penalty factors required to achieve satisfactory results are generally lower. Out of the different Nitsche formulations that we considered, the only formulation for which the critical time-step size is independent of the cut-element size is one with ghost stiffness and ghost mass.

One favorable property of the addition of ghost mass to the problematically cut elements, is that it does not suffer from the severe variational inconsistency issues that plague more conventional mass-scaling methods. 
With numerical experiments, we have demonstrated the efficacy of ghost mass, as it can be utilized without incurring negative impacts on solution quality, despite a potentially substantial increase in the critical time-step size, sometimes by orders of magnitude. For a linear wave equation and a linear Kirchhoff-Love shell model, we have shown that the enforcement of Dirichlet conditions with a variationally consistent Nitsche method, as opposed to a penalty method, may lead to reductions of the solution error by orders of magnitude at the same critical time-step size.

An open research question pertains to the efficient (approximate) inversion of the mass matrix when ghost mass is added to the formulation. Suggestions for developing such a technique have been made in this work.

\enlargethispage{0.5cm}

\vspace{0.4cm}

\textbf{Acknowledgements:} S.K.F.\ Stoter gratefully acknowledges financial support through the Industrial Partnership Program {\it Fundamental Fluid Dynamics Challenges in 
Inkjet Printing\/} ({\it FIP\/}), a joint research program of Canon Production Printing, Eindhoven University of Technology,
University of Twente, and the Netherlands Organization for Scientific Research (NWO). C.V. Verhoosel and S.C. Divi acknowledge the partial support of the European Union’s Horizon 2020 research and innovation programme under Grant Agreement No 101017578 (SIMCor). M. Larson and E.H. van Brummelen gratefully acknowledge the insightful discussions at the special session organized by Prof. Trond Kvamsdal at IGA~2022 in Banff. All simulations have been performed using the open source software package Nutils~\cite{nutils}.

\bibliographystyle{ieeetr}
\bibliography{MyBib}

\end{document}